\documentclass{amsart}

\usepackage{amsmath,amssymb,amsfonts,amscd}
\usepackage[all,cmtip]{xy}
\usepackage{enumerate}
\usepackage{ amssymb, latexsym, amsmath}

\usepackage{lingmacros}
\usepackage{tree-dvips}
\usepackage{amssymb,amsfonts,amsthm,amsmath}

\usepackage{mathrsfs}
\usepackage[english]{babel}
\usepackage[all,cmtip]{xy}
\usepackage{tikz}
\usepackage{tikz-cd}
\usetikzlibrary{cd,quotes}
\usetikzlibrary{calc,babel}
\usepackage{mathtools}

\usepackage{tensor}
\usepackage{csquotes}
\usepackage[breaklinks]{hyperref}
\usetikzlibrary{arrows,chains,matrix,positioning,scopes}

\usepackage{fullpage}
\usepackage{url}
\usepackage{hyperref}
\usepackage{ marvosym }
\usepackage{amssymb}

\numberwithin{equation}{section}

\theoremstyle{plain}
\newtheorem{theorem}{Theorem}[section]

\newtheorem{proposition}[theorem]{Proposition}
\newtheorem{lemma}[theorem]{Lemma}

\newtheorem{corollary}[theorem]{Corollary}

\theoremstyle{definition}
\newtheorem{definition}[theorem]{Definition}

\newtheorem{example}[theorem]{Example}

\newcommand{\CP}{\mathbb{CP}}
\newcommand{\R}{\mathbb{R}}

\newcommand{\C}{\mathbb{C}}
\newcommand{\Z}{\mathbb{Z}}
\newcommand{\RP}{\mathbb{RP}}
\newcommand{\w}{\omega}
\newcommand{\T}{\mathcal{T}}
\newcommand{\F}{\mathcal{F}}
\newcommand{\E}{\mathcal{E}}

\theoremstyle{remark}
\newtheorem{remark}[theorem]{Remark}

\makeatletter

\newcommand{\Rmnum}[1]{\expandafter\@slowromancap\romannumeral #1@}
\makeatother

\newcommand{\Q}{\mathbb{Q}}

\newcommand{\M}{\overline{\mathcal{M}}}

\begin{document}
\title{Fukaya Algebra over $\Z$} 

\author{Mohamad Rabah}
\email{mohamad.rabah@stonybrook.edu, mohamad@amss.ac.cn}

\begin{abstract}
   Given a closed, connected, relatively-spin Lagrangian submanifold in a closed symplectic manifold, we associate to it a curved, gapped, filtered, $A_{n, K}$-algebra over the Novikov ring with integer coefficients, which we extend to an $A_\infty$-algebra. To illustrate our framework, we give a proof of the Quantum Lefschetz Hyperplane Theorem in the K$\ddot a$hler case, and associate virtual fundamental classes to the moduli spaces used in local Gromov-Witten theory, in the symplectic case. 
\end{abstract}
\maketitle

\section{Introduction}
In their book \cite{Fukaya2009}, Fukaya-Oh-Ohta-Ono constructed an $A_\infty$-algebra structure on the singular cohomology of any closed, connected, relatively-spin Lagrangian submanifold $L$ of a closed symplectic manifold $(X, \w)$, over the Novikov ring with rational coefficients. To be more precise, let $\Lambda^\Q_{0, nov}$ be the Novikov ring over the rationals. Namely,
\begin{equation*}
    \Lambda^\Q_{0, nov}:=\{\Sigma_ia_iT^{\lambda_i}e^{n_i}:a_i\in\Q, n_i\in\Z, 0\leq\lambda_1\leq\lambda_2\leq\dots, \lim_i\lambda_i=\infty\}.
\end{equation*}
Where $T, e$ are formal variables of degrees $0,2$ respectively and consider $H^*(L;\Lambda^\Q_{0, nov})\equiv H^*(L;\Q)\otimes_\Q \Lambda^\Q_{0, nov}$.
\begin{theorem}[\cite{Fukaya2009}]
    To each closed, connected, relatively-spin Lagrangian submanifold $L$ of $X$ we can associate an $A_\infty$-algebra structure $\{m_k\}_{k\geq 0}$ on $H^*(L;\Lambda^\Q_{0, nov})$, which is well-defined up to isomorphism.
\end{theorem}
The restriction to rational coefficients in the construction of
Fukaya-Oh-Ohta-Ono is not merely algebraic. Their perturbation theory uses multi-sections on Kuranishi spaces with nontrivial finite isotropy, whose zero sets naturally carry rational weights. A similar restriction appears in Pardon's virtual-cycle construction as in \cite{pardon2016algebraic}. Indeed, the relevant Poincaré-duality statement for orbifolds, viewed through their underlying $\textit{coarse}$ spaces, generally holds only with rational coefficients, since finite isotropy obstructs integral Poincaré duality. Consequently, constructing an $A_\infty$-algebra over the integers requires an integral virtual-cycle theory compatible with the boundary and corner stratifications of the moduli spaces of stable disks. The aim of this paper is to prove an $\textit{integral}$ version of the above result. Our construction applies to arbitrary closed, connected, relatively-spin Lagrangian submanifolds of closed symplectic manifolds, without imposing exactness or monotonicity assumptions on either the Lagrangian or the ambient symplectic manifold. We fix notation and denote by $\Lambda:=\Lambda^\Z_{0, nov}$ the Novikov ring over the integers. That is, 
\begin{equation*}
    \Lambda:=\{\sum_ia_iT^{\lambda_i}e^{n_i}:a_i,n_i\in\Z, 0\leq\lambda_1\leq\lambda_2\leq\dots, \lim_i\lambda_i=\infty\}.
\end{equation*}
Let $(f,g)$ a Morse-Smale pair on $L$ and denote by $\operatorname{Crit}(f)$ the finite set of critical points of $f$ and $CM(f;\Z):=\Z\langle\operatorname{Crit}(f)\rangle$ the Morse complex on $L$ generated by critical points of $f$. We set $CM(f; \Lambda):=CM(f; \Z)\widehat{\otimes}_{\Z}\Lambda$ where $\widehat{\otimes}$ is the completion of the tensor product with respect to the $T$-adic topology.
Our main result is the following.
\begin{theorem}\label{MainResult}
     The Morse complex $CM(f; \Lambda)$ carries the structure of a curved, gapped, filtered $A_{n, K}$-algebra $\{m_k\}_{k\geq0}$. Moreover, such algebra can be extended to a curved, gapped, filtered $A_\infty$-algebra.
\end{theorem}
The proof that we give for our main result requires two main ingredients which we outline here.

\subsection{Transversality of Orbifold Sections}
In general, one cannot generically perturb a section of an orbibundle and make it transverse to the zero-section. This was resolved by Fukaya-Ono in \cite{fukaya1996arnold} by using the notion of $\textit{multi-valued}$ sections. In return, the use of multi-valued sections cannot give us integer-valued cycles or classes. In \cite{fukaya2001floer}, Fukaya-Ono showed us how one can construct an $\textit{integral Euler class}$ from a single-valued orbibundle section, provided that the orbibundle is $\textit{normally complex}$. To this end and following the discussion in \cite{bai2022arnold}, \cite{bai2022integral} we introduce the following notations. Let $\E\rightarrow\mathcal{U}$ be an orbibundle and suppose that $(\Gamma, U, E)$ is an orbibundle chart. For any subgroup $G\subset\Gamma$ we denote by $U^G\subset U$ the $G$-fixed locus of the induced $G$-action and we denote by $NU^G\rightarrow U^G$ the normal bundle of $U^G\subset U$. Moreover, over $U^G$ we can decompose $E|_{U^G}=\dot{E}^G\oplus\check{E}^G$ as a direct sum of trivial and non-trivial $G$-representations, respectively. With this in mind and following \cite{bai2022integral}, we define the following notions.
\begin{definition}[Normal Complex Structure on $\mathcal U$]
    A normal complex structure on $\mathcal U$ is the collection of the following data:
    \begin{enumerate}
        \item For each effective orbifold chart $(\Gamma, U, \psi)$ of $\mathcal U$ and for any subgroup $G\subset\Gamma$, we have a $G$-equivariant complex structure $I_G:NU^G\rightarrow NU^G$.
        \item For $x\in U$ with stabilizer group $\Gamma_x$ and for any pair of subgroups $G\subseteq H\subseteq\Gamma_x$, we identify the fibre $N_xU^H\cong (T_xU^G\cap N_xU^H)\oplus N_xU^G$, where $T_xU^G\cap N_xU^H$ is the trivial $G$-representation $N_xU^H$, while $N_xU^G$ is the non-trivial $G$-representation of $N_xU^H$. We require $I_H|_{N_xU^G}\equiv I_G(x)$.
        \item For a given chart embedding $(V, \Pi, \phi)\hookrightarrow(U, \Gamma, \psi)$ given by an injective group morphism $\Pi\hookrightarrow\Gamma$ and an equivariant embedding $\iota:V\hookrightarrow U$; and for each $y\in V$ such that $\iota(y)=x\in U$ together with an identification of stabilizer groups $\Pi_y\cong\Gamma_x$, such identification gives us an equivariant linear isomorphism $T_y V\cong T_xU$, inducing a linear isomorphism $N_yV^{\Pi_y}\cong N_xU^{\Gamma_x}$. We require that such induced linear isomorphism $N_yV^{\Pi_y}\cong N_xU^{\Gamma_x}$ to be $(I_{\Pi_y}(y), I_{\Gamma_x}(x))$-linear.
    \end{enumerate}
\end{definition}
\begin{definition}[Normal Complex on $\E$] \label{normallyComplexOrbibundle}
A normal complex structure on $\E$ is the collection of the following data:
\begin{enumerate}
    \item For any given orbibundle chart $(\Gamma, U,E)$ of $\E\rightarrow\mathcal U$ and for any subgroup $G\subset\Gamma$, we have a $G$-equivariant complex structure $J_G:\check{E}^G\rightarrow\check{E}^G$ on the bundle $\check{E}^G\rightarrow U^G$.
    \item For each bundle chart $(\Gamma, U, E)$, $x\in U$ with stabilizer $\Gamma_x\subseteq\Gamma$, and each pair of subgroups $H\subseteq G\subseteq\Gamma_x$, we require $J_G|_{\check{E}^H}(x)\equiv J_H(x)$.
    \item For any given orbibundle chart embedding $(\Pi, V, F)\hookrightarrow(\Gamma, U, E)$ given by an injective group morphism $\Pi\hookrightarrow\Gamma$ and an equivariant bundle embedding $F\hookrightarrow E$ covering a base embedding $\iota:V\hookrightarrow U$; for any $y\in V$ such that $\iota(y)=x\in U$ given with the injective linear map on the fibres $F_y\hookrightarrow E_x$ induced from the chart embedding, we require that such injective linear map induces a $(J_{\Pi_y}(y), J_{\Gamma_x}(x))$-linear isomorphism $\check{F}^{\Pi_y}_y\cong E_x^{\Gamma_x}$.
\end{enumerate}
\end{definition} 
In the setting of a normally complex orbibundle $\E\rightarrow\mathcal{U}$, Fukaya-Ono defined the ($C^0$-dense) sections that allow us to extract integer-valued Euler classes. In order to define such sections and also following \cite{bai2022integral}, \cite{bai2022arnold} we introduce the following notations. Suppose that $B$ is a smooth manifold and $G$ is a finite group acting trivially on $B$. Let $\pi_V:V\rightarrow B$ and $\pi_W: W\rightarrow B$ be two smooth complex $G$-equivariant vectorbundles and fix a positive integer $d$. We have an induced vectorbundle $\operatorname{Poly}^G_d(V, W)\rightarrow B$ whose fibre over $b\in B$ is the vectorspace $\operatorname{Poly}^G_d(V_b, W_b)$ of all $G$-equivariant complex polynomial maps $V_b\rightarrow W_b$ of degree at most $d$. Suppose that $V$ is also equipped with a bundle metric and denote by $V_\epsilon$ a $G$-invariant disk subbundle of $V$. Consider $s:V_\epsilon\rightarrow \pi_V^*W$ a smooth $G$-equivariant section.
\begin{definition}[Normally Complex Polynomial Section]
    $s$ is said to be a normally complex polynomial section of degree at most $d$ if the restriction of $s$ to each fibre $V_b\cap V_\epsilon$ is an element of $\operatorname{Poly}^{G}_d(V_b, W_b)$.
\end{definition}

\subsection{Global Kuranishi Charts}
In \cite{abouzaid2021complex}, \cite{abouzaid2024gromov}, Abouzaid-McLean-Smith, showed us how to represent the moduli space of closed, stable pseudo-holomorphic curves by a single Kuranishi chart. To this end let $X$ be any metrizable space.

\begin{definition}[Global Kuranishi Chart]
    A Global Kuranishi chart of $X$ is a tuple $(G, \T, \E, s)$ where:
    \begin{enumerate}
        \item $G$ is a compact Lie group.
        \item $\T$ and $\E$ are smooth $G$-manifolds, possibly with boundaries and corners, where the $G$-action is locally-linear as in Definition \ref{locallyLinearAction}, and with finite stabilizers.
        \item $\E\rightarrow\T$ is a $G$-equivariant vectorbundle.
        \item $s:\T\rightarrow\E$ is a $G$-equivariant continuous section.
    \end{enumerate}
    such that, $s^{-1}(0)/G\cong X$ as topological spaces.
\end{definition}
We will call $\T$ the $\textbf{thickening}$, $\E$ the $\textbf{obstruction bundle}$, $G$ the $\textbf{symmetry group}$ and $s$ the $\textbf{Kuranishi map}$ or section.\\
For instance, any smooth manifold or more generally, any smooth effective orbifold can be given a Global Kuranishi chart as follows.
\begin{example}
\begin{enumerate}
    \item Suppose that $X$ is a smooth manifold. Then the tuple $(\{\mathbf{1}\}, X, \underline{0}, 0)$, where $\{\mathbf{1}\}$ is the trivial group and $\underline{0}\equiv X\times\{0\}$ is the zero-bundle, is a Global Kuranishi chart of $X$.
    \item  Suppose that $X$ is a smooth effective $n$-orbifold. By the $\textit{Orbifold Frame Bundle Theorem}$ as per \cite{adem2007orbifolds}, $X\cong Fr(X)/O(n)$ where $Fr(X)$ is the orthonormal frame bundle of $X$ and the $O(n)$-action is given by change of orthonormal basis, and hence locally-linear. Thus, $(O(n), Fr(X), \underline{0}, 0)$ is a Global Kuranishi chart representation of $X$, where $\underline{0}\equiv Fr(X)\times\{0\}$.\\
    More generally, if $X$ is a smooth $n$-orbifold, then by \cite{pardon2022enough}, there exists a smooth manifold $P$ and a positive integer $n$ such that, $X\cong P/U(n)$ as topological spaces. Moreover, $U(n)$ acts smoothly on $P$ with finite stabilizers. Then $(U(n), P, \underline{0}, 0)$ is a Global Kuranishi chart of $X$.
\end{enumerate}   
\end{example}

\begin{definition}[Oriented Global Kuranishi Chart]
    A Global Kuranishi chart $(G, \T, \E, s)$ is said to be oriented if both $\T$ and $\E$ are oriented, and the action of $G$ preserves such orientation.
\end{definition}
Notice that, for a given metrizable space $X$, we may have more than one Global Kuranishi chart representing $X$. Such Global Kuranishi charts are said to be $\textit{equivalent}$ as per \cite{abouzaid2021complex}, if they are related by the following operations. Namely, let $(G, \T, \E, s)$ be a Global Kuranishi chart of $X$;
\begin{enumerate}
    \item (Germ Equivalence): If $s^{-1}(0)\subseteq U\subseteq\T$ is a $G$-invariant open set of $s^{-1}(0)$ then, $(G, U, \E|_U, s|_U)$ is a Global Kuranishi chart of $X$ equivalent to $(G, \T, \E, s)$.
    \item (Group Enlargement): If $G'$ is a compact Lie group and $\pi:P\rightarrow\T$ is a $G$-equivariant principal $G'$-bundle then, $(G'\times G, P, \pi^*\E, \pi^*s)$ is a Global Kuranishi chart of $X$ equivalent to $(G, \T, \E, s)$.
    \item (Stabilization): If $\pi: V\rightarrow\T$ is a $G$-equivariant vectorbundle then, $(G, V, \pi^*\E\oplus\pi^*V, \pi^*s\oplus\tau_V)$ is a Global Kuranishi chart of $X$ equivalent to $(G, \T, \E, s)$, where $\tau_V:V\rightarrow \pi^*V$ is the tautological section.
\end{enumerate}
Observe that, if $(G, \T, \E, s)$ is a Global Kuranishi chart of $X$, from the assumption that $G$ acts by finite stabilizers, it follows, possibly after stabilizing by a free $G$-representation, that $\E/G\rightarrow\T/G$ is a smooth effective orbifold. With this in mind, we give the following definition.
\begin{definition}[Normally Complex Global Kuranishi Chart] A Global Kuranishi chart $(G, \T, \E, s)$ of $X$ is said to be normally complex, if $\E/G\rightarrow\T/G$ is a normally complex orbibundle as in Definition \ref{normallyComplexOrbibundle}.
\end{definition}

\subsection{Organization}
In $\textbf{section 2}$, we outline the key steps in constructing Global Kuranishi charts for the moduli space of genus-zero pseudo-holomorphic curves. This framework is then used to reprove well-known results in $\textit{Gromov-Witten Theory}$. We refer the reader to Theorem \ref{TwistedGWNegative} and Theorem \ref{QuanLefsch} for precise statements. \\

In $\textbf{section 3}$, we prove that the moduli space of genus-zero bordered pseudo-holomorphic curves with Lagrangian boundary conditions admits a Global Kuranishi chart representation. We refer the reader to Theorem \ref{ChartsDisks} for the precise statement. Additionally, we construct these Global Kuranishi charts in a way compatible with the boundary stratification of the moduli space of genus-zero bordered pseudo-holomorphic curves. This extends the work of Bai-Xu \cite{bai2022arnold} or Rezchikov \cite{rezchikov2022integralarnoldconjecture}, to the case where the boundary of the moduli space is a fibre product of lower dimensional moduli spaces. We refer the reader to Theorem \ref{ChartsforDisksCompatible} for the precise statement. In fact, as a by-product of Theorem \ref{ChartsforDisksCompatible}, one can recover Theorem A of \cite{Fukaya2009}.\\

In $\textbf{section 4}$, we show that the moduli space of $\textit{Floer-Morse Trajectories}$, as defined in Definition \ref{ModuliFloerMorse}, admits a normally complex Global Kuranishi chart representation, which is compatible with its boundary stratification. We refer the reader to Theorem \ref{ThickeningofTheMainComponent} and Theorem \ref{NormalComChartMorse} for the precise statements. We also give the necessary constructions that provide us with an $A_\infty$-algebra structure on Morse cochains, as in Corollary \ref{MorseA-infinity}. Such result is used in the proof of our main result, namely in the proof of Corollary \ref{Ainfty}.\\

We end this paper with two appendices compiling all the necessary definitions and results that are used freely in our arguments. $\textbf{Appendix A}$ contains all the necessary notions and statements used to give a smooth structure on the $\textit{topological Thickenings}$. While $\textbf{Appendix B}$ contains all the necessary notions and results that are employed in our proofs to get integral Euler cycles in our setting.

\section*{Acknowledgements}I am grateful to my Ph.D. advisor, Kenji Fukaya, for a great number of enlightening meetings throughout the years and his generosity in sharing his ideas. Equally, I am also thankful to Mark McLean for his generosity in his time and ideas. Also I would like to thank Jiaji Cai, Shuhao Li, Chris Woodward, Hang Yuan, Yao Xiao and Guangbo Xu for many fruitful discussions. Special thanks go to my girlfriend, Tracy Zhang Ke, for her warmth and unwavering motivation. We are deeply grateful to the referees for their valuable comments, which enabled us to strengthen our main result.
\\
This work is dedicated to the memory of Moussa Daod Sleiman Bazzi.

\section{An Overview on Global Kuranishi Charts for the Moduli Space of Pseudo-holomorphic Spheres} 

\subsection{Construction}
In \cite{abouzaid2021complex}, the authors showed that the moduli space of pseudo-holomorphic spheres in a closed or a geometrically bounded as in \cite{gromov1985pseudo}, symplectic manifold admits a Global Kuranishi chart. Throughout this section, we assume $(X, \w)$ is a closed symplectic manifold and fix an $\w$-compatible almost complex structure $J$ on $X$. Let $\beta\in H_2(X; \Z)$ and $n$ be a non-negative integer. Denote by $\M_{0, n}(X; J, \beta)$ the moduli space of stable genus-zero $J$-holomorphic curves with $n$-marked points representing the class $\beta$.
\begin{theorem}\cite{abouzaid2021complex}\label{AMSI}
    There exists a Global Kuranishi chart representation of $\M_{0, n}(X; J, \beta)$.
\end{theorem}
The purpose of this section is to give an outline of the key ideas that goes behind proving Theorem \ref{AMSI} following \cite{siebert1999gromov}, \cite{abouzaid2021complex}, \cite{Hirschi:2022twt} and \cite{bai2022arnold}.\\
Fix an integral lift $\Omega$ of $\w$ which is compatible with the almost complex structure $J$ on $X$. Namely, let $\Omega\in\Omega^2(X)$ be a symplectic form, compatible with $J$ and $[\Omega]\in H^2(X; \Z)$. Moreover, we fix a Hermitian line bundle $(L, \nabla)\rightarrow X$ of curvature form $F_\nabla=-2\pi i\Omega$. For $u:\Sigma\rightarrow X$ a $W^{k, p}$-map from a nodal genus-zero Riemann surface, representing the class $\beta$, that is $u_*[\Sigma]=\beta$, we denote by $L_u:=u^*L$.
\begin{lemma}\label{ample}
    $L_u\rightarrow\Sigma$ is a holomorphic vectorbundle, which is ample on each irreducible component of $\Sigma$.
\end{lemma}
Lemma \ref{ample}, together with the Serre duality and the Kodaira vanishing theorem, implies that $H^1(\Sigma; L_u)=\{0\}$ and hence by Riemann-Roch, we have $\operatorname{dim}_\C H^0(\Sigma; L_u)=d+1$ where $d:=\Omega(\beta)$.
\begin{definition}
    A framed curve is a tuple $(u, \Sigma, F)$ where:
    \begin{enumerate}
        \item $\Sigma$ is a nodal genus-zero Riemann surface.
        \item $u:\Sigma\rightarrow X$ is a $W^{k, p}$-map such that $kp$ is larger than the dimension of $X$.
        \item $F=\{f_0, \dots, f_d\}$ is a $\C$-basis of $H^0(\Sigma; L_u)$ such that the matrix $(\int_\Sigma\langle f_i, f_j\rangle u^*\Omega)_{i, j}$ has positive eigenvalues, where $\langle., .\rangle$ is the induced Hermitian pairing on $H^0(\Sigma; L_u)$. 
    \end{enumerate}
\end{definition}
Let $(u, \Sigma, F)$ be a framed curve and write $F=\{f_0, \dots, f_d\}$. We have an induced holomorphic map
\begin{align*}
    \phi_F: \Sigma \ & \to \mathbb {CP}^d \\
    z \ & \mapsto [f_i(z)]
\end{align*}
Let $k>>1$ be an integer, and consider the following vectorspace
\begin{equation}\label{AMSIOb}
    H^0(\Sigma; \overline{\operatorname{Hom}}(\phi_F^*T\CP^d, u^*TX)\otimes_\C \phi_F^*\mathcal{O}(k))\otimes_\C \overline{H^0(\Sigma; \phi_F^*\mathcal{O}(k))}.
\end{equation}
\begin{definition}[Topological Thickening]
    We define the (topological) thickening of $\M_{0, n}(X; J, \beta)$ to be the set $\T(\beta):=\{(u, \Sigma, F, \eta)\}$ where:
    \begin{enumerate}
        \item $(u, \Sigma, F)$ is a framed curve.
        \item $\eta$ is an element of \ref{AMSIOb}.
    \end{enumerate}
    such that $u$ represents the class $\beta$ and satisfies the following equation:
    \begin{equation*}
        \Bar{\partial}_J u +\langle\eta\rangle\circ d\phi_F=0
    \end{equation*}
    where $\langle . \rangle$ is the natural Hermitian pairing
    \begin{equation}\label{Pairing}
        H^0(\Sigma; \overline{\operatorname{Hom}}( \phi_F^*T\CP^d, u^*TX)\otimes_\C \phi_F^*\mathcal{O}(k))\otimes_\C \overline{H^0(\Sigma; \phi_F^*\mathcal{O}(k))}\rightarrow \Gamma(\Sigma; \overline{\operatorname{Hom}}( \phi_F^*T\CP^d, u^*TX))
    \end{equation}
\end{definition}
\begin{proposition}\cite{abouzaid2021complex}, \cite{bai2022arnold} \label{AMSIfibre-wise}
For $k>>1$ large enough, $\T(\beta)$ is a topological $U(d+1)$-manifold. Moreover, there exists a $U(d+1)$-smooth manifold $\F$ together with a topological submersion $\pi:\T(\beta)\rightarrow \F$, as in Definition \ref{TopSubmersion}, which is a fibre-wise smooth $U(d+1)$-bundle, as in Definition \ref{Fibre-wise C1}.
\end{proposition}
Before we define the obstruction bundle, we introduce the following notation. Denote by $\mathcal{H}_d$ the space of $(d+1)\times (d+1)$-Hermitian matrices and by $\mathcal{H}_d^+\subset\mathcal{H}_d$ the subset of all Hermitian matrices of positive eigenvalues. Let $\exp:\mathcal{H}_d\rightarrow\mathcal{H}_d^+$ be the exponential map and $\exp^{-1}$ its inverse.
\begin{definition}[Obstruction Bundle]
    We define the obstruction bundle to be the $U(d+1)$-vectorbundle $\E\rightarrow\T(\beta)$ of fibre over $(u, \Sigma, F, \eta)$ given by
    \begin{equation*}
         (H^0(\Sigma; \overline{\operatorname{Hom}}( \phi_F^*T\CP^d, u^*TX)\otimes_\C \phi_F^*\mathcal{O}(k))\otimes_\C \overline{H^0(\Sigma; \phi_F^*\mathcal{O}(k))})\oplus\mathcal{H}_d
    \end{equation*}
\end{definition}
\begin{definition}[Kuranishi Map]
    We define the Kuranishi map $s:\T(\beta)\rightarrow\E$ by 
    \begin{equation*}
        s(u, \Sigma, F, \eta):=(\eta, \exp^{-1}(\int_\Sigma\langle f_i, f_j\rangle u^*\Omega)_{i, j})).
    \end{equation*}
    
\end{definition}
Now using the result of Proposition \ref{AMSIfibre-wise} and Lashof's Theorem \ref{lashof}, it follows using Proposition $5.35$ of \cite{abouzaid2021complex}, possibly after stabilization, $(U(d+1), \T(\beta), \E, s)$ is a (stably complex) Global Kuranishi chart of $\M_{0, n}(X; J, \beta)$.

\subsection{Quantum Lefschetz Principle}
In this section, we are interested in a comparison result between the virtual fundamental classes of the moduli space of stable genus-zero pseudo-holomorphic curves arising from the following geometric setting.\\
Let $(X, \omega)$ be a closed symplectic manifold and $J$ an $\omega$-compatible almost complex strucutre on $X$. Let $\pi:(H, \nabla, h)\rightarrow X$ be a Hermitian line bundle over $X$. We fix on the total space of $H$ the product almost complex structure $J_H:=\pi^*J\oplus i$ where $i$ is the fibre-wise multiplication by $\sqrt{-1}$. For such almost complex structure $J_H$ on the total space of $H$, the fibres of $H$ and the zero-section $0_H\cong X$ are almost complex submanifolds of the total space of $H$. Moreover, $\pi$ is a $(J_H, J)$-holomorphic map. That is, $J\circ d\pi=d\pi\circ J_H$.
We note that, in general the moduli space of $J_H$-holomorphic maps into $L$ representing a fixed homology class, is not compact due to the non-compactness of the total space of $H$. On the other hand, it is a first countable, Hausdorff space when given the Gromov-convergence topology \cite{mcduff2012j}.
\begin{example}\label{example}
    Consider the trivial $\C$-line bundle over $\CP^1$ and let $s$ be a non-vanishing holomorphic section. Then the sequence $\{s_n:=ns\}_{n\in\mathbb{N}}$ has no convergent subsequence. 
\end{example}
\subsubsection{$H$ is Concave}
Following \cite{lian1997mirror}, we recall the following definition.
\begin{definition}[Concave Vectorbundle]\label{concave} 
    Let $V\rightarrow\CP^n$ be a holomorphic $\mathbf{T}$-equivariant vectorbundle. $V$ is said to be concave if 
    \begin{enumerate}
        \item $e_\mathbf{T}(V)\in H_\mathbf{T}^*(\CP^n; \mathbb{Q})$ is invertible.
        \item $H^0(\Sigma; u^*V)=\{0\}$ for any $u:\Sigma\rightarrow\CP^n$ genus-zero holomorphic curve.
    \end{enumerate}
    where $\mathbf{T}=(S^1)^{\operatorname{rk}_\C V}$ is the torus. 
\end{definition}
As any holomorphic vectorbundle over $\CP^1$ splits uniquely as a direct sum of line bundles, we restrict our attention to the case of line bundles as we are only dealing with genus-zero curves. Moreover, we weaken the above definition and say $H\rightarrow X$ is a $\textbf{concave Hermitian Line bundle}$ over $X$ if only the second item of Definition \ref{concave} holds.
\begin{example}
\begin{enumerate}
    \item  Consider the $K_{\CP^n}\rightarrow\CP^n$ the canonical line bundle. Notice that $c_1(K_{\CP^n})=-(n+1)[\w_{FS}]$ where $\w_{FS}$ is the Fubini-Study symplectic form. More generally, if $(X, \w)$ is a rational, monotone symplectic manifold, then for $N>>1$ large enough $K^{\otimes N}_X\rightarrow X$ is concave.
    \item $H\rightarrow X$ is a Hermitian line bundle such that $c_1(H)=-k[\Omega]$ where $k$ is a positive integer and $[\Omega]\in H_2(X; \Z)$ is an integral lift of $\omega\in\Omega^2(X)$.
\end{enumerate}
\end{example}
\begin{lemma}\label{lemmaImage}
    The image of any non-constant $J_H$-holomorphic map in $H$ lies in $X$, where $X$ is viewed as the zero-section of $H$.
\end{lemma}
\begin{proof}
    Let $\Sigma\cong\bigvee_j\CP^1_j$ be a nodal genus-zero Riemann surface and consider a $J_H$-holomorphic map $u:\Sigma\rightarrow H$ to the total space of $H$ of positive degree. That is, $[u]\neq 0\in H_2(H; \Z)\cong H_2(X; \Z)$ where the isomorphism is given by $\pi_*$. Denote by $v:=\pi\circ u$ and note that $v$ is a $J$-holomorphic map since $\pi$ is $(J_H, J)$-holomorphic. Moreover, $[v]=\pi_*[u]$ and thus, $v$ is non-constant.
    Write $u=(v, s)$ where $s\in\Gamma(\Sigma, v^*H)$ is a smooth section and observe that $u$ is $J_H$-holomorphic if and only if $v$ is $J$-holomorphic and $s\in H^0(\Sigma, v^*H)$. Indeed, consider a non-constant irreducible component $v_j:\Sigma_j\cong\CP^1\rightarrow X$ and by abuse of notation, we still denote by $\nabla$ the pullback connection on $v_j^*H$. Then for $x\in\Sigma_j$ and for $w\in T_x\Sigma_j$, we have $d_xu(w)=(d_xv(w), \nabla_w s)$. Upon taking the $(0, 1)$-part of this equation and noting that any complex line bundle over $\Sigma$ is holomorphic with unique holomorphic structure; the result follows as $H^0(\Sigma; v^*H)=\{0\}$.
\end{proof}
\begin{remark}
    In the above proof, we only considered non-constant irreducible components. Indeed, if we have a ghost bubble, then its necessarily glued to a non-constant component as per Proposition $4.7.1$ of \cite{mcduff2012j} and hence its image is in zero-section. 
\end{remark}
Fix a non-zero class $\beta\in H_2(H; \Z)$ and denote by $\M_{0, n}(H; J_H, \beta)$ the moduli space of stable genus-zero $J_H$-holomorphic maps in the total space of $H$ with $n$-marked points, representing the class $\beta$. Similarly, we abuse notation and denote by $\M_{0, n}(X; J, \beta)$ the moduli space of stable genus-zero $J$-holomorphic maps in $X$, representing $\pi_*(\beta)\in H_2(X; \Z)$.
\begin{corollary}
    The moduli space $\M_{0, n}(H; J_H, \beta)$ given the Gromov-convergence topology is compact.
\end{corollary}
\begin{proof}
    Let $u_n:\Sigma\rightarrow H$ be a sequence of $J_H$-holomorphic maps representing the class $\beta$ and write $u_n=(v_n, s_n)$ where $v_n:=\pi\circ u_n$ and $s_n\in H^0(\Sigma; v_n^*H)$. Lemma \ref{lemmaImage} tells us that $s_n\equiv 0$. On the other hand, by Gromov's compactness, we can find a subsequence of $v_n$ which Gromov-converges to $v_\infty$. Thus, $u_n$ has a Gromov-convergent subsequence as $\M_{0, n}(H; J_H, \beta)$ is a Hausdorff space. Therefore, $\M_{0, n}(H; J_H, \beta)$ is sequentially compact and hence compact, as $\M_{0, n}(H; J_H, \beta)$ is a first countable space.
\end{proof}
\begin{remark}
    Here it is necessary to have $\beta\neq0$ as otherwise $\M_{0,n}(H; J_H, 0)$ surjects into the total space of $H\rightarrow X$ via the evaluation map and thus $\M_{0,n}(H; J_H, 0)$ is non-compact.
\end{remark}
\begin{proposition} \label{chartsTwistedGW}
There exists $(G, \T_X, \E_X, s_X)$ and $(G, \T_H, \E_H, s_H)$ Global Kuranishi charts of $\M_{0, n}(X; J, \beta)$ and $\M_{0, n}(H; J_H, \beta)$ respectively such that the following diagram commutes:

\[\begin{tikzcd}
	{\mathcal{E}_X} && {\mathcal{E}_H} \\
	\\
	{\mathcal{T}_X} && {\mathcal{T}_H}
	\arrow["i", hook, from=1-1, to=1-3]
	\arrow[from=1-1, to=3-1]
	\arrow[from=1-3, to=3-3]
	\arrow["{s_X}", bend left = 30pt, from=3-1, to=1-1]
	\arrow[ hook, from=3-1, to=3-3]
	\arrow["{s_H}"',bend right = 30pt, from=3-3, to=1-3]
\end{tikzcd}\]
where the horizontal arrows are inclusion maps.

\end{proposition}
\begin{proof}
    Pick a Hermitian line bundle $\mathcal{L}\rightarrow X$ where $c_1(\mathcal{L})=[\Omega]$ such that $[\Omega]$ is an integral lift of $\w$. We note that, since $H_2(X; \Z)\cong H_2(H; \Z)$, we can use the same compact Lie group $G$, namely $G=U(d+1)$, to define Global Kuranishi charts for both $\M_{0, n}(X; J, \beta)$ and $\M_{0, n}(H; J_H, \beta)$ as per Theorem \ref{AMSI}, where $d:=[\Omega](\beta)$. Now for $u:\Sigma\rightarrow X$ smooth map from a nodal genus-zero Riemann surface, we have a short exact sequence of complex vectorbundles over $\Sigma$
    \begin{equation*}\label{exactseq}
        0\rightarrow u^*TX\rightarrow u^*TH\rightarrow u^*H\rightarrow 0
    \end{equation*}
    Upon fixing a $\C$-basis of holomorphic sections of $u^*\mathcal{L}\rightarrow \Sigma$, we get an induced holomorphic map $\phi:\Sigma\rightarrow\CP^d$. Now noting that tensor product is an exact functor and after choosing $k>>1$ large enough, it follows using Kodaira vanishing theorem that
    \begin{align*}
        0 &\rightarrow H^0(\Sigma; \overline{\operatorname{Hom}}(\phi^* T\CP^d, u^*TX)\otimes_\C\phi^*\mathcal{O}(k))\otimes_\C\overline{H^0(\Sigma; \phi^*\mathcal{O}(k))}\rightarrow \\
        & H^0(\Sigma; \overline{\operatorname{Hom}}(\phi^* T\CP^d, u^*TH)\otimes_\C\phi^*\mathcal{O}(k))\otimes_\C\overline{H^0(\Sigma; \phi^*\mathcal{O}(k))} \rightarrow H^1(\Sigma; v^*H)\rightarrow 0
    \end{align*}
    is an exact sequence and hence the commutativity of the above diagram follows. As all spaces have a $G$-action, and all the vertical arrows are $G$-equivariant maps and the horizontal arrows are inclusions, it follows that the above diagram commutes as $G$-spaces. After applying Theorem \ref{relativeLashof} and noting that the obstruction bundle is a choice of homotopy class of a continuous map into a Grassmannian, which in return can be chosen to be smooth; we get the desired result.
\end{proof}
In order to state and prove the main result of this section, we note the following. Denote by $H^1\rightarrow\T_X$ the quotient bundle $i^*\E_H/\E_X$ over $\T_X$. Following the proof of Lemma \ref{lemmaImage}, for a $J_H$-holomorphic curve $u:\Sigma\rightarrow H$, we write $u=(v, s)$. Then, the proof of Proposition \ref{chartsTwistedGW} implies that we have a splitting
\begin{equation*}
    (\E_H)_u=(\E_X)_v\oplus H^1(\Sigma; v^*H)
\end{equation*}
as $G$-representations. In other words, $H^1\rightarrow\T_X$ is the vectorbundle over $\T_X$ of fibre over $(v,\Sigma, F, \eta)\in\T_X$ given by $H^1(\Sigma; v^*H)$.\\

We recall the following definition.
\begin{definition}[Virtual Fundamental Class]
    Suppose that $(G, \T_X, \E_X,s_X)$ is a Global Kuranishi chart representation of $\overline{\mathcal M}_{0,n}(X;J, \beta)$. We define its virtual fundamental class $[\overline{\mathcal M}_{0,n}(X;J, \beta)]_{vir}:=PD(e_G(\E_X\rightarrow\T_X))\in H^{BM}_{\dim_\R\T-\operatorname{rk}_\R\E-\dim_\R G}(\T/G;\mathbb Q)$, where $e_G$ denotes $G$-equivariant Euler class and $H^{BM}_*$ denote the Borel-Moore homology. Similarly, if $(G, \T_H, \E_H, s_H)$ is a Global Kuranishi chart representation of $\overline{\mathcal M}_{0,n}(H;J, \beta)$, we denote its virtual fundamental class by $[\overline{\mathcal M}_{0,n}(H;J, \beta)]_{vir}:=PD(e_G(\E_H\rightarrow\T_H))$.
\end{definition}
\begin{theorem}\label{TwistedGWNegative}
    $\pi_*[\M_{0, n}(H; J_H, \beta)]_{vir}=e_G(H^1\rightarrow\T_X)\bullet[\M_{0, n}(X; J, \beta)]_{vir}$ where $e_{G}$ denotes the $G$-equivariant Euler class and $\pi_*:H_*(\M_{0, n}(H; J_H, \beta); \mathbb{Q})\rightarrow H_*(\M_{0, n}(X; J, \beta); \mathbb{Q})$.
\end{theorem}
\begin{proof}
    We have the following exact sequence of $G$-equivariant vectorbundles over $\T_X$
    \begin{equation*}
        0\rightarrow\E_X\rightarrow i^*\E_H\rightarrow H^1\rightarrow 0
    \end{equation*}
    and thus we have a splitting
    \begin{equation*}
        i^*\E_H\cong\E_X\oplus H^1
    \end{equation*}
    as $G$-equivariant vectorbundles over $\T_X$. Therefore, upon fixing a $G$-equivariant $\mathbb{Q}$-orientation $\T_X$ and on each of the above $G$-equivariant vectorbundles, and after noting that $[\M_{0, n}(X; J, \beta)]_{vir}=PD(e_G(\E_X\rightarrow\T_X))$ and $[\M_{0, n}(H; J_H, \beta)]_{vir}=PD(e_G(\E_H\rightarrow\T_H))$; the result follows from the functoriality and multiplicative properties of equivariant Euler classes as in Lemma $5.25$ of \cite{abouzaid2021complex} and Appendix $A$ of \cite{hirschi2023properties}.
\end{proof}
\begin{remark}
    Using the results of \cite{abouzaid2021complex}, an immediate corollary of Proposition \ref{chartsTwistedGW} and the proof of Theorem \ref{TwistedGWNegative}, is that we can extend the result of Theorem \ref{TwistedGWNegative}, to have values in any multiplicative complex-oriented generalized cohomology theory. Nevertheless, when considering virtual fundamental classes in generalized cohomology theories, one has to modify the constructed class above by cohomology classes from the Deligne-Mumford space of domains, as in \cite{givental2000wdvv} and \cite{abouzaid2024gromov}, in order to achieve the expected Gromov-Witten axioms. 
\end{remark}

\subsubsection{$H$ is Convex}
We recall the following definition found in \cite{lian1997mirror}.
\begin{definition}[Convex Vectorbundle]\label{convex}
    Let $V\rightarrow\CP^n$ be a holomorphic $\mathbf{T}$-equivariant vectorbundle. $V$ is said to be convex if 
    \begin{enumerate}
        \item $e_\mathbf{T}(V)\in H_\mathbf{T}^*(\CP^n; \mathbb{Q})$ is invertible.
        \item $H^1(\Sigma; u^*V)=\{0\}$ for any $u:\Sigma\rightarrow\CP^n$ genus-zero holomorphic curve.
    \end{enumerate}
    where $\mathbf{T}=(S^1)^{\operatorname{rk}_\C V}$ is the torus. 
\end{definition}
In order to be to give a geometric realization of the results of this subsection, we assume that $(X, \w, J)$ is a compact K$\ddot a$hler manifold. 
For the same reason as above, we restrict our attention to the case of line bundles as we are only dealing with genus-zero curves. Moreover, we weaken the above definition and say $\pi:(H, \nabla)\rightarrow X$ is a $\textbf{convex Holomorphic Line bundle}$ over $X$ if only the second item of Definition \ref{convex} holds.\\

 Let $\sigma\in H^0(X; H)$ be a holomorphic section of $H$ and assume that $\sigma$ is transverse to the zero-section. Namely, 
 \begin{equation} \label{transverseSection}
     (\nabla\sigma)_x:T_xX\rightarrow H_x
 \end{equation}
 is surjective for every $x\in\sigma^{-1}(0)$. Denote by $Y:=\sigma^{-1}(0)$ and note that $(Y, i^*\w, i^*J)$ is a compact K$\ddot a$hler manifold where $i:Y\hookrightarrow X$ is the canonical inclusion. 
 \begin{example}
     For instance, if $X=\CP^4$ and $H=\mathcal{O}(5)\rightarrow X$. Then, $Y$ is a Calabi-Yau three-fold.
 \end{example}
 Fix $\beta\in H_2(Y; \Z)$ and note that $\M_{0 , n}(H; J_H, (\pi_*)^{-1}(i_*\beta))$ is a Hausdorff, non-compact space (c.f. \cite{mcduff2012j} or Example \ref{example} above). Similarly as above and since $i^*\Omega(\beta)=\Omega(i_*\beta)=:d$, we can give both $\M_{0,n}(Y; i^*J, \beta)$ and $\M_{0,n}(X; J, i_*\beta)$ Global Kuranishi charts having the same symmetry group $G=U(d+1)$ using Theorem \ref{AMSI}. \\
 Now we make the following observation. The surjectivity of the operator in \ref{transverseSection} at every point of $Y$, implies that 
\begin{equation} \label{exactSeq2}
    0\rightarrow TY\rightarrow TX|_Y\rightarrow H|_Y\rightarrow 0
\end{equation}
is an exact sequence over $Y$. 
\begin{corollary} \label{chartGW+ve}
    There exists $(G, \T_Y, \E_Y, s_Y)$ and $(G, \T_X, \E_X, s_X)$ Global Kuranishi charts of $\M_{0,n}(Y; i^*J, \beta)$ and $\M_{0,n}(X; J, i_*\beta)$ respectively such that the following diagram commutes:
    \[\begin{tikzcd}
	{\mathcal{E}_Y} && {\mathcal{E}_X} \\
	\\
	{\mathcal{T}_Y} && {\mathcal{T}_X}
	\arrow[hook, from=1-1, to=1-3]
	\arrow[from=1-1, to=3-1]
	\arrow[from=1-3, to=3-3]
	\arrow["{s_Y}", bend left = 30pt, from=3-1, to=1-1]
	\arrow[ hook, from=3-1, to=3-3]
	\arrow["{s_X}"',bend right = 30pt, from=3-3, to=1-3]
\end{tikzcd}\]
where the horizontal arrows are inclusion maps.
\end{corollary}
\begin{proof}
    With the above discussion in mind and using the above exact sequence \ref{exactSeq2}, we can follow the same line of proof as in Proposition \ref{chartsTwistedGW}.
\end{proof}
\begin{definition} \label{H0bundle}
    Denote by $\T:=\{(u, \Sigma, F, \eta)\}$ such that:
    \begin{enumerate}
        \item $u:\Sigma\rightarrow H$ is a smooth map satisfying $\Bar{\partial}_{J_H}u+\pi^*(\langle\eta\rangle\circ d\phi_F)=0$.
        \item $(\pi\circ u,\Sigma, F, \eta)\in\T_X$.
    \end{enumerate}
\end{definition} 
\begin{lemma}
    $\T$ is a topological $G$-manifold.
\end{lemma}
\begin{proof}
    Noting that $\pi_*:H_2(H; \Z)\rightarrow H_2(X;\Z)$ is an isomorphism and with the discussion above in mind, namely we can use the same basis $F$ to frame the domains of $u$ and $\pi\circ u$ simultaneously; it suffice to show that we can always find an $\eta$ such that the linearization of
    \begin{equation} \label{twistedEq}
        \Bar{\partial}_{J_H}u+\pi^*(\langle\eta\rangle\circ d\phi_F)
    \end{equation}
    at $u$ is surjective. Once again, for $u:\Sigma\rightarrow H$ smooth map, we write $u=(\pi\circ u, s)$ where $s\in\Gamma(\Sigma; (\pi\circ u)^*H)$ and in particular
    \begin{equation} \label{EqutoLinear}
        \Bar{\partial}_{J_L}u=(\Bar{\partial}_J(\pi\circ u), \nabla^{0,1}s).
    \end{equation}
    where here we abuse notation and denote $((\pi\circ u)^*\nabla)^{0,1}$ by $\nabla^{0,1}$.
    Therefore, 
    \begin{equation*}
        \operatorname{coker}(D_u\Bar{\partial}_{J_H})\leq\operatorname{coker}(D_{\pi\circ u}\Bar{\partial}_J)\oplus H^1(\Sigma; (\pi\circ u)^*H).
    \end{equation*}
  Now we get the desired result after noting that $H^1(\Sigma; (\pi\circ u)^*H)=\{0\}$ and that the operator $\langle.\rangle\circ d\phi_F$ surjects onto $\operatorname{coker}(D_{\pi\circ u}\Bar{\partial}_J)$ and applying Theorem $6.1$ of \cite{abouzaid2021complex}.
\end{proof}
Notice that, after using Theorem \ref{relativeLashof} twice and arguing as in the proof of Proposition \ref{chartsTwistedGW}, we abuse notation and assume that $\T$ is a smooth $G$-manifold and $\T_X\subseteq\T$ is a smooth submanifold.
Now denote by $N$ the normal bundle of $\T_X\subseteq\T$.
\begin{corollary}
    $N\rightarrow\T_X$ is a $G$-equivariant vectorbundle of fibre $H^0(\Sigma;  u^*H)$ over $(u, \Sigma, F, \eta)\in\T_X$.
\end{corollary}
\begin{proof}
    As $G$ is a compact Lie group, we can find a $G$-invariant metric on $\T$, which in return gives us a $G$-equivariant vectorbundle structure on $N$. As for the second claim, since we used the same $G$-representation to get a $G$-smoothing on both of $\T_X$ and $\T$, and by the second item of Definition \ref{H0bundle}, it suffice to compute $\ker(D_{(u,s)}\Bar{\partial}_{J_H})$ and $\ker(D_{u}\Bar{\partial}_{J})$. Using Equation \ref{EqutoLinear}, it follows that,
    \begin{equation*}
        \ker(D_{(u,s)}\Bar{\partial}_{J_H})=\ker(D_{u}\Bar{\partial}_{J})\oplus H^0(\Sigma; u^*H)
    \end{equation*}
    and hence the desired result.
\end{proof}
 Denote by $\tilde{\sigma}:\T_X\rightarrow N$ the map given by $(u, \Sigma, F, \eta)\mapsto (\sigma\circ u,\Sigma, F, \eta)$.
\begin{lemma}\label{H0section}
    $\tilde{\sigma}$ is a well-defined section of $N$.
\end{lemma}
\begin{proof}
    Let $(u, \Sigma, F, \eta)\in\T_X$ and note that $\tilde{\sigma}(u, \Sigma, F, \eta)=((u, u^*\sigma), \Sigma, F, \eta)$. We calculate as follows: 
    \begin{align*}
       &  \Bar{\partial}_{J_H}(\sigma\circ u)+\pi^*(\langle\eta\rangle\circ d\phi_F) \\
       &  =(\Bar{\partial}_Ju+\langle\eta\rangle\circ d\phi_F, u^*\nabla^{0, 1}\sigma)\\
       & = 0
    \end{align*}
    as $(u, \Sigma, F, \eta)\in\T_X$ and $\sigma\in H^0(X; H)$.
\end{proof}
\begin{theorem} [Quantum Lefschetz Principle]\label{QuanLefsch}
    $[\M_{0, n}(Y; i^*J, \beta)]^{vir}=i^*([\M_{0, n}(X; J, i_*\beta)]^{vir}\cup e_{G} (N))$ where $e_G$ denotes the $G$-equivariant Euler class and $i^*:H^*(\M_{0,n}(X; J, i_*\beta); \mathbb{Q})\rightarrow H^*(\M_{0,n}(Y; i^*J, \beta))$.
\end{theorem}
\begin{proof}
    We first observe that if $u:\Sigma\rightarrow X$ is a smooth map then, the image of $u$ is in $Y$ if and only if $\sigma\circ u(z)=(u(z), 0)$ for all $z\in\Sigma$. Namely, $u^*\sigma(z)=(z, 0)\in u^*L$ for all $z\in\Sigma$. Now from the commutativity of the diagram in Corollary \ref{chartGW+ve}, it follows that $s^{-1}_X(0)\cap\tilde{\sigma}^{-1}(0)=s_Y^{-1}(0)$. Now the result follows from Lemma \ref{H0section}, and the functoriality and multiplicativity properties of Euler classes.
\end{proof}

\section{Global Kuranishi Charts for the Moduli Space of Pseudo-holomorphic Disks}
In this section and throughout the rest of the paper we fix an $\w$-compatible almost complex structure $J$ on $X$ and $L\subset X$ a closed, connected, relatively-spin Lagrangian submanifold. Let $\beta\in\pi_2(X, L)$ be a relative-spherical class and $k, l$ non-negative integers. The goal of this section is to prove the following Theorem.
\begin{theorem}\label{ChartsDisks}
    The moduli space $\M_{k+1, l}(X, L; J, \beta)$ of stable genus-zero bordered $J$-holomorphic maps with $(k+1)$-cyclically ordered boundary marked points, $l$-interior marked points and boundary condition on $L$, representing the class $\beta$; admits a normally complex Global Kuranishi chart $(G, \T, \E, s)$. Together with a smooth submersion $\T\rightarrow X^l\times L^{k+1}$ extending the usual evaluation map.
\end{theorem}
\begin{remark}
    By $\textit{extending the usual evaluation map}$ we mean the following. If we denote by $(G, \T, \E, s)$ the Global Kuranishi chart as in Theorem \ref{ChartsDisks} then we can find $G$-invariant neighborhood of $s^{-1}(0)\subseteq\T$ and a homeomorphism $h$ on such neighborhood, such that $ev\circ h$ is a submersion where $[ev]:\M_{k+1, l}(X, L; J, \beta)\rightarrow X^l\times L^{k+1}$ is the evaluation map induced from $ev$ after taking $G$-quotient. 
\end{remark}

\subsection{Preliminary Constructions}
\begin{definition}
    Let $(\Sigma,\partial\Sigma)$ be a genus-zero bordered Riemann surface and $\mathcal{L}\rightarrow\Sigma$ be a complex line bundle. $\mathcal{L}$ is said to be holomorphic if $\mathcal{L}|_{(\Sigma\setminus\partial\Sigma)}\rightarrow\Sigma\setminus\partial\Sigma$ is a holomorphic line bundle.
\end{definition}
\begin{definition}
    In the setting of the above definition, $\mathcal{L}_{\mathbb{R}}\rightarrow\partial\Sigma$ is said to be totally real line bundle of $\mathcal{L}$ if $\mathcal{L}_{\mathbb{R}}\rightarrow\partial\Sigma$ is a real line bundle such that, $\mathcal{L}_{\mathbb{R}}\otimes_{\mathbb{R}}\mathbb{C}=\mathcal{L}|_{\partial\Sigma}$. 
\end{definition}
\begin{lemma}\label{Integral symplectic form}
There exist a positive integer $N$ and a closed $2$-form $\Omega\in\Omega^{2}(X)$ such that $\Omega$ vanishes on a neighborhood of $L$, and is nondegenerate and compatible with $J$ wherever it is nonzero. Moreover, the relative de Rham class $[\Omega]\in H^{2}(X,L;\mathbb R)$ is integral, and its image under the natural map $H^{2}(X,L;\mathbb R)\longrightarrow H^{2}(X;\mathbb R)$ is equal to $N[\omega]$.
Equivalently, under the canonical identification $H^{2}(X,L;\mathbb R)\cong \widetilde H^{2}(X/L;\mathbb R)$,
if $q:X\to X/L$ denotes the quotient map, then $q^{*}[\Omega]=N[\omega]$.
\end{lemma}
\begin{proof}
Upon setting $g=\omega(. \ , J.)$ as a Riemannian metric on $X$, it follows by Weinstein's neighborhood theorem, that there exists $\epsilon > 0$ so that $\omega =d\lambda$ on $\mathbf{D}_{\epsilon}(L)$, a disc bundle of radius $\epsilon$ over $L$. Now choose a smooth bump function $\phi: X\longrightarrow [0, 1]$ where $\phi\equiv 0$ on $\mathbf{D}_{\frac{\epsilon}{3}}(L)$ and $\phi\equiv 1$ on $\mathbf{D}_\epsilon(L)\setminus\mathbf{D}_{\frac{2\epsilon}{3}}(L)$. Set $\Omega'= d(\phi\lambda)$ on $\mathbf{D}_{\frac{2\epsilon}{3}}(L)$ and $\Omega'=\omega$ on $X\setminus \mathbf{D}_{\frac{2\epsilon}{3}}(L)$.
Note that, the non-degeneracy of $\omega$ and its compatibility with $J$ are open conditions. On the other hand, by construction, $\Omega'$ is a closed $2$-form on $X$ which is non-degenerate and compatible with $J$ on $X\setminus\mathbf{D}_{\frac{2\epsilon}{3}}(L)$. After approximating $\Omega'$ using the $C^{\infty}$-norm on $\Omega^2(X)$, by elements in $H^{2}(X, L; \mathbb{Q})$, it follows that possibly after scaling, there exists $\Omega\in\Omega^2(X)$ a closed $2$-form which is compatible with $J$ and non-degenerate on $X\setminus\mathbf{D}_{\frac{2\epsilon}{3}}$ defining an integral class $[\Omega]\in H^2(X, L; \Z)$. 
\end{proof}
\begin{lemma}\label{LineBundleToFrame}
   There exists a Hermitian line bundle $(\mathcal{L}, \nabla)\rightarrow X$ of curvature form $F_{\nabla}=-2\pi i\Omega$ together with a trivialization $\tau: \mathcal{L}|_{U}\cong U\times\C$, where $U$ is an open neighborhood of $L$.
\end{lemma}
\begin{proof}
    Set $N:= X\setminus L$ and note that $N$ is a non-compact smooth manifold. We denote by $S$ the group of isomorphism classes of pairs $(\mathcal{L}, \tau)$ where $\mathcal{L}\rightarrow N$ is a complex line bundle and $\tau:\mathcal{L}|_{U}\cong U\times\C$ a trivialization at infinity. That is, $U\subset N$ is an open set where $N\setminus U$ is compact. The group structure on $S$ is given by tensor product. Arguing by compactly-supported Cech cohomology, we have $S\cong H^2_c(N; \Z)$ given by the relative first Chern class. By the construction of $\Omega$ in the proof of Lemma \ref{Integral symplectic form}, it follows that, $\Omega\in\Omega_c^2(N)$ such that $[\Omega]\in H^2_c(N; \Z)$. Therefore, there exists $(\mathcal{L}, \tau)\in S$ such that $c_1(\mathcal{L}, \tau)=[\Omega]$.\\
    Now let $\nabla$ be any Hermitian connection on $\mathcal{L}$ which is trivial with respect to $\tau$ at infinity. Then, $\frac{i}{2\pi} F_{\nabla}=\Omega+d\theta$ for some $\theta\in\Omega^1_c(N)$. Setting $\nabla^{'}:=\nabla+i\theta$, we get $F_{\nabla^{'}}=-2\pi i\Omega$. Now the result follows after extending $\mathcal{L}$ over $L$ via $\tau$.
\end{proof}
Let $u:(\Sigma, \partial\Sigma)\rightarrow (X, L)$ be a smooth map from a genus-zero bordered Riemann surface and denote by $\mathcal{L}_u:=u^*\mathcal{L}$ where $\mathcal{L}\rightarrow X$ is as in Lemma \ref{LineBundleToFrame}. Also following the same notation as in the proof of Lemma \ref{LineBundleToFrame} and after pulling back the trivialization $\tau$ by $u$, we get a uniquely determined $\mathcal{L}_\R\rightarrow\partial\Sigma$, totally-real line bundle of $\mathcal{L}_\R$.\\

Now consider the double of $(\Sigma, \partial\Sigma)$, namely let $\Tilde{\Sigma}= \Sigma\cup_{\partial\Sigma}\Bar{\Sigma}$ where $\Bar{\Sigma}$ is the complex conjugate. Thus, $\Tilde{\Sigma}$ is a (nodal) genus-zero orientable nodal surface. Moreover, from this construction we have an anti-holomorphic involution $\iota:\Tilde{\Sigma}\longrightarrow \Tilde{\Sigma}$ with non-empty fixed locus, Fix$(\iota)\cong \partial\Sigma$. In addition, upon pulling-back the trivialization $\tau$ by $u$, we consider $\Tilde{\mathcal{L}}_u\rightarrow\Tilde{\Sigma}$ the double of $\mathcal{L}_u$. We note the following easy Lemma.

\begin{lemma}
Up to biholomorphic conjugacy, there are exactly two anti-holomorphic involutions on $\mathbb{CP}^{1}$: one with nonempty fixed-point set, represented by complex conjugation, and one with empty fixed-point set, represented by the antipodal real structure $z\longmapsto-\frac{1}{\overline z}$.
\end{lemma}
\begin{proof}
Indeed, using brute-force calculations one can see that in affine coordinates, anti-holomorphic involutions on $\mathbb{CP}^{1}$ are of the form $z\mapsto\frac{1}{\Bar{z}}$ or $z\mapsto\frac{-1}{\Bar{z}}$, distinguished by whether it has a non-empty fixed locus or not. This follows from the fact that any anti-holomorphic involution on $\CP^1$ is induced by an invertible complex anti-linear map $\C^2\rightarrow\C^2$, which in return is determined uniquely, up to multiplication by $\C^*$, by an $A\in GL_2(\C)$. With this in mind and for an involution $\iota$, let $A\in GL_2(\C)$ be such that, for $[v]=[v_0:v_1]\in\CP^1$, $\iota([v])\equiv[Av]$. As $\iota$ is an involution, it follows that, $A\Bar{A}v=\lambda I_{\C^2}$ such that $\lambda\in \C^*$. After multiplying from the right by $A$, it follows that $\lambda$ is a non-zero real number. After rescaling $A$ by non-zero real number if necessary, we can assume that $\lambda\in\{-1,1\}$. Complex conjugation, $z\mapsto \Bar{z}$ corresponds to the case where $\lambda=1$, while $z\mapsto\frac{-1}{\Bar{z}}$ corresponds to the case where $\lambda=-1$.   
\end{proof}
\begin{lemma}\label{DegreeOfLineBundle}
    $\Tilde{\mathcal{L}}_u\rightarrow \Tilde{\Sigma}$ is a holomorphic line bundle with a unique holomorphic structure, satisfying $\langle c_{1}(\Tilde{\mathcal{L}}_u), [\Tilde{\Sigma}]\rangle=2(u^{*}\Omega) [\Sigma, \partial\Sigma]$.
\end{lemma}
\begin{proof}
     Note that, $\Tilde{\mathcal{L}}_u$ is a Hermitian line bundle over $\Tilde{\Sigma}$ and thus, it is a holomorphic line bundle, as $H^{0,2}(\Tilde{\Sigma})={0}$, with a unique holomorphic structure since $\Tilde{\Sigma}$ is simply-connected. \\
    To prove the second claim, we give two proofs:\\
    \begin{enumerate}
        \item Notice that, we can represent the relative first Chern number of $\langle c_{1}(\mathcal{L}_{u}, \tau), [\Sigma, \partial\Sigma]\rangle$ with the signed count of number of zeros of a section of $\mathcal{L}_{u}$ transverse to the zero-section and being non-zero on the boundary. On the other hand, $\langle c_{1}(\mathcal{L}_{u}, \tau), [\Sigma, \partial\Sigma]\rangle=(u^{*}\Omega) [\Sigma, \partial\Sigma]$. Thus, upon doubling this section to a section of $\Tilde{\mathcal{L}}_u$, we get the second claim.

        \item Alternatively, we calculate on the double of each disk component of $\Sigma$ as follows:\\
    Denote by $U=\mathbb{D}(1+\epsilon)\setminus\mathbb{D}(1-\epsilon)$ and for $i=1,2$ denote by $F_{h_{i}}$ the Chern curvature of $\mathcal{L}_{i}=(\mathcal{L}_u, h_{i})$. For $v_{i}\in\Gamma (\mathcal{L}_{i})$ we have $h_{i}(v, v)=e^{\psi_{i}}||v||^{2}$ with $F_{h_{i}}=-\partial\Bar{\partial}\psi_{i}$. On $U$, $h_{2}=e^{\phi}h_{1}$ and $v_{2}=z^{k}v_{1}$ where $k:=2(u^*\Omega)[\mathbb{D}, S^1]$ and hence, \begin{equation*}
        h_{2}(v_{2}, v_{2})=e^{\psi_{2}}||v_{2}||^{2}=e^{\psi_{2}}z^{2k}||v_{1}||^{2}=e^{\psi_{2}+k\log|z|^{2}}||v_{1}||^{2}=e^{\psi_{2}-\psi_{1}+k\log|z|^{2}}h_{1}(v_{1}, v_{1})
    \end{equation*}
        Therefore, $\phi=\psi_{2}-\psi_{1}+k \log|z|^{2}$.\\
        Using Stokes' Theorem, we calculate the total curvature of the induced metric $g$ on $\tilde{\mathcal{L}}_u$ as follows: 
        \begin{equation*}
            \int_{S^2}F_{g}=\int_{\mathbb{D}(1-\epsilon)}F_{h_{1}}+\int_{U}F_{g}+\int_{\mathbb{D}(1-\epsilon)}F_{h_{2}}
            =-\int_{S^{1}_{-}}\Bar{\partial}\psi_{1}+\int_{U}F_{g}+\int_{S^{1}_{+}}\Bar{\partial}\psi_{2}
        \end{equation*}
        where the $S^{1}_{+}$ and $S^{1}_{-}$ denotes the boundary circles of the upper and lower $\mathbb{D}(1-\epsilon)$ respectively. Using Stokes' Theorem one more time we have 
        \begin{equation*}
            \int_{U}F_{g}=-\int_{S^{1}_{+}}\Bar{\partial}\log(\phi_{1}h_{1}+\phi_{2}h_{2})+\int_{S^{1}_{-}}\Bar{\partial}\log(\phi_{1}h_{1}+\phi_{2}h_{2})
        \end{equation*}
        Noting that, on $S^{1}_{+}$ we have 
        \begin{equation*}
            -\partial\Bar{\partial}\log h_{2}=-\partial\Bar{\partial}\psi_{2}+k\log|z|^{2}
        \end{equation*}
        it follows that,
        \begin{equation*}
            \int_{S^{1}_{+}}\Bar{\partial}\log(\phi_{1}h_{1}+\phi_{2}h_{2})=-\int_{S^{1}_{+}}\Bar{\partial}\psi_{2}-k\int_{S^{1}_{+}}\Bar{\partial}\log|z|^{2}
        \end{equation*}
        Thus by symmetry, it follows that, 
        \begin{equation*}
            \int_{S^2}F_{g}=-2k\int_{S^{1}_{+}}\Bar{\partial}\log(r)=-2k\int_{S^{1}} \Bar{\partial}\log(r)
        \end{equation*}
        as $\log(r)$ is a harmonic function and hence the above integral is independent of the choice of the radius of $S^{1}_{+}$.
        Thus, to find the total curvature of $g$ we use polar coordinates $z=re^{i\theta}$, $\Bar{z}=re^{-i\theta}$ and thus $d\Bar{z}=e^{-i\theta}dr-ire^{i\theta}d\theta$ and $\partial_{\Bar{z}}=\partial_{r}\frac{\partial r}{\partial\Bar{z}}+\partial_{\theta}\frac{\partial\theta}{\partial\Bar{z}}$. Therefore, 
        \begin{equation*}
            \Bar{\partial}\log(r)=\frac{\partial \log(r)}{\partial\Bar{z}}d\Bar{z}=\frac{\partial \log(r)}{\partial r}\frac{\partial r}{\partial\Bar{z}}(-ire^{i\theta})d\theta
        \end{equation*}
        as $dr=0$ on the unit circle. Noting that, $r^{2}=z\Bar{z}$, it follows that, $2r\frac{\partial r}{\partial\Bar{z}}=z$ and hence, $\frac{\partial r}{\partial\Bar{z}}=\frac{z}{2r}=\frac{1}{2}e^{-i\theta}$. Thus, 
        \begin{equation*}
            \Bar{\partial}\log(r)= \frac{1}{2r}e^{i\theta}(-ire^{-i\theta})d\theta=-\frac{i}{2}d\theta
        \end{equation*}
        and hence,
        \begin{equation*}
            -2k\int_{S^{1}}\Bar{\partial}\log(r)=2k\pi i
        \end{equation*}
        By Chern-Weil Formula, the result follows.
    \end{enumerate}
\end{proof}
\begin{lemma}\label{LiftingInvolutions}
    There exists a lift $\Tilde{\iota}:\Tilde{\mathcal{L}}_u\rightarrow\Tilde{\mathcal{L}}_u$ of $\iota:\Tilde{\Sigma}\rightarrow \Tilde{\Sigma}$. Moreover any two lifts differ by a multiple of a complex number $\alpha\in \mathbb{C}^{*}$ of norm 1.
\end{lemma}
\begin{proof}
    On a single trivializing set, we argue as follows. For existence, we consider $\Tilde{\iota}(x, v)=(\iota(x), \Bar{v})$. Now for $i=1,2$, let $\Tilde{\iota}_{i}$  be two lifts of $\iota$ and consider $\Tilde{\iota}_{1} \circ \Tilde{\iota}_{2}(x, v)=(x, g(v))$ where $g:\Tilde{\Sigma}\rightarrow \operatorname{GL}_{1}(\mathbb{C})$ as $\Tilde{\iota}_{1} \circ \Tilde{\iota}_{2}$ is holomorphic. By compactness of $\Tilde{\Sigma}$, it follows that any lift of $\iota$ to an anti-holomorphic on $\Tilde{\mathcal{L}}_u$ is of the form $(x, v)\rightarrow (\iota(x), \alpha\Bar{v})$, where $\alpha$ is a non-zero complex number. Moreover, as $\Tilde{\iota}_{i}$ are involutions, it follows that $\alpha$ is of norm 1. As for the general case, namely after incorporating transition maps, we refer to Theorem $3.3.13$ of \cite{katz2006enumerative}.
\end{proof}
Denote by $H^{i}(\Tilde{\Sigma}; \Tilde{\mathcal{L}}_u)^{\iota, \Tilde{\iota}}:=\{s\in H^i(\Tilde{\Sigma}; \Tilde{\mathcal{L}}_u):\Tilde{\iota}\circ s\circ\iota=s\}$.
\begin{lemma}\label{DimensionLemma}
We have an isomorphism $H^{i}(\Sigma, \partial\Sigma; \mathcal{L}_u, \mathcal{L}_\R)\cong H^{i}(\Tilde{\Sigma}; \Tilde{\mathcal{L}}_u)^{\iota, \Tilde{\iota}}$ as $\mathbb{R}$-vectorspaces. In particular, $\dim_{\mathbb{R}}H^{i}(\Sigma, \partial\Sigma; \mathcal{L}_u, \mathcal{L}_\R)=\dim_{\mathbb{C}}H^{i}(\Tilde{\Sigma}; \Tilde{\mathcal{L}}_u)$.
\end{lemma}
\begin{proof}
A proof of such result can be found in Section $3.4$ of \cite{katz2006enumerative}, for which we highlight the main construction that goes in such a proof.\\
    We argue using Cech cohomology construction. Indeed, let $\mathcal{A}:=\{U_{i}\}$ be an acyclic cover for the sheaf of local holomorphic sections of $\mathcal{L}_u$ with boundary values in $ \mathcal{L}_{\mathbb{R}}$ on $\Sigma$. We double $\mathcal{A}$ to $\Tilde{\mathcal A}:=\{\Tilde{U}_{i}\}$ where $\Tilde{U}_{i}:= U_{i}\sqcup\iota(U_{i})$ if $U_{i}\cap\partial\Sigma=\phi$ and $\Tilde{U}_{i}:=U_{i}\cup_{U_{i}\cap\partial\Sigma}\iota(U_{i})$ if $U_{i}\cap\partial\Sigma\neq\phi$. Upon possibly further restrictions of $\mathcal A$, we can assume that $\Tilde{\mathcal A}$ is acyclic for the sheaf of local holomorphic sections of $\Tilde{\mathcal{L}}_u$ on $\Tilde{\Sigma}$. Now for $U\subseteq\Sigma$ open set and $s\in\Gamma(U; \mathcal{L}_u, \mathcal{L}_{\mathbb{R}})$, $\sigma(s):=\Tilde{\iota}\circ s\circ\iota$ is a local section $\Tilde{\mathcal{L}}_u$ on $\tilde U$. On the other hand, for a local section $t$ of $\Tilde{\mathcal{L}}_u$, we have $t=\sigma(t)$ if and only if $t\in\Gamma(V; \mathcal{L}_u, \mathcal{L}_{\mathbb{R}})\subseteq\Gamma(V; \Tilde{\mathcal{L}}_u)$ where $V$ is an open set of $\Sigma$. Hence, we get an identification between the $\sigma$-fixed locus of the Cech complex of $\Tilde{\mathcal{L}}_u$ given by $\Tilde{\mathcal A}$ with the Cech complex of $(\mathcal{L}_u, \mathcal L_\R)$ given by $\mathcal A$.\\
    As for the second claim, $\sigma$ induces an involution on $H^i(\tilde\Sigma; \mathcal L_u)$ as an $\R$-vectorspace with $H^{i}(\Tilde{\Sigma}; \Tilde{\mathcal{L}}_u)^{\iota, \Tilde{\iota}}$ being the eigenspace corresponding to eigenvalue $1$.
\end{proof}
\subsection{Space of Framings}
Consider $(\mathbb{CP}^{d}, \omega_{FS}, J_{std})$ with its Fubini-Study symplectic form and its standard integrable almost complex structure and let $\alpha\in\pi_2(\CP^d, \RP^d)$.
\begin{definition}[Space of Framing for Disks]\label{spaceFramings}
    Denote by $\F_{k+1, l}(\alpha):=\{\phi:(\Sigma, \partial\Sigma; \underline{z}, \underline{w})\rightarrow(\CP^d, \RP^d)\}$ satisfying the following conditions:
    \begin{enumerate}
        \item $(\Sigma, \partial\Sigma)$ is a genus-zero bordered Riemann surface.
        \item $\underline{z}=(z_0,\dots,z_k)$ are distinct boundary, non-nodal, marked points which are cyclically ordered via the counter clock-wise orientation on $\partial\Sigma$.
        \item $\underline{w}=\{w_1, \dots, w_l\}$ are distinct, non-nodal, interior marked points.
        \item $\phi$ is a smooth map satisfying $\Bar{\partial}_{J_{std}}\phi=0$ on each irreducible component of $\Sigma$, representing the class $\alpha$.
        \item $\operatorname{Aut}((\Sigma, \partial\Sigma; \underline{z}, \underline{w}); \phi)$ is trivial.
        \item $H^1(\Sigma, \partial\Sigma; \phi^*\mathcal{O}(1), \phi^*_{\partial\Sigma}\mathcal{O}(1)_{\mathbb{R}})=\{0\}$.
    \end{enumerate}
\end{definition}

\begin{remark}
    We note that the domains of elements of $\mathcal{F}_{k+1, l}(\alpha)$ are not fixed. Hence, for $\phi\in\F_{k+1, l}(\alpha)$ and after pulling-back the Euler sequence and using items 5 and 6 of Definition \ref{spaceFramings}, we only get a smooth neighborhood $U\subseteq\F_{k+1, l}(\alpha)$ of $\phi$.
\end{remark}
\begin{proposition}
    The space of framing $\F_{k+1, l}(\alpha)$ is a smooth $O(d+1)$-manifold with corners. 
\end{proposition}
\begin{proof}
     Following \cite{zernik2017moduli}, we have $\M_{k+1, l}(\CP^d, \RP^d; J_{std}, \alpha)$, the stable map compactification of the moduli space of $J_{std}$-holomorphic disks representing the class $\alpha$, is a smooth $U(d+1)^{\Z/2}$-orbifold with corners where $U(d+1)^{\Z/2}=O(d+1)$ is the maximal subgroup of $U(d+1)$ commuting with the complex conjugation on $\CP^d$. Furthermore, the orbifold isotropy group at each point of $\M_{k+1, l}(\CP^d, \RP^d; J_{std}, \alpha)$ is isomorphic to the automorphism group of the underlying curve \cite{zernik2017moduli}. On the other hand, $\mathcal{F}_{k+1, l}(\alpha)\subset \M_{k+1, l}(\CP^d, \RP^d; J_{std}, \alpha)$ is the isotropy-free part of $\M_{k+1, l}(\CP^d, \RP^d;J_{std}, \alpha)$ and hence the result.
\end{proof}

Now set $d:=\Omega(\beta)$, where $\Omega$ is as in Lemma \ref{Integral symplectic form} and $\beta\in\pi_2(X, L)$ as above, and let $\alpha:=d[\CP^1, \RP^1]\in\pi_2(\CP^d, \RP^d)$. Following the same notation as in the above section, for $u:(\Sigma,\partial\Sigma)\rightarrow(X, L)$ smooth map representing the class $\beta$, we have an induced holomorphic line bundle together with a totally-real sub-bundle $(\mathcal{L}_u, \mathcal{L}_\mathbb{R})\rightarrow(\Sigma, \partial\Sigma)$.
\begin{lemma}\label{H1=0}
 $\dim_{\mathbb{R}}H^0(\Sigma, \partial\Sigma; \mathcal{L}_u, \mathcal{L}_\R)=d+1$ and $H^1(\Sigma, \partial\Sigma; \mathcal{L}_u, \mathcal{L}_\R)=\{0\}$.
\end{lemma}
\begin{proof}
Lemma \ref{DegreeOfLineBundle} tells us that $\mathcal{L}_u$ has non-negative degree on each of the irreducible component of $\Sigma$. Now using Lemma \ref{DimensionLemma} and after applying Riemann-Roch, Serre Duality and Kodaira Vanishing Theorem, the result follows.
\end{proof}
\begin{definition}[Framed Marked Disks]\label{FramedDisk}
    A framed marked disk, representing the class $\beta\in\pi_2(X, L)$, is a tuple $((\Sigma, \partial\Sigma; \underline{z}, \underline{w}), u, F)$ where:
    \begin{enumerate}
        \item $(\Sigma, \partial\Sigma; \underline{z}, \underline{w})$ is a marked genus-zero bordered Riemann surface. 
        \item $u:(\Sigma, \partial\Sigma)\rightarrow (X, L)$ smooth map with boundary conditions on $L$, representing the class $\beta$ and $\operatorname{Aut}(u, \Sigma, \underline{z}, \underline{w})$ is of finite order.
        \item $F$ is an $\R$-basis of $H^0(\Sigma, \partial\Sigma; \mathcal{L}_u, \mathcal{L_\R})$ such that the matrix $H_F:=(\int_\Sigma\langle f_i, f_j\rangle u^*\Omega)_{i, j}$ is a positive-definite matrix, where $\langle.,.\rangle$ is the induced pairing, via the isomorphism in Lemma \ref{DimensionLemma}, from the Hermitian pairing on $H^0(\tilde\Sigma; \Tilde{\mathcal L}_u)$.
    \end{enumerate}
\end{definition}

\begin{proposition}\label{FramingMap}
    Suppose that $((\Sigma, \partial\Sigma; \underline{z}, \underline{w}), u, F))$ is a framed marked disk representing the class $\beta$ and write $F:=\{f_0,\dots,f_d\}$. Then the induced holomorphic map 
    \begin{align*}
    \phi_F: (\Sigma,\partial\Sigma;\underline{z},\underline{w}) \ & \to (\mathbb {CP}^d, \RP^d) \\
    z \ & \mapsto [f_i(z)]
\end{align*}
is an element of $\F_{k+1, l}(\alpha)$. On the other hand, given $\phi\in\F_{k+1, l}(\alpha)$ then there exists an $\R$-basis of $H^0(\Sigma,\partial\Sigma; \phi^*\mathcal{O}(1), \phi^*_{\partial\Sigma}\mathcal{O}(1)_\R)$, which is  unique up to conjugation by $O(d+1)$, inducing $\phi$. 
\end{proposition}
\begin{proof}
    By construdction of $\phi_F$ and using Lemma \ref{H1=0}, it suffice to show that $\operatorname{Aut}(\phi_F)\equiv\operatorname{Aut}(\phi_F, \Sigma, \underline{z}, \underline{w})$ is trivial. Suppose that $\sigma\in\operatorname{Aut}(\phi_F)$ and recall that $\sigma:(\Sigma, \partial\Sigma)\rightarrow(\Sigma, \partial\Sigma)$ is a homeomorphism such that on each irreducible component $\Sigma_i$, $\sigma_{i}:\Sigma_i\rightarrow\Sigma_i$ is a biholomorphism such that $\phi_F|_{\Sigma_i}\circ\sigma_i=\phi_F|_{\Sigma_i}$. Thus we can double $\Tilde{\sigma}:\Tilde{\Sigma}\rightarrow\Tilde{\Sigma}$ giving us an element in $\operatorname{Aut}(\Tilde{\phi}_F)$, where $\Tilde{\phi}_F:\Tilde{\Sigma}\rightarrow\Tilde{\Sigma}$ is the double of $\phi_F$. In particular, we get an injective group morphism $\operatorname{Aut}(\phi_F)\hookrightarrow\operatorname{Aut}(\Tilde{\phi}_F)$. Now following \cite{zernik2017moduli}, $\F_{k+1, l}(\alpha)$ can be realized as the pullback by the doubling map $\M_{k+1, l}(\CP^d, \RP^d; J_{std}, \alpha)\rightarrow\M_{k+1, 2l}(\CP^d; J_{std}, \alpha+\Bar{\alpha})$, of all elements $\psi\in\M_{k+1, 2l}(\CP^d; J_{std}, \alpha+\Bar{\alpha})$ such that $H^1(\Tilde{\Sigma}; \psi^*\mathcal{O}(1))=\{0\}$. Using Lemma \ref{DimensionLemma}, it follows that $H^1(\Tilde{\Sigma}, \Tilde{\phi}_F^*\mathcal{O}(1))=\{0\}$ which in return implies that $\operatorname{Aut}(\Tilde{\phi}_F)$ is trivial. Therefore, $\phi_F\in\F_{k+1, l}(\alpha)$ as desired. \\
    As for the second statement, let $\phi\in\F_{k+1, l}(\alpha)$ and denote by $\Tilde{\phi}$ its double. As\\ $H^1(\Sigma,\partial\Sigma; \phi^*\mathcal{O}(1), \phi^*_{\partial\Sigma}\mathcal{O}(1)_\R)$ is trivial, it follows from Lemma \ref{DimensionLemma} that $H^1(\Tilde{\Sigma}; \Tilde{\phi}^*\mathcal{O}(1))$ is also trivial and hence we can find a $\C$-basis of $H^0(\Tilde{\Sigma}; \Tilde{\phi}^*\mathcal{O}(1))$, say $F_\C$, inducing $\Tilde{\phi}$. Note that $\Tilde{\Sigma}$ is equipped with an involution $\iota$ such that $\operatorname{Fix}(\iota)\cong\partial\Sigma$ and using Lemma \ref{LiftingInvolutions} we lift $\iota$ to an involution $\Tilde{\iota}$ on $\Tilde{\phi}^*\mathcal{O}(1)\rightarrow(\Tilde{\Sigma}, \iota)$. Possibly after conjugating $F_\C$ by an element from $U(d+1)$, we write $F_\C=\{f_0,\dots,f_d;g_0,\dots,g_d\}$ where $\Tilde{\iota}\circ f_i\circ\iota=f_i$ and $\Tilde{\iota}\circ g_i\circ\iota=-g_i$ for $i=1,\dots,d$. Now $\{f_0,\dots,f_d\}$ is an $\R$-basis of $H^0(\Tilde{\Sigma};\Tilde{\phi}^*\mathcal{O}(1))^{\Tilde{\iota}, \iota}\cong H^{0}(\Sigma,\partial\Sigma; \phi^*\mathcal{O}(1), \phi^*_{\partial\Sigma}\mathcal{O}(1))$ where the isomorphism is given as in Lemma \ref{DimensionLemma}, giving us the second statement.
\end{proof}

\subsection{Proof of Theorem \ref{ChartsDisks}}
\subsubsection{Using \cite{abouzaid2021complex}, \cite{bai2022arnold}}
$\textit{Gromov's Graph Trick}$:\\
 We start by noting that $\mathcal{F}_{k+1, l}(\alpha)$ is a smooth manifold with corners and not necessarily an algebraic space. Indeed, it may have an odd real dimension. For this reason, we apply Gromov's graph trick on $(X\times\mathbb{CP}^{d}, L\times\mathbb{RP}^{d})$, where $d=\Omega(\beta)$ as above.
\\

Let $(E, \nabla)\rightarrow X\times\mathbb{CP}^{d}$ be a Hermitian vectorbundle and denote by 

\begin{equation*}
J_{X\times\mathbb{CP}^{d}}:=p_{X}^{*}J\oplus p_{\mathbb{CP}^{d}}^{*}J_{std} 
\end{equation*} 

the product almost complex structure on $X\times\mathbb{CP}^{d}$, where 
 \begin{equation*}
        p_{X}:X\times\mathbb{CP}^{d}\rightarrow X 
    \end{equation*}
    \begin{equation*}
        p_{\mathbb{CP}^{d}}:X\times\mathbb{CP}^{d}\rightarrow\mathbb{CP}^{d}
    \end{equation*}
    are the natural projection maps. 
    Consider the $\textbf{shearing map}$
 \begin{equation}
 \Psi:E\oplus p_{X}^{*}TX\oplus p_{\mathbb{CP}^{d}}^{*}T\mathbb{CP}^{d}\rightarrow  p_{X}^{*}TX\oplus p_{\mathbb{CP}^{d}}^{*}T\mathbb{CP}^{d} 
 \end{equation}
 satisfying:
 \begin{equation}\label{Equ1}
     J_{X\times\mathbb{CP}^{d}}(\Psi(e, h))=-\Psi(e, J_{X\times\mathbb{CP}^{d}}(h))
 \end{equation}
\begin{equation}\label{Equ2}
    \Psi(e, \Psi(e, h))=0
\end{equation}
Using the Hermitian connection $\nabla$ on $E$, we define $\Phi:T^{vt}E\oplus T^{h}E\rightarrow T^{vt}E\oplus T^{h}E$ fibre-wise by:
\begin{equation}
\Phi(v, h):=(v, h+\Psi(e, h))
\end{equation}
 for every $e\in E$; where if we denote by $\pi:E\rightarrow X\times\CP^d$ the projection map, then
 \begin{equation*}
     T^{vt}E:=\ker(d\pi)
 \end{equation*}
the vertical tangent bundle of $E\rightarrow X\times\CP^d$ and 
\begin{equation*}
    T^hE:=\{ds(v): v\in T(X\times\CP^d), s\in\Gamma(E), \nabla_vs=0\}
\end{equation*}
the horizontal tangent bundle of $E$ associated to $\nabla$.
 We define the $\textbf{$\Psi$-sheared}$ almost complex structure on the total space of $E$ by 
    \begin{equation}
        J_{\Psi}:=\Phi^{-1}\circ J_{X\times\mathbb{CP}^{d}}\oplus i\circ\Phi
    \end{equation}
    where $i$ denotes the fibre-wise multiplication by $\sqrt{-1}$ on $E$. In particular, for $e=(v, h)\in TE\cong T^{vt}E\oplus T^hE$ we have 
    \begin{equation}\label{Equ3}
        J_{\Psi}(v, h)=(iv, J_{X\times\mathbb{CP}^{d}}h+2\Psi(e,J_{X\times\mathbb{CP}^{d}}h))
    \end{equation}
Now for a smooth map $\Tilde{u}:\Sigma\rightarrow E$ where $\Sigma$ is a genus-zero nodal Riemann surface, we write $\Tilde{u}=(u, s)$ where $u:\Sigma\rightarrow X\times\mathbb{CP}^{d}$ is a smooth map and $s\in \Gamma(\Sigma; u^{*}E)$ a smooth section. Using the splitting of $TE\cong T^{vt}E\oplus T^{h}E$ induced by $\nabla$ as above, we have 
    \begin{equation}
        d\Tilde{u}=(du, u^{*}\nabla s)
    \end{equation}
    In particular, the $(0, 1)$-part of $d\Tilde{u}$ with respect to $(j,J_{\Psi})$, where $j$ is an almost complex structure on $\Sigma$, is 
    \begin{equation}\label{0,1-part}
        (\frac{1}{2}(du+J_{\Psi}\circ du\circ j), (u^{*}\nabla)^{0, 1}s)
    \end{equation}
 From the above discussion, we get immediately the following Corollary.
    \begin{corollary}\label{shearedCurves}
        
        Let $\Tilde{u}:\Sigma\rightarrow E$ be a smooth map and write $\Tilde{u}=(u, s)$ as above. Then, $\Tilde{u}$ is $J_{\Psi}$-holomorphic if and only if
        \begin{equation}
            \Bar{\partial}_{J_{X\times\mathbb{CP}^{d}}}u+\Psi(s, J_{X\times\mathbb{CP}^{d}}\circ du\circ j )=0
        \end{equation}
        \begin{equation}
            (u^{*}\nabla)^{0, 1}s=0
        \end{equation}
        and thus $s\in H^{0}(u^{*}E)$.
    \end{corollary}
\begin{proof}
    Follows from Equation \ref{0,1-part} and the fact that $H^{0, 2}(\Sigma)$ is trivial.
\end{proof}
  From now on we set $E$ and $\Psi$ as follows:
\begin{definition} \label{preOb}
    We define
    \begin{equation}
    E:=p_{X}^{*}TX\otimes_{\mathbb{C}} p_{\mathbb{CP}^{d}}^{*}(\Lambda^{0, 1}\mathbb{CP}^{d}\otimes_{\mathbb{C}}\mathcal{O}(N)\otimes_{\mathbb{C}} \overline{H^{0}(\mathbb{CP}^{d}, \mathcal{O}(N))})
    \end{equation}
    where $N>>1$ is some positive integer (chosen in Proposition \ref{ChoosingN}). Together with the Hermitian connection on $E$ induced from the Chern connection on $\mathcal{O}(N)\rightarrow\CP^d$ and after fixing a Hermitian connection on $(X, \w(., J.), J)$.
    \end{definition}
\begin{definition}
    We define $\Psi:E\oplus p_{X}^{*}TX\oplus p_{\mathbb{CP}^{d}}^{*}T\mathbb{CP}^{d}\rightarrow  p_{X}^{*}TX\oplus p_{\mathbb{CP}^{d}}^{*}T\mathbb{CP}^{d} $ by
    \begin{equation} \label{Equ4}
         \Psi(\eta, v_{1}, v_{2})=(\langle\eta\rangle v_{2}, 0)
    \end{equation}
    where $\langle.\rangle:\mathcal{O}(N)\otimes \overline{H^{0}(\mathbb{CP}^{d}, \mathcal{O}(N))}\rightarrow\mathbb{C}$ is the induced Hermitian pairing from the $N^{th}$-tensor power of the Fubini-Study metric on $\mathcal{O}(1)\rightarrow\mathbb{CP}^{d}$.
\end{definition}
Denote by $\mathcal{O}_\R(N)\rightarrow\RP^d$ the unique real line bundle that is fixed by the conjugation action on $\mathcal{O}(N)\rightarrow\CP^d$ and note that, $\mathcal{O}_\R(N)$ is a totally real line bundle of $\mathcal{O}(N)|_{\RP^d}$. Furthermore, using the trivialization of $\mathcal{L}|_L\cong L\times\C$ we denote by $\mathcal{L}_\R\rightarrow L$ the induced totally real line bundle, that is $\mathcal{L}_\R\cong L\times\R$ under the identification $\mathcal{L}|_L\cong L\times\C$.
\begin{definition}
    We define $E_\R\rightarrow L\times\RP^d$ by
     \begin{equation}
         E_\R:=p_{X}^{*}TL\otimes_{\mathbb{R}} p_{\mathbb{CP}^{d}}^{*}(T^{*}\mathbb{RP}^{d}\otimes_{\mathbb{R}}\mathcal{O}_{\R}(N)\otimes_{\mathbb{R}}H^{0}(\mathbb{CP}^{d}, \mathbb{RP}^{d}; \mathcal{O}(N), \mathcal{O}_\R(N))).
    \end{equation}
\end{definition}
\begin{lemma}
    $E_\R$ is a totally-real sub-bundle of $E$ for $J_\Psi$, that is $E_\R\otimes_\R\C\cong E|_{L\times\RP^d}$.
\end{lemma}
\begin{proof}
    We calculate point-wise as follows. Let $e\in E$, $x\in X$ and $p\in\CP^d$ and let $e'\in E_e$, $v_1\in TX$ and $v_2\in T\CP^d$, then using Equations \ref{Equ3} and \ref{Equ4} we have
    \begin{equation*}
        (J_{\Psi})_{(e, x, p)}(e', v_{1}, v_{2})=(\sqrt{-1}e' ,J_{x}v_{1}+\langle e\rangle v_{2}, (J_{std})_{p}v_{2}).
    \end{equation*}
    Now the claim follows from Equations \ref{Equ1} and \ref{Equ2}.
\end{proof}

Now suppose that $\Tilde{u}:(\Sigma, \partial\Sigma)\rightarrow (E, E_\R)$ smooth map with real-boundary condition, that is $\Tilde{u}(\partial\Sigma)\subseteq E_\R$, where $(\Sigma,\partial\Sigma)$ is a genus-zero bordered Riemann surface. As above, using the Hermitian connection on $E$, we write $\Tilde{u}=(u, s)$ where $u:(\Sigma, \partial\Sigma)\rightarrow (X\times\CP^d, L\times\RP^d)$ and $s\in\Gamma(\Sigma, \partial\Sigma; u^*E, u^*E_\R)$ a smooth section with boundary values in $E_\R$.
\begin{corollary}
    $\Tilde{u}$ is $J_{\Psi}$-holomorphic with real-boundary condition if and only if on each spherical component of $\Sigma$ we have
    \begin{equation}\label{Equ5}
        \Bar{\partial}_{J_{X\times\mathbb{CP}^{d}}}u+\Psi(s, J_{X\times\mathbb{CP}^{d}}\circ du\circ j )=0
    \end{equation}
    \begin{equation}
        s\in H^{0}(u^{*}E)
    \end{equation}
    and on each disk component we have
    \begin{equation}\label{Equ6}
         \Bar{\partial}_{J_{X\times\mathbb{CP}^{d}}}u+\Psi(s, J_{X\times\mathbb{CP}^{d}}\circ du\circ j )=0
    \end{equation}
   \begin{equation}
       u(\partial\Sigma)\subseteq E_{\mathbb{R}}
   \end{equation}
  \begin{equation}
      s\in H^{0}(u^{*}E, u^{*}E_{\mathbb{R}})
  \end{equation}
\end{corollary}
\begin{proof}
   On each spherical component of $\Sigma$, the above statement is exactly Corollary \ref{shearedCurves}. On the other hand, over each disk componenet of $\Sigma$, the same proof of Corollary \ref{shearedCurves} applies after noting that any orientable real vectorbundle over $S^1$ is trivial and hence, $u^*E_\R$ is indeed a totally-real sub-bundle of $u^*E|_{\partial\Sigma}\rightarrow\partial\Sigma$.
\end{proof}
Using the same notation as in the above Corollary, if in particular $u=v\times\phi:(\Sigma, \partial\Sigma)\rightarrow(X\times\CP^d, L\times\RP^d)$ then using Equation \ref{Equ4} we can write Equations \ref{Equ5} and \ref{Equ6} as
\begin{equation}
        \Bar{\partial}_{J}v+\langle\eta\rangle\circ d\phi=0.
    \end{equation}
\subsubsection{Equivariant Hormander's Theorem}
In order to apply the arguments in \cite{abouzaid2021complex} to the proof of Theorem \ref{ChartsDisks}, we need to extend their $\textit{peak-section}$ construction to the equivariant case.\\
 Fix $p, q$ positive integers and let $(Y, J_{Y}, \omega_{Y})$ be a closed K$\ddot a$hler manifold and $(E, \nabla, h)\rightarrow Y$ be a holomorphic Hermitian vectorbundle of curvature form $iR\in C^{\infty}(\bigwedge^{2}T^{*}Y\otimes Herm(E, E))$ where $Herm(E, E)$ is the bundle of Hermitian endomorphisms on $E$. Denote by $\bigwedge_{\omega}:\bigwedge^{p, q}T^{*}Y\otimes Herm(E, E)\rightarrow\bigwedge^{p-1, q-1}T^{*}Y\otimes Herm(E, E)$ to be the metric adjoint of the operator $\omega\wedge .$. Moreover, we assume that there is a finite group $G$ that acts holomorphically and isometrically on $Y$ and $E$ (and hence preserving $iR$ and $\bigwedge_{\omega}$).\\
 \begin{proposition}
     
     If the commutator of the above two operators, $[iR, \bigwedge_{\omega}]:\bigwedge^{p, q}T^{*}Y\otimes Herm(E, E)\rightarrow\bigwedge^{p, q}T^{*}Y\otimes Herm(E, E)$ is positive-definite on each fibre and if we denote by $c:=||A||$ the operator norm of $A$ then, for any $g\in L^{2}(\bigwedge^{p, q}T^{*}Y\otimes E)$ such that $\Bar{\partial}g=0$ and is $G$-invariant, there exists $f\in L^{2}(\bigwedge^{p, q}T^{*}Y\otimes E)$ which is $G$-invariant satisfying $\Bar{\partial}f=g$ and $||f||_{2}\leq c||g||_{2}$.
 \end{proposition}
\begin{proof}
    We give three proofs.
    \begin{enumerate}
        \item (Averaging Argument) By the standard Hormander's Theorem, there exists $f\in L^{2}(\bigwedge^{p, q}T^{*}Y\otimes E)$ such that, $\Bar{\partial}f=g$ and $||f||_{2}\leq c||g||_{2}$. Now define $\Tilde{f}:=\frac{1}{|G|}\Sigma_{a\in G}a.f$.
        \item (Dimension Count) We can argue as in the proof of Lemma \ref{DimensionLemma}. Using elliptic regularity, we can interpret the result of the standard Hormander's theorem as $H^{p, q}(Y, E)=0$. The $G$-action induces an action on the twisted Doulbeaut complex $(\mathcal{E}^{p, *}(Y, V), \Bar{\partial})$ where the invariant part $(\mathcal{E}^{p, *}_{G}(Y, E), \Bar{\partial})$ is also a chain complex as $G$ acts holomorphically. Thus $H^{p, q}_{G}(Y, E)\leq H^{p, q}(Y, E)$ is a real subspace and hence is zero-dimensional by the standard Hormander's Theorem.
        \item  We extend the proof of the standard Hormander's Theorem to the equivariant case. \\
        Using the metric on $TY$ and $E$ we define an $L^{2}$-inner product $(. , .)$ on smooth sections $C^{\infty}(\bigwedge^{p, q}T^{*}Y\otimes E)$ and using $\omega_{Y}$ we define $\Bar{\partial}^{*}$ the adjoint of $\Bar{\partial}$. For $h\in C^{\infty}(\bigwedge^{p, q}T^{*}Y\otimes E)$, we have $||h||_{2}^{2}\leq c(||\Bar{\partial}h||^{2}_{2}+||\Bar{\partial}^{*}h||^{2}_{2})$ and thus the bilinear pairing $(h_{1}, h_{2})_{\mathcal{H}}:=(\Bar{\partial}h_{1}, \Bar{\partial}h_{2})+(\Bar{\partial}^{*}h_{1}, \Bar{\partial}^{*}h_{2})$ defines a metric. Let $\mathcal{H}$ be the Hilbert space completion of  $C^{\infty}(\bigwedge^{p, q}T^{*}Y\otimes E)$ with respect to $||.||_{\mathcal{H}}$. Consider $\psi_{g}(h):=(h, g)_{\mathcal{H}}$ and note that $\psi_{g}\in\mathcal{H}^{*}$ as by Cauchy-Schwartz, we have, $|\psi_{g}(h)|^{2}\leq c||g||^{2}_{\mathcal{H}}||h||_{\mathcal{H}}^{2}=c||g||_{2}^{2}||h||_{\mathcal{H}}^{2}$. Using Riez-Representation Theorem, it follows that there exists $v\in\mathcal{H}$ such that, $||v||_{\mathcal{H}}^{2}\leq c||g||_{2}^{2}$ and $(v, h)_{\mathcal{H}}=\psi_{g}(h)$, for all $h\in\mathcal{H}$. Noting that $g$ is $G$-invariant, it follows that, $a. \psi_{g}(h)=\psi_{g}(a.h)$, for all $ a\in G$ and $h\in\mathcal{H}$. Thus $v$ is $G$-invariant by the non-degeneracy of $(. , .)_{\mathcal{H}}$. Now setting $f:=\Bar{\partial}^{*}v$ solves the equation $\Bar{\partial}f=g$ with $f$ being $G$-invariant.
    \end{enumerate}
\end{proof}
\subsubsection{Construction of the Thickening}
\begin{definition} [Topological Thickening] \label{Ttop}
    We set $\mathcal{T}^{top}:=\{((\Sigma, \partial\Sigma; \underline{z}, \underline{w}), u, F, \eta)\}$ where $((\Sigma, \partial\Sigma;\underline{z},\underline{w}), u, F)$ is a marked framed disk, as in Definition \ref{FramedDisk}, and $\eta\in E$, where $E$ is as in Definition \ref{preOb}. We require $u:(\Sigma,\partial\Sigma)\rightarrow (X, L)$ to satisfy:
        \begin{equation} \label{pertDbar}
        \Bar{\partial}_{J}u+\langle\eta\rangle\circ d\phi_{F}=0
    \end{equation}
    where $\phi_F:(\Sigma, \partial\Sigma)\rightarrow(\CP^d, \RP^d)$ is as in Proposition \ref{FramingMap}.
\end{definition}
We give $\T^{top}$ the topology induced from the Gromov-Hausdorff topology on the closure of the images of $u\times\phi_{F}$ in $X\times\mathbb{CP}^{d}$ and the vectorspace topology of $E$.
\begin{lemma}\label{topThickening}
     The topology on $\mathcal{T}^{top}$ agress with the Gromov's convergence topology.
\end{lemma}
\begin{proof}
    As we are only dealing with the genus-zero case, the marked domain $(\Sigma, \partial\Sigma;\underline{z},\underline{w})$ can be represented by a tree. Now after applying the structure theorems of \cite{kwon2000structure} or \cite{lazzarini2011relative}, the result follows from the argument of Lemma 6.14 in \cite{abouzaid2021complex}.
\end{proof}
\begin{remark}
    We note that, in the case of higher genus bordered curves, the argument of Lemma 6.14 in \cite{abouzaid2021complex} cannot be applied as we may contract a loop, representing a constant genus-one curve which cannot be factored by simple curves as in \cite{mcduff2012j}.
\end{remark}
We now proceed to give $\T^{top}$ a (stable) smooth structure. We will use $\textit{Smoothing Theory}$ to provide us with a smooth structure, after stabilization. We refer the reader to Appendix A for the relevant definitions and results needed in what follows.\\

We start by recalling some notation. For $((\Sigma, \partial\Sigma; \underline{z}, \underline{w}), u, F, \eta)\in\T^{top}$, we have an induced smooth map $u\times\phi_F:(\Sigma, \partial\Sigma)\rightarrow (X\times\CP^d, L\times\RP^d)$. Using the discussion in the last paragraph of $\textit{Gromov's Graph Trick}$ subsection, we have, after composing with the zero-section of $E\rightarrow X\times\CP^d$ and abuse of notation, $u\times\phi_F:(\Sigma, \partial\Sigma)\rightarrow (E,E_\R)$ is a $J_\Psi$-holomorphic map with boundary conditions in $E_\R$. We denote by $\mathcal{M}(E,E_\R; J_\Psi):=\M_{k+1, l}(E, E_\R; J_{\Psi}, \Tilde{\beta})$ the moduli space of genus-zero bordered $J_\Psi$-holomorphic maps with $(k+1)$-boundary marked points and $l$-interior marked points, with boundary conditions in $E_\R$ and representing the class $\Tilde{\beta}:=\beta\oplus\alpha$ after identifying $H_2(E, E_\R;\Z )\cong H_2(X\times\CP^d, L\times\RP^d; \Z)$.
\begin{proposition} \label{ChoosingN}
    There exists a positive integer $N_0$ so that for all $N\geq N_0$, the linearization $D\Bar{\partial}_{J_\Psi}$ of $\Bar{\partial}_{J_\Psi}$ is surjective at any point of $\T^{top}\hookrightarrow\mathcal{M}(E, E_\R; J_\Psi)$. 
\end{proposition}
\begin{proof}
    We first observe that after choosing $N>>1$ large enough, we can assume that $H^1(X\times\CP^d; E)=\{0\}$. Thus after using Lemma \ref{topThickening}, it follows that the map $((\Sigma, \partial\Sigma; \underline{z}, \underline{w}), u, F, \eta)\in\T^{top}\mapsto ((\Sigma, \partial\Sigma; \underline{z}, \underline{w}), u\times\phi_F)\in\mathcal{M}(E, E_\R; J_\Psi)$ is a well-defined topological embedding. Now we follow the proof of Proposition 6.26 in \cite{abouzaid2021complex}.\\

     Let $(\Sigma, \partial\Sigma; \underline{z}, \underline{w}), u, F, \eta)\in\mathcal{T}^{top}$ and assume for a contradiction that, $\operatorname{coker}(D_{J_{\Psi}})\neq\{0\}$. We identify $\operatorname{coker}(D_{J_{\Psi}})\cong \ker(D_{J_{\Psi}}^{*})$ where $D_{J_{\Psi}}^{*}$ is the $L^{2}$-adjoint of $D_{J_{\Psi}}$ and let $e\in \ker(D_{J_{\Psi}}^{*})$ be a non-zero element and assume that the support of $e$ is on a disk-component (as otherwise we argue as in \cite{abouzaid2021complex} Proposition 6.26). Notice that, if we assume that $e$ is of $W^{k, p}$-regularity, the equation $D^{*}_{J_\Psi}e=0$ together with elliptic regularity implies that $e$ is smooth and thus we assume everything that follows is smooth. \\
     Noting that, $\phi_{F}$ is a regular curve in $(\mathbb{CP}^{d}, \mathbb{RP}^{d})$, that is $\operatorname{coker}(D_{\phi_F}\bar \partial_{J_{std}})=\{0\}$, it follows that $e\in C^{\infty}(u^{*}TX, u^{*}TL)\hookrightarrow C^{\infty}((u\times\phi_{F})^{*}(TX\times T\mathbb{CP}^{d}), (u\times\phi_{F})^{*}(TL\times T\mathbb{RP}^{d}))$. Moreover, notice that the linearization of $\Bar{\partial}_{J_{\Psi}}$ is
    \begin{equation} \label{LinearSheared}
        D_{J_{\Psi}}=D_{u}\Bar{\partial}_{J}+ \langle ,\rangle\circ d \phi_{F}
    \end{equation} 
    as $\langle ,\rangle\circ d \phi_{F}$ is a linear operator.
    \begin{enumerate}
        \item Assume that $u$ is non-constant and has a smooth, interior, non-branched point $z\in\Sigma$.
    Now denote by $\Tilde{e}\in C^{\infty}(\Tilde{\Sigma}, u^{*}TX)$ the double of $e$ and fix a smooth vector $w\in T_{z}\Sigma$ and let $f$ be a complex anti-linear map $T_{\phi_{F}(z)}\mathbb{CP}^{d}\rightarrow T_{u(z)}X$ satisfying
    \begin{equation}\label{non-zero}
        \langle f(d\phi_{F}(j(w))), \Tilde{e}(z)\rangle\neq 0
    \end{equation}
    where $j:T\Sigma\rightarrow T\Sigma$ is the complex structure on $\Sigma$. Using Equivariant Hormander's Theorem and for $N>>1$ large enough, we can find holomorphic section $s_{N}$ of $\overline{\operatorname{Hom}}(\Tilde{\phi}_{F}^{*}T\CP^{d}, u^{*}TX)\otimes_\C\tilde{\phi}_{F}^{*}\mathcal{O}(N)\rightarrow\Tilde{\Sigma}$ which are $\Z/2$-equivariant with respect to conjugation. In addition, using the proof in \cite{abouzaid2021complex} Lemma 6.24, we can find holomorphic sections of $\Tilde{s}_{N}$ of $\tilde\phi_{F}^{*}\mathcal{O}(N)\rightarrow\Tilde{\Sigma}$ that converges to $\delta_{z}$. Using Lemma \ref{LiftingInvolutions}, we can lift the conjugation action on $\Tilde{\Sigma}$ to a $\Z/2$-action on $\tilde\phi_{F}^{*}\mathcal{O}(N)$ and hence using Equivariant Hormander's Theorem we deduce that for $N>>1$ and upon restricting the above sections to $\Sigma\subset\Tilde{\Sigma}$ we get $\langle s_{N}, \Tilde{s}_{N}\rangle$ converges to the Dirac delta section $\delta_{f}\in\bigwedge^{0, 1}\Sigma\otimes u^{*}TX$. Now the result follows from the proof of Proposition 6.26 in \cite{abouzaid2021complex}. Indeed, using Equation \ref{non-zero}, it follows that for $N$ large enough, we have,
    \begin{equation}
         \langle\langle s_{N}\otimes\Tilde{s}_{N}\rangle\circ d\tilde\phi_{F}(j(\zeta)), \tilde e\rangle=
        \langle\langle s_{N}\otimes\Tilde{s}_{N}\rangle\circ d\phi_{F}(j(\zeta)), e\rangle\neq 0
    \end{equation}
    where $\zeta$ is a smooth vector field on $\Sigma$ so that $\zeta(z)=w$. Contradicting the assumption that, $e\in \ker(D_{J_{\Psi}}^{*})$.
    \item  Assume that $u$ is constant. We note that, $u^{*}TX$ is a trivial holomorphic bundle with a trivial totally real sub-bundle $u^*_{\partial\Sigma}TL$ and thus, upon doubling $u$ to $\Tilde{u}$ (which we can do as $u$ is constant), we get $\bigwedge^{0, 1}\Sigma\otimes u^{*}TX\rightarrow\bigwedge^{0, 1}\Tilde{\Sigma}\otimes\Tilde{u}^{*}TX$ which is injective. Now the result follows by elliptic regularity and Riemann-Roch.
    \item For the general case, we apply the above arguments after using the structure theorems of \cite{kwon2000structure} or \cite{lazzarini2011relative}. 
    \end{enumerate}
    To finish the proof, we "shrink" $\T^{top}$, and abuse notation and still denote it by $\T^{top}$, so that $\T^{top}\hookrightarrow\mathcal{M}(E, E_\R; J_\Psi)$ has a compact closure. Namely, we consider $((\Sigma, \partial\Sigma; \underline{z}, \underline{w}), u, F, \eta))\in\T^{top}$ where using the Hermitian metric on $E$, we have $|\eta|< C$ for some fixed real number $C>0$. Now as the surjectivity of $D\Bar{\partial}_{J_\Psi}$ is an open condition, the result follows as desired.
\end{proof}

To avoid notational clutter, denote by $\F:=\F_{k+1, l}(\alpha)$ the space of framings as in Definition \ref{spaceFramings}. Using Lemma \ref{topThickening}, it follows that the projection map 
\begin{align*} \label{projection}
    \T^{top} \ & \to\F \\
   \ ((\Sigma, \partial\Sigma; \underline{z}, \underline{w}), u, F, \eta) \ & \to ((\Sigma, \partial\Sigma; \underline{z}, \underline{w}), \phi_F)
\end{align*}
is (sequentially) continuous. 

\begin{corollary} \label{IsoOrbiSpace}
$\mathcal{T}^{top}$ is a topological $O(d+1)$-manifold with a fibre-wise smooth structure on $\mathcal{T}^{top}\rightarrow\mathcal{F}$. In addition, the stabilizer group of any point $((\Sigma, \partial\Sigma; \underline{z}, \underline{w}), u, F, \eta)\in\T^{top}$ is a finite group isomorphic to $\operatorname{Aut}(u)$, the automorphism group of the underlying $J$-holomorphic curve.    
\end{corollary}
\begin{proof}
    The first statement is precisely Proposition \ref{ChoosingN}, where the local Euclidean charts of $\T^{top}$ are given by the open neighborhoods of $\mathcal{M}(E, E_\R; J_\Psi)$ where $D\Bar{\partial}_{J_\Psi}$ is surjective. As for the second statement, we argue as in \cite{abouzaid2021complex} Lemma 6.4 on the dual graph of the curve $((\Sigma, \partial\Sigma;\underline{z},\underline{w}), u)$ after using the structure theorems of \cite{kwon2000structure} or \cite{lazzarini2011relative} or as in the proof of Theorem $4.30$ in \cite{abouzaid2024gromov}.
\end{proof}
\begin{proposition} \label{C1-loc structure}
    There exists $V$ a real $O(d+1)$-representation so that $\T^{top}\times V$ can be given an $O(d+1)$-equivariant smooth structure. Moreover, such smooth strucuture is unique up to stable $O(d+1)$-isotopy of smooth structures (c.f. Definition \ref{stableIsoDefn}).
\end{proposition}
\begin{proof}
     This is an application of Lashof's Theorem \ref{lashof} together with section 4.1 of \cite{abouzaid2021complex}, provided that the projection map $\pi:\mathcal{T}^{top}\rightarrow\mathcal{F}$ is an $O(d+1)$-equivariant $C^{1}_{loc}$-topological submersion as in Definition \ref{TopSubmersion}. Which in return we obtain using the $C^1$-estimates in \cite{fukaya2021construction} or \cite{fukaya2024corrigendum}.
     To avoid notational clutter, we write the above projection map $\pi$ as $(\Sigma,\partial\Sigma; u, F, \eta )\mapsto \phi_F$, with the understanding that $(\Sigma, \partial\Sigma)$ is a marked genus-zero bordered Riemann surface. We also denote by $\M:=\M_{k+1, l+l'}$ the Deligne-Mumford moduli space of genus-zero bordered Riemann surfaces with $(k+1)$-boundary marked points and $(l+l')$-interior marked points. Here $l'\geq \dim_\R\F +3$ is a fixed positive integer. \\
     
     For $\phi_F\in\F$ we can identify a neighborhood of $\phi_F$ by a neighborhood of $(\Sigma,\partial\Sigma; w_1,\dots,w_l')\in\M$ as follows: \\
     Fix $L_1, \dots, L_{l'}$ linear hypersurfaces of $\CP^d$ such that they are transverse to $\phi_F$. The existence of such linear hypersurfaces follows from Proposition $6.5$ of \cite{abouzaid2021complex} after noting that $\phi_F$ is realized as the pull-back by the doubling map of an element $\psi\in\overline{\mathcal M}_{k+1, 2l}(\CP^d;J_{std}, \alpha+\bar\alpha)$ as per \cite{zernik2017moduli}. Alternatively, using the fact that elements of $\F$ are automorphism-free, we can construct such linear hypersurfaces as follows. Fix $U_{1}, \dots, U_{l'}\subset\mathbb{CP}^{d}$ open subsets whose complement is of codimension 1, intersecting the image of $\phi_F$ transversally at $l'$-distinct interior points. The existence of such open sets, is done upon noting that $\phi_F$ is a simple curve (as in \cite{mcduff2012j}) and thus we can apply Frobenious Theorem locally on $L_{1}, \dots, L_{l'}\leq T_{\phi(w_{i})}\mathbb{CP}^{d}$ where $w_{i}\in\Sigma$ is an injective point (as in \cite{mcduff2012j}), and thus in particular $\zeta\mapsto d_{w_{i}}\phi_F(\zeta)$ is injective and $L_{i}$ is chosen to be transverse to this map. Denote by $\mathcal{F}(L_{*})\subset\mathcal{F}$ the open subset of all elements transverse to $U_{1}, \dots, U_{l'}$. Then by Proposition $6.5$ of \cite{abouzaid2021complex}, the map $\mathcal{F}(L_{*})\rightarrow\overline{\mathcal{M}}$ sending $\phi\in\F(L_*)$ to its domain with $w_{1}, \dots, w_{l'}$ as extra interior marked points, is the desired local diffeomorphism.\\
     Now denote by $\Pi:\pi^{-1}(\F(L_*))\rightarrow\M$ the composition of $\pi$ with the above map. For a fixed $(\Sigma,\partial\Sigma; u,  F, \eta)\in\T^{top}$, let $U\subseteq\overline{\mathcal{M}}$ be an open set so that $\mathcal{F}(L_{*})\rightarrow\overline{\mathcal{M}}$ is a diffeomorphism onto $U$, and let $K\subset\Pi^{-1}(U)\subseteq\mathcal{T}$ be an open set. Denote by $g:U\rightarrow g(U)=\F(L_*)\subseteq\F$ be the local inverse of $\F(L_*)\rightarrow\M$. In order use to use the estimates found in \cite{fukaya2021construction}, we observe that for $N>>1$ large enough, the linear operator $\langle.\rangle\circ d\phi_F:E\rightarrow\Gamma(\Sigma; \Lambda^{0, 1}\Sigma\otimes u^*TX)$ contains the subspaces found in Definition 9.2 of \cite{fukaya2021construction}. Indeed, for increasing positive integers $N$, the subspaces $\{\operatorname{Im}(\langle.\rangle\circ d\phi_F)\leq \}_N$ will have a direct limit equal to an open and dense subset of $\Gamma(\Sigma; \Lambda^{0, 1}\Sigma\otimes u^*TX)$. Now denote by $G$ the composition of $g$ with the map $\textit{Glue}$ as constructed in Theorem 3.13, 8.16 of \cite{fukaya2016exponential}. After possibly taking smaller $U$ and $K$, the following diagram commutes:
     \[\begin{tikzcd}
	{U\times K} && {\pi^{-1}(\mathcal{F}(L_*))} \\
	\\
	U && {\overline{\mathcal{M}}}
	\arrow["G", from=1-1, to=1-3]
	\arrow["{pr_{U}}"', from=1-1, to=3-1]
	\arrow["\Pi", from=1-3, to=3-3]
	\arrow[hook, from=3-1, to=3-3]
\end{tikzcd}\]
where the restriction $G|_{\{a\}\times K}$ is obtained by Newton-Iteration as in Section 5 of \cite{fukaya2016exponential}. This implies that $G|_{\{a\}\times K}:K\rightarrow\pi^{-1}(\F(L_*))$ is a smooth map. Now using Theorem 6.4 of \cite{fukaya2016exponential} and noting that the $\textit{pre-glued}$ curve, as in Equation 5.4 of \cite{fukaya2016exponential} depends smoothly on the gluing parameters. On the other hand, by the estimates $(9.13)$ of \cite{fukaya1710construction} it follows that the comparison of two different gluing parameters are locally $C^1$-smooth for a fixed domain. Therefore, $\pi$ is a $C^1_{loc}$-topological submersion. Finally, after using Proposition \ref{microBundleLift}, we apply Lashof's Theorem \ref{lashof} to get our desired result.
\end{proof}
\subsubsection{Proof of Theorem \ref{ChartsDisks}}
\begin{proof}
    Let $G=O(d+1)$ and set $\T:=\T^{top}\times V$ where $V$ is as in Proposition \ref{C1-loc structure}. Denote by $\mathcal{H}_d$ the vectorspace of $(d+1)\times(d+1)$-symmetric matrices and define $\E:=E\oplus\mathcal{H}_d\oplus V$, where $E$ is as in Corollary \ref{ObwithVanishing}. Finally, we define $s:\T\rightarrow\E$ by $((\Sigma, \partial\Sigma; \underline{z}, \underline{w}), u, F, \eta; v)\mapsto (\eta, \exp^{-1}(H_F), v)$ where $H_F$ is as in the third item of Definition \ref{FramedDisk} and $\exp:\mathcal H_d\rightarrow\mathcal{H}^+_d$ is the matrix exponential map and $\mathcal{H}_d^+$ denotes the convex cone of all symmetric positive-definite matrices.\\
    For $[(\Sigma, \partial\Sigma; \underline z, \underline w), u, I_d, 0; 0]\in s^{-1}(0)/G$, where $I_d$ is the $(d+1)\times(d+1)$ identity matrix, the map $[(\Sigma, \partial\Sigma; \underline z, \underline w), u, I_d, 0; 0]\mapsto [(\Sigma, \partial\Sigma; \underline z, \underline w), u]\in\M_{k+1, l}(X, L; J, \beta)$ serves as a homeomorphism after using the second statement of Corollary \ref{IsoOrbiSpace}.
    Therefore, $(G, \T, \E, s)$ is a Global Kuranishi chart of $\M_{k+1, l}(X, L; J, \beta)$ as desired.\\ 
    
    Now we show that, possibly after stabilization, such Global Kuranishi chart admits a normally complex structure. Indeed, let $((\Sigma, \partial\Sigma; \underline{z}, \underline{w}), u, F, \eta)\in\mathcal{T}^{top}$ and note that, the automorphism group of each disk domain with at least one boundary marked point is torsion-free and hence disk-components have trivial isotropy groups. With this in mind and by abuse of notation, assume that $\Sigma$ is a closed Riemann surface and argue as follows. Note that, the term $\langle ., .\rangle\circ d\phi_{F}$ in Equation \ref{LinearSheared} is a compact operator. Moreover, after fixing a $J$-linear connection on $X$, we can write the operator $D_u\Bar{\partial}_J=O+K$ where $O$ is a complex linear operator and $K$ is a compact operator. With Corollary \ref{IsoOrbiSpace} in mind and using the linear homotopy, $t\rightarrow O+ t(K+\langle,\rangle\circ d \phi_{F})$ of Fredholm operators, we get an isomorphism $T^{vt}\T^{top}\oplus\operatorname{coker}(O)\cong\ker(O)$, as the operator in Equation \ref{LinearSheared} is surjective. Finally, we fix an $O(d+1)$-equivariant connection and bundle metric on $T\mathcal{T}$ where the fibres of $T^{vt}\T^{top}$ are totally-geodesic, and a non-zero vector $v\in T\mathcal{F}$. Taking parallel transport of $\operatorname{coker}(O)$ along the horizontal lift of $v$, gives the desired claim; after further stabilization by $iV$ where $V$ is as in Proposition \ref{C1-loc structure}.\\
Lastly, we further stabilize by $ev^{*}(TX^l\oplus TL^{k+1})$ and argue as in \cite{abouzaid2021complex} Lemma 4.5. 
\end{proof}
\subsection{Using \cite{abouzaid2024gromov}}
In the above proof, we have used the \textit{framings}, namely the associated holomorphic map $\phi:(\Sigma, \partial\Sigma)\rightarrow(\CP^d, \RP^d)$ to perturb the $\Bar{\partial}$-equation itself and make the linearization of the perturbed equation surjective. Following \cite{abouzaid2024gromov}, we can also give a construction that proves Theorem \ref{ChartsDisks}, using the space $\F\equiv\mathcal{F}_{k+1, l}(\alpha)$. 
We start by recalling the definition of a $\textbf{finite-dimensional approximation scheme}$ as in \cite{abouzaid2024gromov}.
\begin{definition}[Finite-Dimensional Approximation Scheme]\label{finite dimensional}
    Let $G$ be a compact Lie group and $\pi:V\rightarrow B$ be a smooth $G$-equivariant vectorbundle. A finite dimensional approximation scheme $(V_{\mu}, \lambda_\mu)_{\mu\in\mathbb{N}}$ for $C^\infty_c(V)$, the space of compactly-supported smooth sections of $V$, is a sequence of finite-dimensional $G$-representations $(V_{\mu})$ together with a sequence of $G$-equivariant linear maps $\lambda_\mu:V_{\mu}\rightarrow C^\infty_c(V)$, for all $\mu\in\mathbb{N}$ satisfying the following conditions:
    \begin{enumerate}
        \item $V_\mu\leq V_{\mu+1}$ is a sub-representation, for every $\mu\in\mathbb{N}$.
        \item $\lambda_{\mu}|V_{\mu-1}=\lambda_{\mu-1}$, for every $\mu\in\mathbb{N}$.
        \item $\bigcup_{\mu}\lambda_{\mu}(V_{\mu})\subset C^{\infty}(V)$ is dense with respect to the $C^\infty_{loc}$-topology.
    \end{enumerate}
\end{definition}
\begin{lemma}\cite{abouzaid2024gromov}\label{finiteAppScheme}
    Finite dimensional approximation schemes exists.
\end{lemma}
Denote by $\mathcal{C}\rightarrow\F$ to be the \textbf{universal curve} over $\F$. Following the notation in Proposition \ref{spaceFramings}, $\mathcal{C}$ can be identified by $\F_{k+1, l+1}(\alpha)$ where $\mathcal{C}\rightarrow\F$ is the forgetful map. In particular, the fibre of $\mathcal{C}$ over a point $[(\Sigma,\partial\Sigma; \underline{z}, \underline{w}); \phi]\in\F$ comes with an identification with $(\Sigma,\partial\Sigma)$, the domain of $\phi$. Moreover, as in Proposition \ref{FramingMap}, $\mathcal{C}$ is a smooth $O(d+1)$-manifold. 
\begin{definition}[Fibre-wise Metric]
    A fibre-wise metric on $\mathcal{C}$ is a continuous map $\mu:\mathcal{C}\times_\F\mathcal{C}\rightarrow\R_{\geq 0}$ such that the restriction of $\mu$ to each irreducible component of the fibre over a point in $\F$, is a metric induced from a smooth Riemannian metric on such irreducible component.
\end{definition}
\begin{definition}[Consistent Domain Metric]
    A consistent domain metric on $\F$ is an $O(d+1)$-invariant fibre-wise metric on $\mathcal{C}$.
\end{definition}
\begin{lemma}\label{consMetric}
    $\mathcal{F}$ admits a consistent domain metric.
\end{lemma}
\begin{proof}
    This is a simpler version of Lemma 4.17 in \cite{abouzaid2024gromov} as $O(d+1)$ is a compact Lie group.\\ Indeed, upon fixing an $O(d+1)$-invariant cover of $\mathcal{F}$ formed of trivializing open sets of $\mathcal{C}\rightarrow\mathcal{F}$ and upon fixing a Haar measure on $O(d+1)$, we get the result by an averaging argument and partition of unity on such a cover.
\end{proof}
Now note that $\mathcal C^*\subset\mathcal C$, the complement of nodes and marked points is an $O(d+1)$-smooth manifold. With this in mind and after denoting $O(d+1)$ by $G$, we consider the following vectorbundle.
\begin{definition}\label{0,1 forms}
Let $V:=\Lambda_{\mathcal{C}^*/\mathcal{F}}^{0,1}\otimes_{\mathbb{C}}TX\longrightarrow \mathcal{C}^*\times X$
be the $G$-equivariant vector bundle whose fiber over a point $(c,x)\in \mathcal{C}^*\times X$ is the space of all complex anti-linear maps $T_c\bigl(\mathcal{C}^*|_{\pi(c)}\bigr)\longrightarrow T_xX$.
\end{definition}
Using Lemma \ref{finiteAppScheme}, we fix $(V_\mu, \lambda_\mu)_{\mu\in\mathbb{N}}$ a finite-dimensional approximation scheme of $V$.
\begin{definition}
    We define the pre-thickened moduli space $\T^{pre}:=\mathcal{T}^{pre}(\beta, (V_{\mu}, \lambda_\mu))$ to be the set of tuples $((\Sigma,\partial\Sigma, \underline{z},\underline{w}), u, ,F, e)$ where $((\Sigma,\partial\Sigma, \underline{z},\underline{w}), u, ,F)$ is a framed marked disk, as in Definition \ref{FramedDisk}, representing the class $\beta$ and $e\in V_{\mu}$ for some $\mu\in\mathbb{N}$. We require $u$ to satisfy
    \begin{equation}\label{pertEquAMSII}
        \Bar{\partial}_{J}u|_{\mathcal{C}^*_{\phi_F}}+(\lambda_{\mu}(e))\circ\Gamma_{u}=0
    \end{equation}
    where $\mathcal{C}^*_{\phi_F}$ is the complement of all nodes in the fibre of $\pi:\mathcal{C}\rightarrow\mathcal{F}$ over $\phi_F$, where $\phi_F$ is given as in Proposition \ref{FramingMap}, and $\Gamma_u:\mathcal{C}_{\phi_F}\rightarrow \mathcal{C}\times X$ is the graph map.
\end{definition}
We give $\mathcal{T}^{pre}$ the Gromov-Hausdorff topology on the graphs given by the closure of the image of $\Gamma_{u}$ and the vectorspace topology on $\varinjlim V_{\mu}$.
\begin{remark}
    In order to avoid notational clutter and as we are keeping the identification $\iota:\mathcal{C}_{\phi_F}\cong(\Sigma, \partial\Sigma)$ implicit, we write $\Gamma_u(z)=(z, u(z))$ instead of $\Gamma_u(z)=(z, u(\iota(z)))$.
\end{remark}
\begin{lemma}\label{choosingMu}
    There exists a $\mu\in\mathbb{N}$ so that the linearization of Equation \ref{pertEquAMSII} is surjective at every point of $\mathcal{T}^{pre}$.
\end{lemma}
\begin{proof}
    This follows from Gromov's compactness and openness condition of having the linearization surjective. Indeed, noting that the cokernel of the linearization of the unperturbed $\Bar{\partial}_J$-equation is a finite-dimensional subspace of $C^{\infty}_c(V)$, the result follows from the third item of Definition \ref{finite dimensional}. 
\end{proof}
From now on, we fix $\mu\in\mathbb{N}$ as in Lemma \ref{choosingMu} and consider the following topological space.
\begin{definition}[Topological Thickening]
    We define the topological thickening $\T^{top}$ to be the set of all tuples $((\Sigma,\partial\Sigma, \underline{z},\underline{w}), u, ,F, e)\in\T^{pre}$ where $e\in V_\mu$ for the fixed $\mu$ as above. 
\end{definition}
Similary, $\T^{top}$ is given the Gromov-Hausdorff topology on the graphs given by the closure of the image of $\Gamma_{u}$ and the vectorspace topology on $V_{\mu}$. Moreover, by the same proof of Lemma \ref{topThickening}, the topology given on $\T^{top}$ agrees with the Gromov's convergence topology. Now consider the projection map $\T^{top}\rightarrow\F$ also given by $((\Sigma,\partial\Sigma; \underline{z},\underline{w}), u,F,e)\mapsto\phi_F$ and observe that using Lemma \ref{choosingMu} and the proof of Corollary \ref{IsoOrbiSpace}, we get the following Corollary.

\begin{corollary} \label{IsoOrbiSpaceII}
$\mathcal{T}^{top}$ is a topological $O(d+1)$-manifold with a fibre-wise smooth structure on $\mathcal{T}^{top}\rightarrow\mathcal{F}$. In addition, the stabilizer group of any point $((\Sigma, \partial\Sigma; \underline{z}, \underline{w}), u, F, e)\in\T^{top}$ is a finite group isomorphic to $\operatorname{Aut}(u)$, the automorphism group of the underlying $J$-holomorphic map.    
\end{corollary}  
Furthermore, using the same line of proof as in Proposition \ref{C1-loc structure} and after noting that from the third item of Definition \ref{finite dimensional}, it follows that for $\mu$ large enough, we can assume that $\lambda_\mu(V_\mu)$ contains all the subspaces found in Definition 9.2 of \cite{fukaya2021construction}, we get that $\T^{top}$ admits a $G$-smoothing. For completeness we state the following Proposition.
\begin{proposition} \label{Smoothing}
 There exists $R$ a real $O(d+1)$-representation so that $\T^{top}\times R$ can be given an $O(d+1)$-equivariant smooth structure. Moreover, such smooth strucuture is unique up to stable $O(d+1)$-isotopy of smooth structures.    
\end{proposition}
\begin{definition}[Thickening]
    We define the thickening to be $\T:=\T^{top}\times R$ with its $O(d+1)$-equivariant smooth structure, where $R$ is as in Proposition \ref{Smoothing}.
\end{definition}
\begin{definition}
    Fix a consistent domain metric on $\F$ and let $((\Sigma, \partial\Sigma;\underline{z}, \underline{w}), u, F, e; v)\in\T$. We define $H_F:=(\int_{\mathcal{C}|_{\phi_F}}\langle f_{i}, f_{j}\rangle u^*\Omega)_{i,j}$, where the pairing $\langle.,.\rangle$ is the $O(d+1)$-equivariant metric on $\mathcal{C}|_{\phi_F}$ induced by Lemma \ref{consMetric}.
\end{definition}
\begin{remark}
    In the above definition, we are abusing notation and writing $u^*\Omega\in\Omega^{1, 1}(\mathcal{C}|_{\phi_F})\cong\Omega^{1,1}(\Sigma)$, keeping the identification $\iota$ implicit.
\end{remark}
\begin{definition}[Obstruction Bundle]
    We define the obstruction bundle $\E$ to be the trivial vectorbundle $V_\mu\oplus\mathcal H_d\oplus V$ over $\T$.
\end{definition}
\begin{definition}[Kuranishi Section]\label{KuranishiSecAMSII}
    We define the Kuranishi map $s:\T\rightarrow\E$ by $((\Sigma, \partial\Sigma;\underline{z}, \underline{w}), u, F, e; v)\mapsto (e, \exp^{-1}(H_F); v)$ where $\exp:\mathcal H_d\rightarrow\mathcal H_d^+$ is the matrix exponential map.
\end{definition}

\begin{proof}[Proof of Theorem \ref{ChartsDisks}]
   Notice that, $s^{-1}(0)=\{(\Sigma, \partial\Sigma;\underline{z}, \underline{w}), u, I_d, 0; 0)\}$ where $I_d$ is the $(d+1)\times(d+1)$ identity matrix, and $u$ represents the class $\beta$ and  satisfies $\Bar{\partial}_Ju=0$ . Thus, the map $((\Sigma, \partial\Sigma;\underline{z}, \underline{w}), u, I_d, 0; 0)\mapsto((\Sigma, \partial\Sigma;\underline{z}, \underline{w}), u)$ is a homeomorphism. On the other hand, by Corollary \ref{IsoOrbiSpaceII}, it follows that the this map induces a homeomorphism between $s^{-1}(0)/G$ and $\M_{k+1, l}(X, L;J, \beta)$ as desired.
\end{proof}
\subsection{Compatibility with Co-Dimension 1 Corners}

As for now, for a given $\beta\in\pi_2(X, L)$ and $k, l$ non-negative integers, we have a constructed a Global Kuranishi chart for $\M_{k, l}(X, L; J, \beta)$, say $(G, \T, \E, s)$. In order to use such charts for further applications, we have to check compatibility with the charts constructed for the co-dimension $1$ corner of $\M_{k+1, l}(X, L; J, \beta)$. To be more precise, if for $i=1,2$; $k_i, l_i$ are non-negative integers, and $\beta_i\in\pi_2(X, L)$ such that 
\begin{enumerate}
    \item $\beta_1+\beta_2=\beta$.
    \item $k_1+k_2=k+1$.
    \item $l_1+l_2=l$.
\end{enumerate}
 Moreover, if $(G_i, \T_i, \E_i, s_i)$ are Global Kuranishi charts of $\M_{k_i+1, l_i}(X, L,\beta_i; J)$ and we denote by $\F_i$ the corresponding framing space as in Definition \ref{spaceFramings}, then we require:
\begin{enumerate} 
    \item \label{GroupEmbedding }$G_1\times G_2\hookrightarrow G$.
    \item \[\begin{tikzcd}
	{\mathcal{T}_1\times_{ev}\mathcal{T}_2} && {\mathcal{T}} \\
	\\
	{\mathcal{F}_{1}\times_{ev}\mathcal{F}_2} && {\mathcal{F}}
	\arrow[hook, from=1-1, to=1-3]
	\arrow[from=1-1, to=3-1]
	\arrow[from=1-3, to=3-3]
	\arrow[hook, from=3-1, to=3-3]
\end{tikzcd}\]
to be a commutative diagram of smooth manifolds.
\item \[\begin{tikzcd}
	{\iota^*(\mathcal{E}_1\oplus\mathcal{E}_2)} && {\mathcal{E}} \\
	\\
	{\mathcal{T}_1\times_{ev}\mathcal{T}_2} && {\mathcal{T}}
	\arrow[hook, from=1-1, to=1-3]
	\arrow[from=1-1, to=3-1]
	\arrow[from=1-3, to=3-3]
	\arrow["{\iota^*(s_1\oplus s_2)}", bend left = 30pt, from=3-1, to=1-1]
	\arrow[hook, from=3-1, to=3-3]
	\arrow["s"', bend right = 30pt, from=3-3, to=1-3]
\end{tikzcd}\]
to be a commutative diagram, where $\iota:\T_1\times_{ev}\T_2\hookrightarrow\T_1\times\T_2$ is the canonical inclusion.
\item The fibre-product Global Kuranishi, $(G_1\times G_2, \T_1\times_{ev}\T_2, \E_1\oplus\E_2,s_1\oplus s_2)$ is equivalent to\\ $(G, \T|_{\pi^{-1}(\iota(\F_1\times_{ev}\T_2))}, \iota^*\E,\iota^*s)$, the restriction of $(G, \T, \E, s)$ to the boundary corresponding to degeneration of the form $\F_1\times_{ev}\F_2$, where $\pi:\T\rightarrow \F$ is the forgetful map.
\item $\iota^*ev=(ev_1, ev_2)$, where $ev$ is the evaluation map on $\T$ and $ev_i$ is the evaluation map on $\T_i$ for $i=1,2$.
\end{enumerate}

We will refer to the above conditions as the $\textbf{compatibility conditions}$.
The first and second items are used to ensure compatibility of the abstract smooth structures on the thickenings. While the third item will be used in the inductive arguments found in section $4$.\\

To this end and for a given smooth map $u:(\Sigma, \partial\Sigma; \underline{z}, \underline{w})\rightarrow (X, L)$ representing the class $\beta$, we abuse terminology and define the following.
\begin{definition}[Framed Marked Disk with Constraint]\label{framedCurvewithCons}
    A framed marked disk, representing the class $\beta\in\pi_2(X, L)$, is a tuple $((\Sigma, \partial\Sigma; \underline{z}, \underline{w}), u, F)$ as in Definition \ref{FramedDisk} with the extra condition, that if we write $F=\{f_0, \dots, f_d\}$ then we require $[f_i(z_0)]_{i=0, \dots, d}=[1:0:\dots:0]$.
\end{definition}
Following the same notation as in Definition \ref{spaceFramings}, we have the following Lemma. 
\begin{lemma} \label{SpaceofFramingConstraint}
    $\F':=\{\phi\in\F:\phi(z_0)=[1:0:\dots:0]\}\subseteq\F$ is a smooth $O(d)$-manifold.
\end{lemma}
\begin{proof}
    The key point of the proof is the fact $PO(d+1)$ acts transitively on $\RP^d$. Moreover, the evaluation map is equivariant with respect to this action. Indeed, denote by $ev:\F\rightarrow\RP^d$ the evaluation map at $z_0$. Observe that $ev$ is a smooth submersion. Indeed, for $\phi\in\F$ we have $d_\phi ev:H^0(\Sigma, \partial\Sigma; \phi^*T\CP^d, \phi|_{\partial\Sigma}^*T\RP^d)\rightarrow T_{\phi(z_0)}\RP^d$ is given by $\zeta\mapsto d_{z_0}\phi(\zeta)$. Now the result follows from the fact, that $T\CP^d$ is globally generated, and that $H^0(\Sigma, \partial\Sigma; \phi^*T\CP^d, \phi|_{\partial\Sigma}^*T\RP^d)\leq H^0(\tilde\Sigma, \tilde\phi^*T\CP^d)$ and $T_{\phi(z_0)}\RP^d\leq T_{\phi(z_0)}\CP^d$ are totally-real subspaces, where $\tilde\phi:\tilde\Sigma\rightarrow\CP^d$ is the double of $\phi$ as in Lemma \ref{DimensionLemma}.\\
    As for the second statement, it follows from the fact that the group of isometries of $\RP^d$ that fixes a point is $O(d)$.
\end{proof}
Denote by $\pi:\mathcal{C}'\rightarrow\F'$ the universal curve over $\F'$ and note that, $\pi$ is a smooth $O(d)$-equivariant map. Similarly, as in Definition \ref{0,1 forms}, let $V\rightarrow\mathcal{C}'\times X$ be the $O(d)$-equivariant vectorbundle of fibre over a point $(c, x)\in\mathcal{C}'\times X$, the space of all anti-holomorphic maps $T_c(\mathcal{C}'|_{\pi(c)})\rightarrow T_x X$. Using Lemma \ref{finiteAppScheme} we can find and fix $(V_\mu ,\lambda_\mu)_{\mu\in\mathbb{N}}$ a finite-dimensional approximation scheme of $V$. Using Lemma \ref{choosingMu}, we fix a $\mu\in\mathbb{N}$ so that for a given framed marked disk $((\Sigma,\partial\Sigma, \underline{z},\underline{w}), u, ,F)$ as in Definition \ref{framedCurvewithCons}, representing the class $\beta\in\pi_2(X, L)$, the linearization, with respect to the $u$-variable, of Equation \ref{pertEquAMSII} is surjective at every point $\phi_F\in\F'$.
\begin{corollary} \label{ChartsDiskswithRightSymmetry}
    The moduli space $\M_{k+1, l}(X, L, \beta; J)$ can be given a Global Kuranishi chart $(O(d), \T, \E, s)$ together with a smooth submersion $ev:\T\rightarrow L$ extending the evaluation map at the zero-th boundary marked point.
\end{corollary}
\begin{proof}
    Set $\T':=\{((\Sigma, \partial\Sigma, \underline{z}, \underline{w}), u, F, e)\}$ where $((\Sigma,\partial\Sigma, \underline{z},\underline{w}), u, ,F)$ as in Definition \ref{framedCurvewithCons} and $e\in V_\mu$. For such choice of $V_\mu$ and using the proof of Proposition \ref{C1-loc structure}, we get that $\T'$ is a topological $O(d)$-manifold with fibre-wise smooth structure on $\T'\rightarrow\F'$. Now after stabilizing by $ev^*TL$ and arguing as in Lemma $4.5$ of \cite{abouzaid2021complex}, it follows by Theorem \ref{lashof}, that we can find $V$ a real $O(d)$-representation so that $\T:=\T'\times V$ is a smooth $O(d)$-manifold together with a smooth submersion $ev:\T\rightarrow L$. Similarly as above, we set $\E$ to be the trivial vectorbundle  $V_\mu\oplus\mathcal{H}_d\oplus V$ over $\T$ and $s$ is given analogously to Definition \ref{KuranishiSecAMSII}.
\end{proof}

For applications we have to compare different Global Kuranishi charts coming from boundary components of moduli spaces of $J$-holomorphic maps. We choose to take the Global Kuranishi chart representation as given in Corollary \ref{ChartsDiskswithRightSymmetry}, ensuring that the first item of the above $\textbf{compatibility conditions}$ is satisfied. Yet, we still have to modify the construction of these Global Kuranishi charts to ensure that the second and third items are also satisfied. We do this in a coherent way using the Lexographic partial order $\precsim$ on the monoid $\mathfrak{G}:=\{(\Omega(\beta), k, l)\}$ where $\beta\in\pi_2(X, L)$, $k,l$ non-negative integer, and $\Omega$ is as in Lemma \ref{Integral symplectic form}. In what follows, we fix notation and for $\beta_i\in\pi_2(X, L)$, $k_i,l_i$ as above, denote by:
\begin{enumerate}
    \item $d_i:=\Omega(\beta_i)$ and $d=\Omega(\beta)$.
    \item $G_i:=O(d_i)$ and $G=O(d)$.
    \item $\F_i\equiv\F'_{k_i+1, l_i}(\alpha_i)$ and $\F\equiv\F'_{k+1, l}(\alpha)$ be the space of framing, with a point constraint, as given in Lemma \ref{SpaceofFramingConstraint}.
    \item $\mathcal C_i\equiv\mathcal{C}'_{k_i+1, l_i}(\alpha_i)\rightarrow\F'_{k_i+1, l_i}(\alpha_i)$ and $\mathcal C\equiv\mathcal{C}'_{k+1, l}(\alpha)$ the corresponding universal curves.
\end{enumerate}

Let $\phi_i\in\F_i$ for $i=1,2$, and observe that we can write $\phi_i(z)=[f^i_j(z)]_{j=0, \dots, d_i}$ where $f^i_0, \dots, f_{d_i}^i$ is an $\R$-basis of $H^0(\Sigma_i,\partial\Sigma_i;\mathcal{L}_{u_i}, \mathcal{L}_{u_i, \R})$ as in the proof of Proposition \ref{FramingMap}. We first extend $f^i_j$ to $\Sigma_1\vee\Sigma_2$, where for simplicity, we define $\Sigma_1\vee\Sigma_2$ by gluing the last boundary marked point $z^1_{k_1}\in\Sigma_1$, to the zeroth boundary marked point $z^2_0\in\partial\Sigma_2$ as follows:
\[
\tilde f^1_j(z)=
\begin{cases}
    f^1_j(z) & \text{if } z\in\Sigma_1\\
    f^1_j(z^1_{k_1}) & \text{if } z\in\Sigma_2
\end{cases}
\]
\[
\tilde f^2_j(z)=
\begin{cases}
    f^2_j(z^2_0) & \text{if } z\in\Sigma_1\\
    f^2_j(z) & \text{if } z\in\Sigma_2
\end{cases}
\]
Now let $j=0,\dots, d_1$ be the first index so that $f^1_{j}(z^1_{k_1})\neq 0$ and let $\lambda\in\R^*$ be such that $\lambda f^2_0(z_0^2)=f^1_{j}(z^1_{k_1})$. We define $\phi_1\vee\phi_2:\Sigma_1\vee\Sigma_2\rightarrow \CP^{d_1+d_2}$ by
\begin{equation*}
    \hspace{-1em}\phi_1\vee\phi_2(z):=[\tilde f^1_j(z)+\lambda\tilde f_0^2(z)-\lambda f^2_0(z_0^2):\tilde f_0^1(z):\dots:\tilde f^1_{j-1}(z):\tilde f^1_{j+1}(z):\dots:\tilde f^1_{d_1}(z): \lambda \tilde f^2_{1}(z):\dots:\lambda\tilde f^2_{d_2}(z)]
\end{equation*}

Such construction provides us with an embedding of $\F_1\times_{ev}\F_1\hookrightarrow\F$ compatible with the diagonal embedding $G_1\times G_2\hookrightarrow G$ describing one component of the possible $k_1(k_2+1)+k_2(k_1+1)$ boundary components of $\F$ represented by $(\beta_1, \beta_2)$; which all can be defined in an analogous manner.

In the view of Lemma \ref{choosingMu}, we abuse terminology and call such $V:=V_\mu=: V_{k+1, l}(\beta)$ a finite-dimensional approximation scheme for $\M_{k+1, l}(X, L; J, \beta)$. Similarly for $i=1,2$, we denote by $V_i$ a finite-dimensional approximation scheme of $\M_{k_i+1, l_i}(X, L;J, \beta_i)$ as given in Lemma \ref{choosingMu}.

\begin{lemma}\label{ChoiceFiniteDimApp}
    There exists a finite-dimensional approximation scheme $V$ for the moduli space $\M_{k+1, l}(X, L; J, \beta)$ so that if for $i=1,2$ and for any $\beta_i, k_i, l_i$ as above, and $V_i$ is a finite-dimensional approximation scheme of $\M_{k_i+1, l_i}(X, L;J, \beta_i)$, then $V_1\oplus V_2\leq V$ is a $G$-subrepresentation.
\end{lemma}
\begin{proof}
    We argue by upward induction on the discrete monoid $\mathfrak{G}$.
    The base case follows from Lemma \ref{finiteAppScheme} and Lemma \ref{consMetric}.\\
    As for the inductive step, denote by $V_i$ the finite-dimensional approximation scheme for $\M_{k_i+1,l_i}(X, L;J,\beta_i)$ and using Lemma \ref{finiteAppScheme}, let $V'$ be a finite-dimensional approximation scheme for $\M_{k+1, l}(X, L; J, \beta)$. Denote by $\iota:\F_1\times_{ev}\F_2\hookrightarrow\F$ and $\iota^*\mathcal C\rightarrow\F_1\times_{ev}\F_2$. Then, using $\iota$ and the definition of differential forms on nodal Riemann surfaces, we get the following commutative diagram:
    \[\begin{tikzcd}
	{\bigwedge^{0,1}_{\iota^*\mathcal{C}/\mathcal{F}_1\times_{ev}\mathcal{F}_2}\otimes_{\mathbb{C}}TX} && {\bigwedge^{0,1}_{\mathcal{C}/\mathcal{F}}\otimes_{\mathbb{C}}TX} \\
	\\
	{\iota^*\mathcal{C}\times X} && {\mathcal{C}\times X}
	\arrow[hook, from=1-1, to=1-3]
	\arrow[from=1-1, to=3-1]
	\arrow[from=1-3, to=3-3]
	\arrow[hook, from=3-1, to=3-3]
\end{tikzcd}\]
    Inspired by Lemma $2.6$ of \cite{rezchikov2022integralarnoldconjecture}, we give the following construction:\\
Recall that, $G_1\times G_2:=O(d_1)\times O(d_2)$, $G:=O(r_1+r_2)$ and denote by $\lambda:G_1\times G_2\hookrightarrow G$ the diagonal embedding. We have, $\iota:N:=\iota^*\mathcal C\times X\hookrightarrow M:=\mathcal C\times X$ is a $(G_1\times G_2)$-equivariant embedding, where $G_1\times G_2$ acts on $M$ through $\lambda$. Now let, $E_{12}
:=
\operatorname{pr}_1^*V_1
\oplus
\operatorname{pr}_2^*V_2
\longrightarrow N$ and denote by $r_i:=\operatorname{rk}_\R(V_i)$ and set $r:=r_1+r_2$. As $O(d_1)$ and $O(d_2)$ are compact, we can choose and fix equivariant metrics on $V_1$ and $V_2$. Let, $F_i:=\operatorname{Fr}(V_i)\longrightarrow N_i:=\mathcal C_i\times X$ be their orthonormal frame bundles. Then, $F_H
:=
\operatorname{pr}_1^*F_1
\times_N
\operatorname{pr}_2^*F_2
\longrightarrow N$ is a $(G_1\times G_2)$-equivariant principal
$O(r_1)\times O(r_2)$-bundle. Under the block-diagonal inclusion $O(r_1)\times O(r_2)\hookrightarrow O(r)$, there is a canonical $(G_1\times G_2)$-equivariant isomorphism $\operatorname{Fr}(E_{12})
\cong
F_H\times_{O(r_1)\times O(r_2)}O(r)$, given by $[(u_1,u_2),A]
\longmapsto
(u_1\oplus u_2)\circ A$. Consider a finite-dimensional
orthogonal $(G_1\times G_2)$-representation $W$ and an $(G_1\times G_2)$-equivariant morphism of principal $O(r)$-bundles $\Theta:
\operatorname{Fr}(E_{12})
\longrightarrow
N\times\operatorname{St}_r(W)$, where $\operatorname{St}_r(W)$ denotes the Stiefel manifold of orthonormal $r$-frames in $W$. Note that, $\Theta(u)=(y,\theta_u),
\theta_u:\mathbb{R}^r\hookrightarrow W$, and hence we can define, $\alpha:E_{12}\hookrightarrow N\times W$ given by, $\alpha_y(u\xi):=(y,\theta_u(\xi))$. $\alpha$ is a $(G_1\times G_2)$-equivariant vector-bundle embedding. Indeed, as $G_1\times G_2\leq G$ is a closed subgroup of a compact orthogonal group, there exist a finite-dimensional orthogonal $G$-representation $U$ and an $(G_1\times G_2)$-equivariant isometric embedding $j:W\hookrightarrow U$ such that, $j((g_1,g_2)\cdot w)=\lambda((g_1,g_2))\cdot j(w)$ for any $(g_1, g_2)\in G_1\times G_2$ and $w\in W$. Now consider the $G$-equivariant principal $G$-bundle $P_G:=M\times G\longrightarrow M$, with actions $(x,g)\cdot a=(x,ga)$, $a\cdot(x,g)=(ax,ag)$. Let, $E_U:=P_G\times_GU\longrightarrow M$ and $V:=E_U\oplus V'\longrightarrow M$ be the associated vectorbundles. There is a $(G_1\times G_2)$-equivariant vector-bundle embedding $\Phi:
\operatorname{pr}_1^*V_1
\oplus
\operatorname{pr}_2^*V_2
\oplus
\iota^*V'
\hookrightarrow
\iota^*V$ given by, $\Phi_y(v_1,v_2,v')
=
\left(
\left[(\iota(y),1_G),
j\!\left(\alpha_y(v_1,v_2)\right)\right],
v'
\right)$ for $y\in N$. Therefore, $\operatorname{pr}_1^*V_1
\oplus
\operatorname{pr}_2^*V_2
\oplus
\iota^*V'
\subseteq
\iota^*V$. Now set $V=
\bigl((X\times G)\times_GU\bigr)\oplus V'$.\\
    After repeating the above argument finitely many times, corresponding to all possible ways $\beta=\beta_1+\beta_2$, $k+1=k_1+k_2$ and $l=l_1+l_1$, we get the desired result.
\end{proof}
\begin{theorem}\label{ChartsforDisksCompatible}
    There exists a Global Kuranishi chart $(G, \T, \E, s)$ of $\M_{k+1, l}(X, L; J, \beta)$ such that for $i=1,2$ and for any $\beta_i, k_i, l_i$ as above, $\M_{k_i+1, l_i}(X, L; J, \beta_i)$ has a Global Kuranishi chart $(G_i, \T_i, \E_i, s_i)$ satisfying all the $\textbf{compatibility conditions}$.
\end{theorem}
\begin{proof}
    We argue by upward induction on $\mathfrak{G}$. For the base case, we apply Corollary \ref{ChartsDiskswithRightSymmetry}. As for the inductive step, denote by $\T^{top}_i:=\T_{k_i+1, l_i}^{top}(\beta_i):=\{((\Sigma, \partial\Sigma; \underline{z}, \underline{w}), u, F, e)\}$ to be the set of tuples where $((\Sigma, \partial\Sigma; \underline{z}, \underline{w}), u, F)$ is a framed marked disk as in Definition \ref{framedCurvewithCons} representing the class $\beta_i$ and $e\in V_i\equiv V_{k_i+1, l_i}(\beta_i)$ where $V_i$ is as in Lemma \ref{ChoiceFiniteDimApp}. Similarly, $\T^{top}:=\T^{top}_{k+1, l}(\beta):=\{((\Sigma, \partial\Sigma; \underline{z}, \underline{w}), u, F, e)\}$, where $e\in V\equiv V_{k+1, l}(\beta)$. Using Lemma \ref{ChoiceFiniteDimApp} and the proof of Proposition \ref{ChoosingN}, it follows that $\T, \T_i$ are topological $G, G_i$-manifolds respectively. Using the proof of Proposition \ref{C1-loc structure}, it follows that there exists $R_i$ a real $O(d_i)$-representation  so that $\T^{top}_i\times R_i$ is a smooth $G_i$-manifold together with a smooth submersion $ev:\T_i^{top}\times R_i\rightarrow X^{l_i}\times L^{k_i+1}$. With such choices, we get a commutative diagram:
    \[\begin{tikzcd}
	{(\mathcal{T}^{top}_1\times R_1)\times_{ev}(\mathcal{T}^{top}_2\times R_2)} && {\mathcal{T}^{top}\times(R_1\oplus R_2)} \\
	\\
            	{\mathcal{F}_1\times_{ev}\mathcal{F}_2} && {\mathcal{F}}
	\arrow[hook, from=1-1, to=1-3]
	\arrow[from=1-1, to=3-1]
	\arrow[from=1-3, to=3-3]
	\arrow[hook, from=3-1, to=3-3]
\end{tikzcd}\]
compatible with the diagonal embedding $G_1\times G_2\hookrightarrow G$, where the lower horizontal arrow is as constructed above. Now apply Theorem \ref{relativeLashof} to find $R$ a $G$-representation so that $\T_i:=\T^{top}_i\times (R_i\oplus R)$ and $\T:=\T^{top}\times (R_1\oplus R_2\oplus R\oplus R)$ are smooth $G_i, G$-manifolds respectively. From the commutativity of the diagram, it follows that the diagram
\[\begin{tikzcd}
	{\mathcal{T}_1\times_{ev}\mathcal{T}_2} && {\mathcal{T}} \\
	\\
	{\mathcal{F}_1\times_{ev}\mathcal{F}_{2}} && {\mathcal{F}}
	\arrow[hook, from=1-1, to=1-3]
	\arrow[from=1-1, to=3-1]
	\arrow[from=1-3, to=3-3]
	\arrow[hook, from=3-1, to=3-3]
\end{tikzcd}\]
commutes and compatible with the group embedding $G_1\times G_2\hookrightarrow G$.\\

As for the obstruction bundles, we set $\E_i$ to be the trivial bundle $V_i\oplus\mathcal H_{d_i}\oplus(R_i\oplus R)$ over $\T_i$ and $\E$ to be the trivial bundle $V\oplus\mathcal{H}_{d}\oplus(R_1\oplus R_2\oplus R\oplus R)$ over $\T$. After defining the Kuranishi sections $s_i:\T_i\rightarrow\E_i$ and $s:\T\rightarrow\E$ as in Definition \ref{KuranishiSecAMSII}, the third item of the compatibility conditions above, follows from Lemma \ref{ChoiceFiniteDimApp}. The fourth item of the compatibility conditions follows after observing that the Global Kuranishi chart $(G,\T,\E,s)$ as constructed for $\overline{\mathcal M}(\beta)$ when restricted to its boundary corresponding to degeneration of the form $\F_1\times_{ev}\F_2\hookrightarrow\F$ is identically equal to the fibre-product Global Kuranishi chart after stabilizing by $\iota^* T\T/(T\T_1\oplus T\T_2)\oplus\iota^*\E/(\E_1\oplus\E_2)\rightarrow \T_1\times_{ev}\T_2$ and hence the desired equivalence. As for the last compatibility condition, we refer to the Lemma below.
\end{proof}

\begin{lemma}
Let, $(G,\mathcal{T},\mathcal{E},s)$
be a Global Kuranishi chart for $\overline{\mathcal{M}}_{k+1, l}(\beta)$, and, for $i=1,2$, let $(G_i,\mathcal{T}_i,\mathcal{E}_i,s_i)$ be Global Kuranishi charts for $\overline{\mathcal{M}}_{k_i+1, l_i}(\beta_i)$, where
$\beta_1+\beta_2=\beta$. Suppose that the $r_1$-th marked point on the first component is identified with the $r_2$-th marked point on the second component. Assume that, after stabilizations and group enlargements, there is an equivariant isomorphism of Global Kuranishi charts $\Phi\colon
\left(
G_1\times G_2,\,
\mathcal{T}_1
\,{}_{ev^{1}_{r_1}}\!\times_{ev^{2}_{r_2}}
\mathcal{T}_2,\,
\mathcal{E}_1\oplus\mathcal{E}_2,\,
s_1\oplus s_2
\right)
\longrightarrow
\left.
(G,\mathcal{T},\mathcal{E},s)
\right|_{\mathcal{S}_{\beta_1,\beta_2}}$,
where $\mathcal{S}_{\beta_1,\beta_2}$ is the corresponding boundary stratum. If $\Phi$ is induced by an identification of the universal curves and universal maps, then $\Phi$ intertwines the evaluation maps. More precisely, if the $j$-th marked point of the glued curve comes from the $\rho_i(j)$-th marked point of the $i$-th component, then $ev_j\circ\Phi
=
ev^{i}_{\rho_i(j)}\circ pr_i$, for $i=1,2$.
\end{lemma}

\begin{proof}
For $i=1,2$, let $q_i\colon \mathcal{T}_i\longrightarrow\mathcal{F}_i$
be the natural map to the corresponding framing moduli space. Let $\pi_i\colon\mathcal{C}_i\longrightarrow\mathcal{F}_i$
be the universal curve, let $\sigma^i_j\colon\mathcal{F}_i\longrightarrow\mathcal{C}_i$
be the section corresponding to the $j$-th marked point, and let $U_i\colon\mathcal{C}_i\longrightarrow X$
be the universal map. The evaluation map on the thickening
$\mathcal{T}_i$ is therefore $ev^i_j
=
U_i\circ\sigma^i_j\circ q_i$. Similarly, let
$
q\colon\mathcal{T}\longrightarrow\mathcal{F}$
be the map associated to the ambient global Kuranishi chart, and let $\pi\colon\mathcal{C}\longrightarrow\mathcal{F}$, $\sigma_j\colon\mathcal{F}\longrightarrow\mathcal{C}$,
$U\colon\mathcal{C}\longrightarrow X$
denote the universal curve, its marked sections, and the universal map,
respectively. Thus $ev_j=U\circ\sigma_j\circ q$.
Consider the fiber-product thickening $\mathcal{T}_{12}
:=
\mathcal{T}_1
\,{}_{ev^1_{r_1}}\!\times_{ev^2_{r_2}}
\mathcal{T}_2$ and the corresponding fiber product of framing moduli spaces $\mathcal{F}_{12}
:=
\mathcal{F}_1
\,{}_{ev^{\mathcal F_1}_{r_1}}\!\times_
        {ev^{\mathcal F_2}_{r_2}}
\mathcal{F}_2$, where $ev^{\mathcal F_i}_j:=U_i\circ\sigma^i_j$.
The maps $q_i$ induce a map $q_{12}\colon\mathcal{T}_{12}\longrightarrow\mathcal{F}_{12}$, $q_{12}(t_1,t_2)=\bigl(q_1(t_1),q_2(t_2)\bigr)$.
Let $\Phi\colon
\mathcal{T}_{12}\longrightarrow
\left.\mathcal{T}\right|_{\mathcal{S}_{\beta_1,\beta_2}}$
be the identification of global Kuranishi charts obtained after
stabilization and group enlargement. By the compatibility of this
identification with the framing data, there is a map $\phi\colon
\mathcal{F}_{12}\longrightarrow
\left.\mathcal{F}\right|_{\mathcal{S}_{\beta_1,\beta_2}}$ such that $q\circ\Phi=\phi\circ q_{12}$. The map $\phi$ sends a pair of framed maps to the corresponding nodal
framed map.

Suppose that the $j$-th marked point of the glued curve comes from the
$\rho_i(j)$-th marked point of the $i$-th component, where $i\in\{1,2\}$.
The identification of the universal curves over $\phi$ preserves the
universal maps and the marked sections. Consequently, $ev^{\mathcal F}_j\circ\phi
=
ev^{\mathcal F_i}_{\rho_i(j)}\circ pr_i$.
It follows that
\begin{align*}
ev_j\circ\Phi
&=
ev^{\mathcal F}_j\circ q\circ\Phi\\
&=
ev^{\mathcal F}_j\circ\phi\circ q_{12}\\
&=
ev^{\mathcal F_i}_{\rho_i(j)}
   \circ pr_i\circ q_{12}\\
&=
ev^{\mathcal F_i}_{\rho_i(j)}\circ q_i\circ pr_i\\
&=
ev^i_{\rho_i(j)}\circ pr_i.
\end{align*}
This proves the compatibility of the evaluation maps.

At the nodal marked points, the equality $ev^1_{r_1}\circ pr_1
=
 ev^2_{r_2}\circ pr_2$ holds by the definition of the fiber product $\mathcal{T}_{12}$. Hence the evaluation map at the node is also compatible with the boundary identification.
\end{proof}

\section{Global Kuranishi Charts for the Moduli Space of Floer-Morse Trajectories}
\subsection{Morse Trajectories}
We start by investigating the classical Morse theory case and prove results in this setting that we will be used in the next subsection. We follow closely the discussion in \cite{fukaya1997zero}, \cite{Fukaya2009}.
\subsubsection{Preliminaries}
    \begin{definition}
    A ribbon tree is a tuple $t:=(T, \iota, rt)$ where, 
    \begin{enumerate}
        \item $T$ is planar, connected cycle-free rooted tree.
         \item $\iota:T\rightarrow\mathbb{D}$ is a topological embedding.
        \item $\iota^{-1}(\partial\mathbb{D})$ is the set of vertices of valency 1.
        \item $rt\in\iota^{-1}(\partial\mathbb{D})$ is a distinguished point called the root of $t$.
    \end{enumerate}
\end{definition}
The choice of root of $T$ is equivalent to a choice of ribbon structure as then we enumerate the vertices and cyclically order the edges after pulling-back, via $\iota$, the counter-clockwise orientation on $\mathbb{D}\subset\C$.
\begin{definition}
    For a ribbon tree $(T, \iota, rt)$ we denote by $C^{i}(T)$ the $i^{th}$-cell of $T$ where $i=0, 1$.
    \begin{enumerate}
        \item We denote by $C^{0}_{ext}(T):=\iota^{-1}(\partial\mathbb{D})$ the set of external vertices.
        \item We denote by $C^{0}_{int}(T):=C^{0}(T)\setminus C^{0}_{ext}(T)$ the set of internal vertices.
         \item We have a distinguished element $rt\in C^{0}_{ext}(T)$ called the root of $T$.
        \item We denote by $C^{1}_{ext}(T)$ the set of all edges adjacent to an external vertex.
        \item We denote by $C^{1}_{int}(T):=C^{1}(T)\setminus C^{1}_{ext}(T)$ the set of all internal edges.
        \item Elements of $C^{1}_{ext}(T)$ are called the leaves of $T$, with a distinguished leaf $e_0$ adjacent to the root.
    \end{enumerate}
\end{definition}
\begin{remark}
    For a ribbon tree $t=(T,\iota, rt)$ we sometimes abuse notation and refer to $C^i_*(t)\equiv C^i_*(T)$, for $i=0,1$ and $*=ext, int$.
\end{remark}
\begin{definition}
    A ribbon tree $(T, \iota, rt)$ is said to be stable if every internal vertex has valency at least $3$.
\end{definition}
In this section, we will only work with stable ribbon trees.
\begin{definition}
    Two ribbon trees $(T_{i}, \iota_{i}, rt_{i})$ for $i=1,2$, are said to be equivalent or of the same combinatorial type, if $T_{1}\cong T_{2}$ are homeomorphic and $\iota_{1}$, $ \iota_{2}$ are isotopic.
\end{definition}
\begin{definition}
    For $k$ non-negative integer, we denote by $G_{k+1}$ the set of all  isomorphism classes of (stable) ribbon trees with $(k+1)$-leaves. 
\end{definition}
\begin{definition}
    We denote by $G_0$ the unique tree with two external verices, one of them the root, no internal vertices, and one edge connecting the vertices denoted by $e_0$.
\end{definition}
In order to give $G_{k+1}$ a topology, namely the Gromov-Hausdorff topology, we give each ribbon tree a metric.
\begin{definition}
\begin{enumerate}
    \item Let $t=(T, \iota, rt)$ be a ribbon tree, a metric on $t$ is a map $l:C^{1}_{int}(T)\rightarrow (0, \infty)$.
    \item We denote by $Gr_{k+1}(t):=\{l:C^{1}_{int}(T)\rightarrow (0, \infty)\}$ the set of all metrics on $t$.
    \item We denote by $Gr_{k+1}:=\cup_{G_{k+1}}Gr_{k+1}(t)$.
\end{enumerate}
\end{definition}
We give $Gr_{k+1}$ a sequential topology as follows:\\
Let $l_{i}\in Gr_{k+1}(t)$ be a sequence and suppose that $l_{i}(e)\rightarrow l_{\infty}(e)$, for every $e\in C^{1}_{int}(T)$. Let $T^{'}$ be tree constructed from $T$ by shrinking all edges $e$ so that $l_{\infty}(e)=0$. Then $l^{'}:=l_{\infty|T'}\in Gr_{k+1}(t')$ where $t':=(T', \iota|_{T'}, rt)$. \\
A theorem of Stasheff says that $Gr_{k+1}$ is homeomorphic to $\R^{k-2}$ and $Gr_{k+1}:=\cup_{G_{k+1}}Gr_{k+1}(t)$ is a cellular decomposition (c.f. Lemma $1.3$ of \cite{fukaya1997zero} and the discussion above it). In fact, Fukaya-Oh upgraded this result as follows:
\begin{theorem}[Theorem $14.1$ \cite{fukaya1997zero}]\label{StasheffSmooth}
    $Gr_{k+1}$ is a smooth manifold with corners of $\dim_{\R}=k-2$ and each of its strata $Gr_{k+1}(t)$ is a smooth manifold with corners of $\dim_{\R}=k+1-\Sigma_{C^{0}_{int}(T)}|v|$, where $|v|$ is the valency of $v$.
\end{theorem}

We will also allow internal edges to acquire infinite length. In other words, we will compactify $Gr_{k+1}$ by adding $\textit{broken trees}$.
\begin{definition}
    A broken tree of $(k+1)$-leaves is a tuple $(T, \iota, rt; l)$ where,
    \begin{enumerate}
        \item $(T, \iota, rt)\in G_{k+1}$.
        \item $l: C^1_{int}(T)\rightarrow (0, \infty]$. 
    \end{enumerate}
   with the extra decoration on every $e\in C^1_{int}(T)$ given by a non-negative integer $b(e)$, such that $b(e)>0$ if $l(e)=\infty$. We call $b(e)$ the number of breakings at edge $e$ and $\Sigma_{C^1_{int}(T)}b(e)$ the total number of breakings of the tree $T$.
\end{definition}
\begin{definition}
    Denote by $\overline{Gr}_{k+1}$ the space, for the Gromov-Hausdorff topology, of all broken metric, ribbon trees with $(k+1)$-leaves.
\end{definition}

Broken trees of $(k+1)$-leaves can be constructed inductively as follows:\\
Suppose that $k+1=k_1+k_2$ and for $i=1,2$, let $t_i=(T_i, \iota_i, rt_i; l_i)$ be a broken tree of total breakings equal to $b_i$. We define $t=(T, \iota, rt; l)$ by setting $T:=T_1\sqcup T_2/\sim$ where $\sim$ identifies the root of $t_1$ with one of the $k_2$-leaves of $t_2$ other than its root. The resulting tree will be a broken tree of $(b_1+b_2+1)$-breakings. Such construction provides us with one of the possible $k_2$-$\textit{gluing}$ maps $\overline{Gr}_{k_1+1}\times \overline{Gr}_{k_2+1}\rightarrow\overline{Gr}_{k+1}$.
\begin{theorem}[Section $1.9$ \cite{boardman2006homotopy}, Theorem $2.8$ \cite{ma2010geometric}]
     $\overline{Gr}_{k+1}$ is a compact topological space homeomorphic to the $\textit{Stasheff}$ associahedra as in \cite{stasheff1963homotopy}.
\end{theorem}

Now following the discussion found in \cite{abouzaid2010geometric} and Section $7.1$ of \cite{abouzaid2011topological}, we can give $\overline{Gr}_{k+1}$ a smooth structure of a manifold with corners, realizing the $\textit{Stasheff}$ associahedra. 
Let, $t=(T,\iota,rt_1),
t^\infty=(T^\infty,\iota^\infty,rt^\infty)$
be elements of $G_{k+1}$, and let $\pi\colon T\longrightarrow T^\infty$ be a surjective map of rooted ribbon trees which preserves the ribbon structures and the distinguished exterior vertices, and for which the inverse image of every vertex of $T^\infty$ is connected. Define $E_{\mathrm{fin}}(\pi)
:=
\left\{
e\in C_{\mathrm{int}}^1(T)
\;\middle|\;
\pi(e)\in C^0(T^\infty)
\right\}$.
Namely, $E_{\mathrm{fin}}(\pi)$ is the set of internal edges of $T$ that are collapsed by $\pi$. In particular, $\#E_{\mathrm{fin}}(\pi)
=
\#C_{\mathrm{int}}^1(T)
-
\#C_{\mathrm{int}}^1(T^\infty)$.

The cell associated with $\pi$ is $Gr(t\longrightarrow t^\infty)
:=
(0,\infty)^{E_{\mathrm{fin}}(\pi)}$,
and its closure is canonically identified with $\overline{Gr}(t\longrightarrow t^\infty)
:=
[0,\infty]^{E_{\mathrm{fin}}(\pi)}$.
We denote the corresponding characteristic map by $\kappa_\pi\colon
\overline{Gr}(t\rightarrow t^\infty)
\longrightarrow
\overline{Gr}_{k+1}$.
Under this map, an edge $e$ for which $l(e)=0$ is contracted,
whereas an edge for which $l(e)=\infty$ becomes a breaking edge. We equip $[0,\infty]$ with the smooth structure induced by the
bijection $\beta\colon[0,\infty]\longrightarrow[0,1], \beta(l)=e^{-l}$, where $e^{-\infty}:=0$.
\begin{definition}
Let $F$ be a Fechet space. A continuous map $f\colon\overline{Gr}_{k+1}\longrightarrow F$
is called smooth if, for every
$\pi\colon t\to t^\infty$, the map $f_\pi:=f\circ\kappa_\pi\colon
[0,\infty]^{E_{\mathrm{fin}}(\pi)}
\longrightarrow F$ is smooth with respect to the above smooth structure. 
\end{definition}
As for the case of the universal curve $\widetilde{\overline{Gr}}_{k+1}
:=
\left\{
(p, (t,l))
\;\middle|\;
(t,l)\in\overline{Gr}_{k+1},
p\in T
\right\}$ and $\widetilde f\colon
\widetilde{\overline{\mathrm{Gr}}}_{k+1}\times L
\rightarrow\mathbb{R}$ a continuous function. Fix a cell $\pi\colon t\to t^\infty$ and a flag $(e,v)$ of $T$,
that is, $e\in C^1(T)$ and $v\in\partial e$. For every metric $l$, let $\vartheta_{e,v}^{\,l}\colon
(e,l)\longrightarrow I_{e,v}(l)\subset\mathbb{R}$
be the unique length-preserving parametrization that sends $v$ to $0$. Define $P_{\pi,e,v}
:=
\left\{
(l,s)\in
[0,\infty]^{E_{\mathrm{fin}}(\pi)}\times\mathbb{R}
\;\middle|\;
s\in I_{e,v}(l)
\right\}$.
The function induced by $\widetilde f$ in this flag parametrization is $f_{\pi,e,v}\colon P_{\pi,e,v}\times L\longrightarrow\mathbb{R}$, given by $f_{\pi,e,v}(l,p,x)
:=
\widetilde f\left(
\kappa_\pi(l),
(\vartheta_{e,v}^{\,l})^{-1}(p),
x
\right)$.
\begin{definition}
We say that $\widetilde f$ is smooth if, for every
$\pi\colon t\to t^\infty$ and every flag $(e,v)$, there exists a
function $\widehat f_{\pi,e,v}\colon
[0,\infty]^{E_{\mathrm{fin}}(\pi)}
\times\mathbb{R}\times L
\longrightarrow\mathbb{R}$ such that $\left.
\widehat f_{\pi,e,v}\right|_{P_{\pi,e,v}\times L}
=
f_{\pi,e,v}$, and such that the function $(q,p,x)\longmapsto
\widehat f_{\pi,e,v}
\left(
(-\log q_e)_{e\in E_{\mathrm{fin}}(\pi)},
p,x
\right)$ is smooth on $[0,1]^{E_{\mathrm{fin}}(\pi)}
\times\mathbb{R}\times L$.
\end{definition}
\begin{remark}
    In particular, if $f$ is constant in the fibre-direction of some small enough disk bundle of $\overline{Gr}_{k_1+1}\times\overline{Gr}_{k_2+1}\rightarrow\overline{Gr}_{k+1}$ and smooth on each strata of $\overline{Gr}_{k+1}$ then, $f$ is a smooth map. 
\end{remark}

\subsubsection{Domain-Dependent Morse Trajectories}
Suppose that $L$ a closed smooth manifold and $(f, g)$ a Morse-Smale pair on $L$. 
\begin{example}
     Let $\underline{x}=\{x_{0}, x_{1}, x_{2}\}$ be critical points of $f$ and consider the unique trivalent tree $t\in G_3$. Then by the Morse-Smale condition, we have that the unstable manifold $W^{-}_{x_{1}}$ and the stable manifold $W^{+}_{x_{0}}$ intersect transversaly, while by the uniqueness part of the fundamental theorem of ODEs, it follows that $W^{-}_{x_{1}}$ and $W^{-}_{x_{2}}$ only intersect if $x_{1}=x_{2}$ in which case their intersection is not transversal.
\end{example}

In order to achieve transversality, we will perturb the $\textit{Morse trajectory equations}$ away from neighborhoods of the vertices as in \cite{charest2015floer}.
We start by fixing a coherent coordinate system on edges of trees.
\begin{definition} \label{coordinates}
    For $(t, l)\equiv(T, \iota, rt; l)\in Gr_{k+1}$, we give the edges of $T$ coordinates as follows:
    \begin{enumerate}
        \item If $e\in C^{1}_{ext}(T)$ and does not contain the root, we give coordinates on $e$ so that $e\cong(-\infty, 0]=: I_e$.
        \item If $e\in C^{1}_{ext}(T)$ and contains the root, we give coordinates on $e$ so that $e\cong[0, \infty)=: I_e$.
        \item If $e\in C^1_{int}(T)$, we give coordinates on $e$ so that $e\cong[0, l(e)]=:I_e$.
    \end{enumerate}
    We denote the coordinates on the edges by $s$ and write $s\in T$.
\end{definition}
\begin{definition}
\begin{enumerate} \label{universalFamily}
    \item We denote by $\widetilde{Gr}_{k+1}:=\{(t, l, s): t=(T, \iota, rt)\in G_{k+1}, l\in Gr_{k+1}(t), s\in T\}\rightarrow Gr_{k+1}$ the universal family of metric ribbon trees, together with the projection map $(t, l, s)\mapsto (t, l)$.
    \item For $t\in G_{k+1}$, we denote by $\widetilde{Gr}_{k+1}(t)$ the pre-image of the open cell $Gr_{k+1}(t)\subseteq Gr_{k+1}$ containing $t$ by the above projection map.
\end{enumerate}
     
\end{definition}
\begin{definition} \label{T-depMorse}
    Let $t=(T, \iota, rt)\in G_{k+1}$ and $\tilde f\in C^{\infty}(\widetilde{Gr}_{k+1}(t)\times L)$. $\tilde f$ is said to be a Morse-function associated to $f$ if
    \begin{enumerate}
        \item $\tilde f(s,x)\equiv f(x)$ for every $s\in T$ and for every $x$ on some neighborhood of each critical point of $f$.
        \item $\tilde f(s, .)\equiv f(.)$ for every $s$ in a neighborhood of each vertex of $T$. 
    \end{enumerate}
\end{definition}
Let $t=(T, \iota, rt)\in G_{k+1}$,  $\underline{x}:=\{x_0,\dots,x_k\}$ a finite set of critical points of $f$ and $\tilde f\in C^\infty(\widetilde{Gr}_{k+1}(t)\times L)$ as in Definition \ref{T-depMorse}.  

\begin{definition} \label{ModuliMorseFixedDomain}
    We denote by $\mathcal{M}(t, \underline{x}; \tilde f)\equiv\mathcal{M}(t, \underline{x}; \tilde f,g):=\{(l, u): l\in Gr_{k+1}(t), u:T\rightarrow L\}$ where $u$ is a continuous map such that for each $e\in C^{1}(T)$ we require $u_{e}:I_e\rightarrow L$, the restriction of $u$ to the edge $e$, to solve the $\textit{domain-parametrized}$ Morse trajectory equation:
        \begin{equation} \label{MorseTraj}
            \frac{d}{ds}u_{e}+\nabla^g \tilde f(s, u_{e}(s))=0
        \end{equation}
such that
\begin{enumerate}
    \item If $e\in C^1_{ext}(T)$ and $e\neq e_0$ then, $u_e(-\infty)=x_i$ for some $i=1,\dots,k$.
    \item If $e=e_0$ then, $u_{e_0}(\infty)=x_0$.
\end{enumerate}
\end{definition}
\begin{remark}
\begin{enumerate}
    \item In the above definition we have abused notation and wrote $u_{e}(\pm\infty)$, since solutions of the Morse trajectory equation have exponential decay convergence near critical points and hence such limit always exists.
    \item Moreover, from the continuity of $u:T\rightarrow L$, it follows that we do not have to specify the incident conditions on the internal vertices.
    \item As we require internal edges to have finite non-zero length and hence in particular, $t$ is an $\textit{unbroken}$ tree, it follows that none of the internal vertices are mapped by $u_e$ to critical points of $f$.
\end{enumerate}
\end{remark}
\begin{lemma}\label{genericFixedComb}
    For a generic choice of $\tilde{f}\in C^\infty(\widetilde{Gr}_{k+1}(t)\times L)$ as in Definition \ref{T-depMorse}, $\mathcal{M}(t, \underline{x}; \tilde f)$ is a smooth manifold.
\end{lemma}
\begin{proof}
    In this proof we will do two levels of perturbations, one along external edges so that unstable manifolds intersect transversely and the second along internal ones. To avoid notational clutter, for $e\in C^1(T)$ we denote by $f_e:I_e\rightarrow L$ the smooth map given by $f(\gamma_e(s))$ where, $s\in I_e$ is one of the intervals as in Definition \ref{coordinates} and $\gamma_e:I_e\rightarrow L$ is a smooth map such that $\frac{d}{ds}\gamma_e+\nabla^gf(\gamma_e)=0$. We start by perturbing along external edges. \\
    
    Using the ribbon structure on $t$ we enumerate the internal vertices $v_{1}, \dots, v_{m}$ which are adjacent to external edges. Moreover, for each such vertex $v_{i}$, we have a cyclic ordering, say $e^{1}_{i}, \dots, e^{k_{i}}_{i}$, on all the external edges adjacent to it. By Sard's Theorem, we can find $h\in C^{\infty}(L)$ so that $h_{e^{j}_{i}}\equiv f_{e^{j}_{i}}$ for every $s\in(-\infty, 2]\cup (1, 0]$ or $s\in[0,1)\cup[2, \infty)$ if ${e^{j}_{i}}=e_0$, $h_{e}\equiv f_{e}$ on each internal edge $e\in C^{1}_{int}(T)$ of the underlying tree $T$, so that after integrating the vectorfield $\nabla^g h$, the images of all unstable manifolds have images intersecting transversely at $v_{i}$. For $a=0,\dots, k$, we denote by $\Tilde{W}^{-}_{x_{a}}(e^{j}_{i})$ the image of $W^{-}_{x_{a}}$, the unstable manifold of $f$ at $x_a$, by the flow of $\nabla^g h|_{e^{j}_{i}}$ and thus we can assume that $\cap_{j=1, }^{k_{i}}\Tilde{W}^{-}_{x_{a}}(e^{j}_{i})$ is a smooth submanifold of $L$ with the understanding that if $e^{j}_{i}=e_0$ then $\Tilde{W}^{-}_{x_{a}}(e^{j}_{i})=W^{+}_{x_{0}}$. Define $\tilde h\in C^\infty(\widetilde{Gr}_{k+1}(t)\times L)$ by $\tilde h(s, .)\equiv h(.)$ for every $s\in T$ and observe that, by construction, such $\tilde h$ satisfies the conditions of Definition \ref{T-depMorse}.\\
    
    Now we realize $\mathcal{M}(t, \underline{x}; \tilde h)$ as a fibre product as follows:\\
    For each $v_{i}$ as above, let $e_{i, 1}, \dots, e_{i, d}$ be the internal edges forming the unique minimal path joining $v_{1}$ to $v_{i}$. Such path exists as $T$ is connected, of finite type and unique as $T$ is cycle-free. For each such $e_{i, j}$, denote by $W_{i, j}(s):=-\nabla^g\tilde h(s, .)|_{e_{i, j}}=-\nabla^g f|_{e_{i, j}}$ and consider the smooth map $\exp_{i}: L\times \widetilde{Gr}_{k+1}(t)\rightarrow L$ given by 
    \begin{equation*}
        \exp_{i}(p, (s,l)):=\exp(W_{i, 1}(s))\circ\dots\circ \exp(W_{i, d}(s))(p),
    \end{equation*}
    with the understanding that $\exp(W_i(s))\circ\exp(W_j(s))(p)\equiv\exp_{\exp_p(W_j(s))}(W_i(s))$.
    We define $Exp_{\tilde h}:L\times \widetilde{Gr}_{k+1}(t)\rightarrow L^{m}$ by 
    \begin{equation} \label{Exponentiation}
        Exp_{\tilde h}(p, (s,l)):=(\exp_{i}(p, (s,l))_{i=1, \dots, m}
    \end{equation}
    As $\tilde h (s, .)\equiv f$ for every $s$ on an internal edge of $T$, it follows that 
    \begin{equation*}
        \mathcal{M}(t, \underline{x}; \tilde h)=Exp_{\tilde h}^{-1}(\prod_{i=1}^{m} \cap_{j=1}^{k_{i}}\Tilde{W}^{-}_{x_{a}}(e^{j}_{i})).
    \end{equation*}
    
    Now for each $v_{i}$ as above, let $O_{i}$ be an open subset in $e_{i, 1}$ the edge adjacent to $v_{i}$ and not intersecting $v_{i}$. We require that, if $v_{i}=v_{1}$ then $O_{1}$ is empty. Moreover, we require that each $O_{i}$ is small enough so that, if $\gamma_e$ is a Morse trajectory for $f$ as above then $\gamma_{e_{i,1}}(O_i)$ is contained in a small enough open set of $L$ so that by Poincare Lemma and the non-degeneracy of the metric on $L$, the choice of a vectorfield supported on $\gamma_{e_{i,1}}(O_i)$ is interchangeable with a choice of function supported on $\gamma_{e_{i,1}}(O_i)$ of gradient vector equal to such vectorfield. With this in mind, denote by $C^{N}(O_{i})$ the Banach space of $C^{N}$-functions on $\widetilde{Gr}_{k+1}(t)\times L$ supported on $O_{i}\times L$. Namely, if $s\notin O_i$ and $f_i\in C^N(O_i)$ then $f_i(s, .)\equiv 0$. Consider the universal moduli space $\mathcal{M}(t, \underline{x}):=\cup_{\prod_{i}^{m}C^{N}(O_{i})}\mathcal{M}(t, \underline{x},  \tilde h+f_{1}+\dots+f_{m})$. In order to use Sard-Smale Theorem, we show that $\mathcal{M}(t, \underline{x})$ is a Banach manifold. To this end, we first realize $\mathcal{M}(t, \underline{x})$ as a fibre product. Indeed, observe that $\mathcal{M}(t, \underline{x})=\overline{Exp}^{-1}(\prod_{i=1}^{m} \cap_{j=1}^{k_{i}}\Tilde{W}^{-}_{x_{a}}(e^{j}_{i}))$ where $\overline{Exp}:L\times \widetilde{Gr}_{k+1}(t)\times\prod_{i}^{m}C^{N}(O_{i})\rightarrow L^{m}$ is given by
    \begin{equation*}
        \overline{Exp}(p, (s, l); \Bar{f}):=Exp_{\tilde h+f_{1}+\dots+f_{m}}(p, (s,l))
    \end{equation*}
    Thus, it suffice to show that $\overline{Exp}$ is a submersion.\\
    
    Indeed, let $(p, (s,l), \Bar{f})\in L\times \widetilde{Gr}_{k+1}(t)\times\prod C^{N}(O_{i})$ and suppose that $\overline{Exp}(p, (s,l), \Bar{f})=(p_{i})_{i=1, \dots, m}$. Fix a vector $(X_{i})\in\oplus_{i}T_{p_{i}}L$. Note that, by our choice of $O_i$ as in the above paragraph, namely that using the Poincare Lemma we are identifying vectorfields of support on $\gamma_e(O_i)$ and functions supported on $\gamma_e(O_i)$ whose gradient is equal to such vectorfield, it suffice to find $\Bar{g}\in\prod C^{N}(O_{i})$, a vector $V\in T_{p}L$ and a vector $\zeta\in T_l Gr_{k+1}(t)$ so that $d_{(p, (s, l), \Bar{f})}Exp_{\tilde h+(f_{1}+g_1)+\dots+(f_{m}+g_m)}(V, \zeta)=(X_{i})$. We give an ordering on $v_{1}, \dots, v_{m}$ by $v_{i}\leq v_{j}$ if and only if, the minimal path joining $v_{1}$ to $v_{j}$ passes through $v_{i}$. Notice that the choice of $g_{i}$ does not affect the $X_{j}$ for $j<i$ this is due to the choice of $O_{i}$. Moreover for $v_{1}$, we have that $p_{1}=p$ and thus we can take $V=X_{1}$ and $\zeta=0$. 
    Now suppose that we have chosen $g_{1}, \dots, g_{i}$ so that the first $i^{th}$-components of $d_{(p, (s,l), \Bar{f})}Exp_{\tilde h+(f_{1}+g_1)+\dots+(f_{i}+g_i)+f_{i+1}+\dots+f_m}(V,0)$ is $X_{1}, \dots, X_{i}$ and the $(i+1)$-component is $Y_{i}$. We set $g_{j}=0$ for $j>i+1$ and choose $g_{i+1}$ so that the differential becomes $X_{i}-Y_{i}$. Continuing in this fashion, it follows that $\overline{Exp}$ is a submersion and hence, $\overline{Exp}^{-1}(\prod_{i=1}^{m} \cap_{j=1}^{k_{i}}\Tilde{W}^{-}_{x_{a}}(e^{j}_{i}))$ is a smooth Banach manifold. Now as the projection map $\mathcal{M}(t, \underline{x})\rightarrow\prod_{i}^{m}C^{N}(O_{i})$ is Fredholm and countable intersection of Baire sets is a Baire set, the result follows.
\end{proof}
We introduce the following partial order on $G_{k+1}$.
\begin{definition}
    Let $t, t'\in G_{k+1}$. We write $t>t'$ if $t'$ is obtained from $t$ by shrinking some interior edge.
\end{definition}
\begin{definition}
    For $t\in G_{k+1}$ we denote by $Gr_{k+1}^{+}(t):=\{l:C^{1}_{int}(T)\rightarrow[0, \infty)\}$.
\end{definition}
Thus, from the (sequential) topology given on $Gr_{k+1}$, it follows that $Gr_{k+1}^{+}(t)=\cup_{t>t'}Gr_{k+1}(t')$.
Moreover, using Theorem \ref{StasheffSmooth} and after fixing some metric on $Gr_{k+1}$, we can identify a neighborhood of $Gr_{k+1}(t')\subset Gr_{k+1}(t)$ by an $\epsilon$-disk bundle. Namely, for $\epsilon>0$ small enough we denote by 
\begin{equation} \label{normalBundleofStasheffStrata}
    \mathbb{D}_{\epsilon}(Gr_{k+1}(t')):=\{\exp_{a}(tV):a\in Gr_{k+1}(t'), V\in\nu_{a}(Gr_{k+1}(t')), 0\leq t<\epsilon\}
\end{equation}
where $\nu_{a}(Gr_{k+1}(t'))$ is the fibre of the normal bundle of $Gr_{k+1}(t')\subset Gr^+_{k+1}(t)$ at $a$. In particular, $Gr_{k+1}(t')\subset\mathbb{D}_{\epsilon}(Gr_{k+1}(t'))\subset Gr^+_{k+1}(t)$.
\begin{definition}
    Let $t=(T, \iota, rt)\in G_{k+1}$ and $\tilde f\in C^\infty(\widetilde{Gr}^+_{k+1}(t)\times L)$. $\tilde f$ is said to be is said to be a Morse-function associated to $f$, if $\tilde f$ satisfies the same conditions as in Definition \ref{T-depMorse}, with the extra condition that $\tilde f(s, .)\equiv f$ for every $s$ on an edge $e$ of length $l(e)=0$.
\end{definition}
\begin{definition}\label{MorseModuliEdgeZero}
    We denote by $\mathcal M^+(t, \underline{x};\tilde f):=\{(l, u):l\in Gr_{k+1}^+(t), u:T\rightarrow L\}$ where $u$ is as in Definition \ref{ModuliMorseFixedDomain}, with the understanding that $u_e$ is constant on each edge $e$ of length $l(e)=0$.
\end{definition}
\begin{definition}
    For $\delta>0$, $t=(T,\iota,\mathrm{rt})\in G_{k+1}$, and
$\tilde f\in C^\infty(\widetilde{Gr}^+_{k+1}(t)\times L)$, we denote by $\mathcal{M}^{\delta}(t,\underline{x};\tilde f)
:=\left\{(l,u)\in\mathcal{M}^+(t,\underline{x};\tilde f)
\;\middle|\;l(e)<\delta\text{ for some }e\in C_{\mathrm{int}}^{1}(T)
\right\}$. In particular, if $T$ is the unique tree with no interior edges then, $\mathcal M^{\delta}(t, \underline x;\tilde f)\equiv\mathcal M(t, \underline x;\tilde f)$.
\end{definition}
\begin{lemma} \label{ExtensionMorseTraj}
    There exists $\epsilon, \delta>0$ where $\delta$ depends on $\epsilon$ and the injectivity radius of the metric on $Gr_{k+1}$, such that if $\Tilde{f}\in C^{\infty}(\widetilde{Gr}^+_{k+1}(t')\times L)$ so that $\mathcal{M}^+(t', \underline{x}; \Tilde{f})$ is a smooth manifold with corners then, $\mathcal{M}^{\delta}(t, \underline{x};\Tilde{f}_t)$ is also a smooth manifold corners. With the understanding that $t'$ is formed from $t$ after shrinking the internal edge of length less than $\delta$ and $\tilde f_t(s, .)\equiv f(.)$ for every $s$ on such edge. Moreover, we have a diffeomorphism $\mathcal{M}^+(t',\underline{x}; \Tilde{f})\times\mathbb{D}_{\epsilon}(Gr_{k+1}(t'))\cong\mathcal{M}^{\delta}(t, \underline{x};\Tilde{f}_t)$.
\end{lemma}
\begin{proof}
    The result follows from the implicit function theorem and the proof of Lemma \ref{genericFixedComb} and so we follow the same notation as above.\\
    In order to prove the first statement, let $\epsilon>0$ and $\Tilde{f}\in C^{\infty}(\widetilde{Gr}_{k+1}(t')\times L)$ be as in Lemma \ref{genericFixedComb}. 
    We note that by the choice of $O_{i}$ as above, namely $O_{i}$ is disjoint from interior
    vertices and only intersects the edge adjacent to $v_{i}$ in its minimal path to $v_{1}$.
    It follows that, $Exp_{\Tilde{f}}$ as given in Equation \ref{Exponentiation}, is independent of the normal variable of 
    $\mathbb{D}_{\epsilon}(Gr_{k+1}(t'))$. In other words, it is independent of the Stasheff 
    deformation parameter in the contracted edge in $t$ forming $t'$. Thus we can extend it
    smoothly in that direction to a smooth map  $\widetilde{Exp}_{\Tilde{f}}$ on
    $L\times\mathbb{D}_{\epsilon}(Gr_{k+1}(t'))$ by adding an extra component corresponding to the gradient vectorfield of $f$ in the normal direction. As $\widetilde{Exp}_{\Tilde{f}}$ is a
    submersion on $L\times \widetilde{Gr}_{k+1}(t')$, it follows that for $\epsilon_1>0$ small enough, $\widetilde{Exp}_{\Tilde{f}}$ is a submersion on $L\times\mathbb{D}_{\epsilon}(Gr_{k+1}(t'))$. Therefore, using the proof of Lemma \ref{genericFixedComb} and after identifying $\mathbb{D_{\epsilon}}(Gr_{k+1}(t'))\subset Gr_{k+1}(t)$ by metric ribbon trees of the same combinatorial type as $t$ with one of their interior edges $e$ having $0\leq l(e)<\delta_1$, where $\delta_1>0$ only depends on $\epsilon_1$ and injectivity radius of the metric on $Gr_{k+1}(t)$; it follows that $\mathcal{M}^{\delta_1}(t, \underline{x}; \Tilde{f}_t)$ is also a smooth manifold with corners, with the understanding that $\tilde f_t\equiv f$ on the fibres of $\mathbb{D}_{\epsilon_1}(Gr_{k+1}(t')) $.\\

    To prove the second statement, we note that $Gr_{k+1}(t')\subset Gr_{k+1}^{+}(t)$ and compare $Exp_{\Tilde{f}}$ and $Exp_{\Tilde{f}_t}$, where $\Tilde{f}\in C^{\infty}(\widetilde{Gr}_{k+1}(t')\times L)$ is as in Lemma \ref{genericFixedComb}. In order to do so, we extend continuously $Exp_{\Tilde{f}}$ to $Exp_{\Tilde{f}}^{+}:L\times Gr_{k+1}^{+}(t)\rightarrow L^{m'}$.
    As $t>t'$ it follows that, $m'\leq m$ where $m$ and $m'$ are the number of interior vertices, $v_{1}, \dots, v_{m}; v_{1}', \dots, v_{m'}'$ in $t, t'$ respectively adjacent to an external edge. Moreover it follows that we have a map $c:\{1, \dots, m\}\rightarrow\{1, \dots, m'\}$ given by $c(i)=j$ if and only if $v_{j}'$ is obtained from $v_{i}$. Similarly, we define an embedding $e_{t', t}: L^{m'}\rightarrow L^{m}$ given by $e_{t' t}(p_{j})_{j=1, \dots, m'}=(q_{i})_{i=1, \dots, m}$ where $p_{j}=q_{i}$ if and only if $c(i)=j$. Then by construction we have 
    \begin{equation*}
        Exp^+_{\Tilde{f}}|_{L\times \widetilde{Gr}_{k+1}(t')}=e_{t', t}\circ Exp_{\Tilde{f}_t}
    \end{equation*}
    Once again by our choice of $O_{i}$, namely each is supported away from any interior vertex and only intersects the edge adjacent to $v_{i}$ in its minimal path to $v_{1}$, and since $e_{t', t}$ is an embedding and $Exp_{\Tilde{f}}$ is transverse to the embedding of $\prod_{i=1}^{m'} \cap_{j=1}^{k'_{i}}\Tilde{W}^{-}_{x_{a}}(e^{j}_{i})$; it follows by openness of the transversality condition, that there exists an $\epsilon_2>0$ small enough such that $Exp^+_{\Tilde{f}}|_{L\times \mathbb{D}_{\epsilon_2}(Gr_{k+1}(t'))}$ is transverse to the embedding of $\prod_{i=1}^{m} \cap_{j=1}^{k_{i}}\Tilde{W}^{-}_{x_{a}}(e^{j}_{i})$. As above, after identifying $\mathbb{D}_{\epsilon_2}(Gr_{k+1}(t'))\subset Gr_{k+1}(t)$ by metric ribbon trees with one of there interior edges $e$ has $0\leq l(e)<\delta_2$, it follows by the implicit function theorem that we have a diffeomorphism $\mathcal{M}^+(t',\underline{x}; \Tilde{f})\times\mathbb{D}_{\epsilon_2}(Gr_{k+1}(t'))\cong\mathcal{M}^{\delta_2}(t, \underline{x};\Tilde{f}_t)$.\\
    Now the result follows after taking $0<\epsilon<\min(\epsilon_1, \epsilon_2)$.
\end{proof}

\begin{definition}\cite{charest2015floer}\label{domDepMorse}
    A domain-dependent Morse function associated to $f$ is a smooth map $\Tilde{f}\in C^\infty(\widetilde{Gr}_{k+1}\times L)$ such that for any $(t, l, s)\in\widetilde{Gr}_{k+1}$ we have $\tilde f_t\in C^\infty (\widetilde{Gr}_{k+1}(t)\times L)$ satisfying the conditions of Definition \ref{T-depMorse}.
\end{definition}
\begin{definition}
    For $\tilde f\in C^\infty(\widetilde{Gr}_{k+1}\times L)$, we denote by $\mathcal{M}(\underline{x}, \tilde f)$ the set of all triples $(t, l, u)$ where:
    \begin{enumerate}
        \item $t\in G_{k+1}$.
        \item $l\in Gr_{k+1}(t)$.
        \item $u\in \mathcal{M}^+(t, \underline{x}; \tilde f_t)$ as in Definition \ref{MorseModuliEdgeZero}.
    \end{enumerate}
\end{definition}
\begin{theorem} \label{ModuliSpaceofMorseAsFamily}
    For a generic choice of a domain-dependent Morse function $\tilde f$ associated to $f$, $\mathcal{M}(\underline{x}; \tilde f)$ is a smooth $Gr_{k+1}$-family.
\end{theorem}
\begin{proof}
    This is an easier version of the proof of Theorem \ref{ThickeningAsFamily} and analogously we argue by induction over the cellular-decomposition of $Gr_{k+1}=\cup_{G_{k+1}}Gr_{k+1}(t)$ given by the combinatorial-type of the trees.\\
    For the base case, consider the unique tree $t_0$ with one internal vertex and $(k+1)$-leaves. Then, by Lemma \ref{genericFixedComb}, we get the desired claim.\\
    As for the inductive step, assume that we have constructed $\tilde f:\widetilde{Gr}_{k+1}(t)\rightarrow C^\infty(L)$ such that $\mathcal{M}(t', \underline{x}; \tilde f_{t'})$ is a smooth manifold with corners for every $t'<t$. Then, by Lemma \ref{ExtensionMorseTraj}, Whitney Extension Theorem and Sard-Smale, it follows that for a generic choice of $\tilde f\in C^\infty(\widetilde{Gr}_{k+1}\times L)$ as in Definition \ref{domDepMorse}, $\mathcal{M}(t, \underline{x}; \tilde f_t)$ is a smooth manifold with corners for every $t\in Gr_{k+1}$.\\
    After arguing as in the proof of Theorem \ref{ThickeningAsFamily}, the result follows from the Implicit Function Theorem.
\end{proof}
To finish this subsection, we observe that, after using standard gluing results of Morse trajectories, as in the proof of Theorem \ref{NormalComChartMorse} for instance, and using the same orientation convention as in Appendix $A$ of \cite{mescher2016perturbed}, the following result follows.
\begin{corollary}[Theorem $5.25$ \cite{mescher2016perturbed}]\label{MorseA-infinity}
    The Morse complex $CM(f; \Z)$ admits a structure of an $A_\infty$-algebra.
\end{corollary}
\subsection{Floer-Morse Trajectories}
In this subsection we fix $(X, \w)$ a closed symplectic manifold together with an $\w$-compatible almost complex structure $J$ and $L\subset X$ a closed, connected, relatively-spin Lagrangian submanifold. In addition, we fix $(f, g)$ a Morse-Smale pair on $L$.
\subsubsection{Thickening}
In order to define what we mean by a Floer-Morse trajectory in $X$, we need to introduce more notation on ribbon trees. Let $(t, l)=(T, \iota, rt; l)\in Gr_{k+1}$ be a rooted metric ribbon tree. Using the ribbon structure, we have an induced orientation on the edges of $T$.
\begin{definition}
    We define $\operatorname{head}, \operatorname{tail}:C^{1}(T)\rightarrow C^{0}(T)$ to be:
    \begin{enumerate}
        \item $\operatorname{head}(e):=v_{+}$ where $v_{+}$ is the vertex adjacent to $e$ where the orientation on $e$ is pointing outward from $v_{+}$.
        \item $\operatorname{tail}(e):=v_{-}$ where $v_{-}$ is the vertex adjacent to $e$ where the orientation on $e$ is pointing inward from $v_{-}$.
    \end{enumerate}
\end{definition}
\begin{definition}[Domain of Floer-Morse Trajectories]
    We denote by $D_{t}:=\bigsqcup_{C^{0}_{int}(T)}\Sigma_{v}\sqcup\bigsqcup_{C^{1}(T)}I_{e}/\sim$ the compact, Hausdorff, connected space, where:
    \begin{enumerate}
        \item For each $v\in C^{0}_{int}(T)$, $(\Sigma_{v}, \partial\Sigma_{v}; z_{1}, \dots, z_{|v|})$ is a bordered, nodal, genus-zero Riemann surface with cyclically-ordered boundary marked points,
        \item For each $e\in C^{1}(T)$, $I_{e}$ is given as in Definition \ref{coordinates},
    \end{enumerate}
    which we construct as follows: 
    using the ribbon structure on $t$ we have a cyclic-ordering of all edges $e_{v}^{1}, \dots, e_{v}^{|v|}$ adjacent to $v\in C^0_{int}(T)$ then, we identify:
    \begin{enumerate}
        \item $z_{i}\sim l(e_{v}^{i})$ if $e_{v}^{i}\in C^{1}_{int}(T)$ and $\operatorname{head}(e_{v}^{i})=v$.
        \item $z_{i}\sim 0$, otherwise.
    \end{enumerate}
\end{definition}
\begin{definition}[Floer-Morse Trajectories of Fixed Combinatorial Type]
    Let $(\beta_v)_{v\in C^0_{int}(T)}$ be a sequence of relative-spherical classes, $\underline{x}=\{x_0, \dots, x_k\}$ a sequence of critical points of $f$ and $\tilde f_t\in C^{\infty}(\widetilde{Gr}_{k+1}(t)\times L)$ be as in Definition \ref{T-depMorse}. A Floer-Morse Trajectory of type $(t, (\beta_v)_{v\in C^0_{int}(T)})$ associated to $\underline{x}$ is a continuous map $u:D_t\rightarrow X$ such that:
    \begin{enumerate}
        \item For each $v\in C^0_{int}(T)$, $u_v:(\Sigma_v, \partial\Sigma_v)\rightarrow (X, L)$ is a $J$-holomorphic map with boundary condition on $L$, representing the class $\beta_v$.
        \item For each $e\in C^1(T)$, $u_e:I_e\rightarrow L$ solves Equation \ref{MorseTraj} so that
        \begin{enumerate}
            \item If $\operatorname{head}(e)=rt$ then $u_e(\infty)=x_0$.
            \item If $\operatorname{tail}(e)\in C^0_{ext}(T)$ then $u_e(-\infty)=x_i$ for some $i=1,\dots, k$.
        \end{enumerate}
    \end{enumerate}
    We denote by $\mathcal{M}(t, \underline{x}, (\beta_v)_{v\in C^0_{int}(T)}; \tilde f_t)\equiv\mathcal{M}(t, \underline{x}, (\beta_v)_{v\in C^0_{int}(T)}; \tilde f_t, g, J)$ the moduli space of all tuples $(l, u)$ where $l\in Gr_{k+1}(t)$ and $u$ is a Floer-Morse trajectories of combinatorial type $(t, (\beta_v)_{v\in C^0_{int}(T)})$ associated to $\underline{x}$.
\end{definition}
\begin{remark}
    We reiterate that from the continuity condition on $u:D_t\rightarrow X$, it follows that the incident condition on all internal Morse trajectories are uniquely determined. Thus, to avoid notational clutter, we do not specify them.
\end{remark}
Before we proceed, we need to fix some more notations. For $\tilde f\in C^\infty(\widetilde{Gr}_{k+1}\times L)$ as in Definition \ref{domDepMorse}, we denote by:
\begin{enumerate}
    \item $\psi_\tau$ the time-$\tau$ flow of $-\nabla^g\tilde f_t(s, .)$.
    \item $\tilde W_{x_0}^+=\{x\in L: \lim_{\tau\rightarrow\infty}\psi_\tau(x)=x_0\}$.
    \item $\tilde W_{x_i}^-=\{x\in L: \lim_{\tau\rightarrow-\infty}\psi_\tau(x)=x_i\}$.
    \item $Q$ the image of the embedding
    \begin{equation*}
        (L\setminus\operatorname{Crit}(f))\times\widetilde{Gr}_{k+1}(t)\times[0,\infty)\hookrightarrow L\times L
    \end{equation*}
    given by given by $(x, (s, t, l);\tau)\mapsto (x, \psi_{\tau}(x))$ and by $Q_{l(e)}$ the image of the above embedding for $\tau=l(e)$.
\end{enumerate}
\begin{proposition} \label{ThickeningFixedDomain}
    For a generic choice of $\tilde f_t\in C^\infty(\widetilde{Gr}_{k+1}(t)\times L)$ as in Definition \ref{T-depMorse}, $\mathcal{M}(t, \underline{x}, (\beta_v)_{v\in C^0_{int}(T)}; \tilde f_t)$ admits a Global Kuranishi chart representation.
\end{proposition}
\begin{proof}
    The result follows from Theorem \ref{ChartsDisks} and Lemma \ref{genericFixedComb}. 
    We start by numbering the internal vertices of the underlying tree $T$ of $t$ as follows: \\
    For each $i=1,\dots, k$ consider the unique minimal path joining $rt$ the root to $v_i\in C^0_{ext}(T)$. Following such path, we number the internal vertices accordingly, with the understanding that if $v\in C^0_{int}(T)$ is contained in the path from $rt$ to $v_i$ and from $rt$ to $v_j$ with $i<j$, then we number $v$ using the path joining $rt$ to $v_i$. Such ordering provides us with an ordering of the relative-spherical classes, say $\beta_1, \dots,\beta_{|C^0_{int}(T)|}$.\\
    
     Now for each $\beta_i$ consider the associated Global Kuranishi chart $(G_i, \T_i, \E_i, s_i)$ as in Theorem \ref{ChartsDisks} and denote by $ev_i:\T_i\rightarrow L^{|v_i|}$ the submersive evaluation map as in Theorem \ref{ChartsDisks}. 
     Let $ev_t:=(ev_1, \dots, ev_{|C^0_{int}(T)|}):\T_1\times\dots\times\T_{|C^0_{int}(T)|}\rightarrow L^{|C^1(T)|+|C^1_{int}(T)|}$ and define $\T:=\{(l, u)\}$ where $l\in Gr_{k+1}(t)$ and $u\in ev^{-1}_t(\tilde W^+_{x_0}\times Q_{l(e_1)}\times\dots\times Q_{l(e_{|C^1_{int}(T)|})} \times \tilde W^-_{x_1}\times\dots\times \tilde W^-_{x_k})$ where $\tilde W^-_{x_i}$ and $\tilde W^+_{x_0}$ are as in the proof of Lemma \ref{genericFixedComb}. By Lemma \ref{genericFixedComb} and the fact that $ev_t$ is a submersion, it follows that for a generic choice of $\tilde f_t$, $\T$ is a smooth manifold with corners. Now set $G\equiv G_{t, (\beta_v)}:=G_1\times\dots\times G_{|C^0_{int}(T)|}$, $\E\equiv\E_{t, (\beta_v)}:=i^*(\E_1\oplus\dots\oplus\E_{|C^0_{int}(T)|})$ and $s\equiv s_{t, (\beta_v)}:=i^*(s_1\oplus\dots\oplus s_{|C^0_{int}(T)|})$ where $i:\T\hookrightarrow\T_1\times\dots\times\T_{|C^0_{int}(T)|}$ is the canonical inclusion. Then $(G, \T, \E, s)$ is a Global Kuranishi chart of $\mathcal{M}(t, \underline{x}, (\beta_v)_{v\in C^0_{int}(T)}; \tilde f_t)$ as desired.
\end{proof}
\begin{remark}
    We note that, the need to use Lemma \ref{genericFixedComb} in the above proof is necessary as we assume neither that $\beta_v\neq 0$ for $v\in C^0_{int}(T)$ nor that the critical points $x_i$ are distinct. 
\end{remark}
\begin{definition}
    For $\beta\in\pi_2(X, L)$, we denote by $\mathcal{M}(t, \underline{x}, \beta; \tilde f_t):=\bigsqcup_{\Sigma\beta_v=\beta} \mathcal{M}(t, \underline{x}, (\beta_v)_{v\in C^0_{int}(T)}; \tilde f_t)$ the moduli space of Floer-Morse trajectories of type $(t, \beta)$ associated to $\underline{x}$.
\end{definition}
Noting that, by Gromov's compactness, the set of $\beta_i\in\pi_2(X, L)$ that be represented by a $J$-holomorphic map with $\w(\beta_i)\leq\w(\beta)$, is finite. We get the following corollary.
\begin{corollary} \label{ThickeningMorseFloerFixedDom}
    For a generic choice of $\tilde f_t\in C^\infty(\widetilde{Gr}_{k+1}(t)\times L)$ as in Definition \ref{T-depMorse}, $\mathcal{M}(t, \underline{x}, \beta; \tilde f_t)$ admits a Global Kuranishi chart representation $(G_{t, \beta}, \T(t, \underline{x}, \beta; \tilde f_t), \E_{t, \beta}, s_{t, \beta})$.
\end{corollary}
 \begin{definition}
     For $\delta>0$, $t=(T, \iota, rt)\in G_{k+1}$, $\beta\in\pi_2(X, L)$ and $\tilde f_t\in C^\infty(\widetilde{Gr}_{k+1}(t)\times L)$, we denote by $\mathcal{M}^\delta(t, \underline{x}, \beta; \tilde f_t)\subset\mathcal{M}(t, \underline{x}, \beta; \tilde f_t)$ the set of all Floer-Morse trajectories of type $(t, \beta)$ associated to $\underline{x}$, where $l(e)<\delta$ for some $e\in C^1_{int}(T)$.
 \end{definition}
 Let $t, t'\in G_{k+1}$ such that $t>t'$ and $t'$ is formed from $t$ by shrinking a single interior edge.
 \begin{proposition} \label{ExtensionofDomDepMorse}
     There exists $\epsilon,\delta>0$ where $\delta$ depends on $\epsilon$ and the injectivity radius of the metric on $Gr_{k+1}$, such that if $\tilde f_{t'}\in C^\infty(\widetilde{Gr}_{k+1}(t')\times L)$ is as in Proposition \ref{ThickeningFixedDomain} then, $\mathcal{M}^\delta(t, \underline{x}, \beta; \tilde f_t)$ admits a Thickening $\T^\delta(t, \underline{x}, \beta; \tilde f_t)$ with the understanding that:
     \begin{enumerate}
         \item $\tilde f_t|_{Gr_{k+1}(t')}\equiv \tilde f_{t'}$.
         \item $\tilde f_t(s, .)\equiv f$ on the shrinked edge $e$ of length $l(e)<\delta$.
     \end{enumerate}
     Moreover, we have a diffeomorphism $\T^\delta(t, \underline{x}, \beta; \tilde f_t)\cong \T(t', \underline{x}, \beta; \tilde f_{t'})\times \mathbb{D}_{\epsilon}(Gr_{k+1}(t'))$, where $\T(t', \underline{x}, \beta; \tilde f_{t'})$ is the Thickening of $\mathcal{M}(t', \underline{x}, \beta; \tilde f_{t'})$ as given in Proposition \ref{ThickeningFixedDomain}.
 \end{proposition}
\begin{proof}
   The result follows from Proposition \ref{ThickeningFixedDomain} and Lemma \ref{ExtensionMorseTraj}. Indeed, let $\tilde f_{t'}$ be as in the statement and work on a single component of $\mathcal{M}(t, \underline{x}, \beta; \tilde f_t)$ where $\tilde f_t$ satisfies the same conditions as in the statement. Following the notation in the proof of Proposition \ref{ThickeningFixedDomain}, we consider $ev_t^{-1}(\tilde W^+_{x_0}\times Q^{|C^1_{int}(T')|}\times\{(x, \phi^f_\tau(x)):x\in L, \tau\geq 0\}\times \tilde W^-_{x_1}\times\dots\times \tilde W^-_{x_k})$, where $\phi^f_\tau$ is the time $\tau$-flow of $-\nabla^gf$. By Lemma \ref{ExtensionMorseTraj}, it follows that for $\delta>0$ small enough, $\tilde W^+_{x_0}\times Q^{|C^1_{int}(T')|}\times\{(x, \phi^f_\tau(x)):x\in L, 0\leq\tau<\delta\}\times \tilde W^-_{x_1}\times\dots\times \tilde W^-_{x_k}$ is a smooth submanifold of $L^{|C^1(T)|+|C^1_{int}(T)|}$ intersecting transversally the diagonal. Now the result follows from Theorem \ref{ChartsDisks}, namely from the fact that $ev_t$ is a submersion. \\

   As for the second claim, we need to compare $ev_t$ and $ev_{t'}$. To do so, we make use of the ordering of internal vertices as in the proof of Proposition \ref{ThickeningFixedDomain} and write 
   \begin{equation*}
    ev_t:\T(\beta_1)\times\dots\times\T(\beta_i)\times\T(\beta_{i+1})\times\dots\times\T(\beta_{|C^0_{int}(T)|})\rightarrow L^{|C^1(T)|+|C^1_{int}(T)|} 
   \end{equation*}
   \begin{equation*}
       ev_{t'}:\T(\beta_1)\times\dots\times\T(\beta_i+\beta_{i+1})\times\dots\times\T(\beta_{|C^0_{int}(T')|})\rightarrow L^{|C^1(T')|+|C^1_{int}(T')|} 
   \end{equation*}
   Then, $ev_{t'}=e_{t', t}\circ ev_t$ where $e_{t', t}$ is as in the proof of Lemma \ref{ExtensionMorseTraj}. Now the result follows from the implicit function theorem, provided that $\T(\beta_i)\times_{ev}\T(\beta_{i+1})$ is a smooth submanifold of $\T(\beta_i+\beta_{i+1})$. To this end, we consider the Thickening as given in Theorem \ref{ChartsforDisksCompatible}. 
\end{proof}

\begin{definition}[Moduli Space of Floer-Morse Trajectories] \label{ModuliFloerMorse}
For $\underline{x}:=\{x_0,\dots,x_k\}$ finite set of critical points of $f$, $\beta\in\pi_2(X, L)$ and $\tilde f\in C^\infty(\widetilde{Gr}_{k+1}\times L)$ as in Definition \ref{domDepMorse}. We denote by $\mathcal{M}(\underline{x}, \beta; \tilde f)\equiv\mathcal{M}(\underline{x}, \beta; \tilde f, g, J):=\{((t, l), u)\}$ such that:
\begin{enumerate}
    \item $(t, l)\in Gr_{k+1}$.
    \item $u\in \mathcal{M}(t, \underline{x}, \beta; \tilde f_t)$ with the understanding that $\tilde f_t\equiv \tilde f|_{\widetilde{Gr}_{k+1}(t)}$.
\end{enumerate}
\end{definition}
\begin{theorem} \label{ThickeningAsFamily}
    For a generic choice of $\tilde f\in C^\infty(\widetilde{Gr}_{k+1}\times L)$ as in Definition \ref{domDepMorse}, $\mathcal{M}(\underline{x}, \beta; \tilde f)$ admits a Thickening $\T(\underline{x}, \beta; \tilde f)$ together with a smooth submersion $\T(\underline{x}, \beta; \tilde f)\rightarrow Gr_{k+1}$.
\end{theorem}
\begin{proof}
    We define $\T(\underline{x}, \beta; \tilde f):=\bigsqcup_{t\in G_{k+1}}\T(t, \underline{x}, \beta; \tilde f_t)$, where $\T(t, \underline{x}, \beta; \tilde f_t)$ is as in Corollary \ref{ThickeningMorseFloerFixedDom} and argue by induction on the cell decomposition $Gr_{k+1}=\bigcup_{t\in Gr_{k+1}}Gr_{k+1}(t)$.\\
    For the base case, we consider the unique tree $t_0$ of a single interior vertex and $(k+1)$-exterior edges and use Corollary \ref{ThickeningMorseFloerFixedDom}. 
    As for the inductive step, assume that we have constructed $\tilde f_t\in C^\infty(\widetilde{Gr}_{k+1}(t)\times L)$ so that for every $t>t'$, $\T(t', \underline{x}, \beta; \tilde f_{t'})$ is a smooth manifold with corners and let $t''>t$ where $t$ is formed from $t''$ by shrinking a single interior edge. Now by Proposition \ref{ExtensionofDomDepMorse}, it follows that the Banach space of all functions $f\in C^N(\widetilde{Gr}_{k+1}(t'')\times L)$ so that $f|_{\widetilde{Gr}_{k+1}(t)\times L}\equiv \tilde f_t$ and $\T^\delta(t'', \underline{x}, \beta; f)$ is a smooth manifold with corners, is non-empty. Now the inductive step follows from Corollary \ref{ThickeningMorseFloerFixedDom} and the proof of Lemma \ref{ExtensionMorseTraj}, and after noting that countable intersection of Baire sets is a Baire set.\\ 

    Using the implicit function theorem, it follows from the proof of Proposition \ref{ExtensionofDomDepMorse} that for $((t, l), u)\in\mathcal{T}(\underline{x}, \beta; \tilde f)$, $V_\epsilon\times U$, where $V_\epsilon$ is a neighborhood of zero in the normal fibre of the $Gr_{k+1}(t)\subset Gr_{k+1}$ and $u\in U\subset \mathcal{T}(t, \underline{x}, \beta; \tilde f_t)$ is a smooth coordinate chart, forms an atlas of $\mathcal{T}(\Bar{x}, \beta; \tilde f)$ and hence the first statement. Indeed, let $((t,l),u)\in
\mathcal T(\underline{x},\beta;\widetilde f)$, and let $t'$ denote the limiting combinatorial type obtained by shrinking the edge whose length is less than $\delta$. After choosing $\epsilon>0$ as in Proposition \ref{ExtensionofDomDepMorse}, there is a diffeomorphism $\mathcal T^\delta(t,\underline{x},\beta;\widetilde f_t)\cong\mathcal T(t',\underline{x},\beta;\widetilde f_{t'}) \times\mathbb D_\epsilon\left(Gr_{k+1}(t')\right)$. Choose a smooth coordinate chart $u\in U\subset\mathcal T(t',\underline{x},\beta;\widetilde f_{t'})$ and a local trivialization $V_\epsilon$ of the normal disk bundle $\mathbb D_\epsilon\left(Gr_{k+1}(t')\right)$ near the point corresponding to $(t, l)$. Restricting the above diffeomorphism gives a coordinate chart $U\times V_\epsilon
\longrightarrow\mathcal T^\delta(t,\underline{x},\beta;\widetilde f_t)$. Now it follows by the uniqueness part of parametric implicit-function theorem that the transition maps are smooth as they are induced by the diffeomorphisms of Proposition \ref{ExtensionofDomDepMorse}. \\
With this in mind, it follows that locally the forgetful map $\T(\underline{x}, \beta; \tilde f)\rightarrow Gr_{k+1}$ is a projection map and hence it is a smooth submersion. 
\end{proof}

\subsubsection{Global Kuranishi Charts}
In the above subsection we have constructed a thickening of the moduli space of Floer-Morse trajectories making sure it is a smooth $Gr_{k+1}$-family, provided that the evaluation map on each Thickened moduli space of pseudo-holomorphic curves is a submersion. Moreover, in Theorem \ref{ChartsforDisksCompatible} we have constructed inductively Global Kuranishi charts of the Moduli space of pseudo-holomorphic disks, with the property that if $\beta=\beta_1+\dots+\beta_j$ where each $\beta_i$ can be represented by a pseudo-holomorphic curve and $k+1=k_1+\dots+k_j$ then, the obstruction bundles satisfies $\E_{k_1+1}(\beta_1)\oplus\dots\oplus\E_{k_j+1}(\beta_j)\hookrightarrow\E_{k+1}(\beta)$.\\
Following the same notation as in the above subsection we introduce the following definition.

\begin{definition}[Obstruction Bundle for the moduli of Floer-Morse Trajectories]\label{ObsMorseFloerTraj}
     We define $\E(\underline{x}, \beta)$ to be the pullback of the trivial vectorbundle $\E_{k+1}(\beta)$ over $Gr_{k+1}$ by the projection map $\T(\underline{x}, \beta)\rightarrow Gr_{k+1}$.
\end{definition}
\begin{definition}[Kuranishi Section for the moduli of Floer-Morse Trajectories]
We define $s_{\underline{x}, \beta}:\T(\underline{x}, \beta)\rightarrow\E(\underline{x}, \beta)$ to be the section given by $s((t, l), u):=((t, l), s_{t,(\beta_v)}(u), 0, \dots, 0)$ where $s_{t,(\beta_v)}$ is as in the proof of Proposition \ref{ThickeningFixedDomain} and the number of zeroes is equal to the rank of the normal bundle of $Gr_{k+1}(t)\subseteq Gr_{k+1}$ plus the rank of the quotient bundle $\E_{k_1+1}(\beta_1)\oplus\dots\oplus\E_{k_{|C^0_{int}(T)|}+1}(\beta_{|C^0_{int}(T)|})\hookrightarrow\E_{k+1}(\beta)$.
\end{definition}
\begin{lemma}
    $\E(\underline{x}, \beta)\rightarrow\T(\underline{x}, \beta)$ is an $O(d)$-equivariant vectorbundle and $s_{\underline{x},\beta}:\T(\underline{x}, \beta)\rightarrow\E(\underline{x}, \beta)$ is a well-defined $O(d)$-eqivariant section, where $d=\Omega(\beta)$ and $\Omega$ is as in Lemma \ref{Integral symplectic form}.
\end{lemma}
\begin{proof}
    By construction, $\E(\underline{x}, \beta)\rightarrow\T(\underline{x}, \beta)$ is an $O(d)$-equivariant vectorbundle, where the $O(d)$-action on each strata $\T(t, \underline{x}, (\beta_v); \tilde f_t)$ is given by the induced action from the group embedding $O(d_1)\times\dots\times O(d_{|C^0_{int}(T)|})\hookrightarrow O(d)$ where $t=(T, \iota, rt)$ and $d_i=\Omega(\beta_i)$.\\
    $s_{\underline{x}, \beta}$ is well-defined, as the following diagram commutes:
    \[\begin{tikzcd}
	{\mathcal{E}_{k_1+1}(\beta_1)\oplus\dots\oplus\mathcal{E}_{k_{|C^0_{int}(T)|}+1}(\beta_{|C^0_{int}(T)|})} && {\mathcal{E}_{k+1}(\beta)} \\
	\\
	{Gr_{k+1}(t)} && {Gr_{k+1}}
	\arrow[hook, from=1-1, to=1-3]
	\arrow[from=1-1, to=3-1]
	\arrow[from=1-3, to=3-3]
	\arrow[hook, from=3-1, to=3-3]
\end{tikzcd}\]
To show that $s_{\underline{x}, \beta}$ is continuous, we take $((t, l), u)\in\T(\underline{x}, \beta)$ and neighborhood $((t, l), u)\in V_\epsilon\times U\subset\T(\underline{x}, \beta)$ as in the proof Theorem \ref{ThickeningAsFamily}. On such neighborhood, $s_{\underline{x}, \beta}$ is continuous as it is the pullback of a sum of continuous sections.
\end{proof}
\begin{theorem} \label{ThickeningofTheMainComponent}
    The tuple $(O(d), \T(\underline{x}, \beta), \E(\underline{x}, \beta), s_{\underline{x}, \beta})$ is stably normally complex Global Kuranishi chart of $\mathcal{M}(\underline{x}, \beta;\tilde f)$.
\end{theorem}
\begin{proof}
    At this stage and having either Lemma \ref{IsoOrbiSpace} or Lemma \ref{IsoOrbiSpaceII} in mind, it suffice to show that such chart $(O(d), \T(\underline x, \beta), \E(\underline{x}, \beta), s_{\underline{x}, \beta})$ is normally complex, possibly after further stabilization.\\
    Let $((t, l), u)\in\T(\underline{x}, \beta)$ and suppose that its isotropy subgroup is non-trivial. Observe that, as all Morse trajectories have prescribed domains with initial values and all pseudo-holomorphic disks have at least one marked point, it follows that all the non-trivial representations of $T_{((t,l), u)}\T(\underline x, \beta)$ are from the sphere components. That is, the non-trivial representations are of the form $\oplus_i\ker(D_{u_i^{sphere}}\bar\partial_J+K)$, where the sum is taken over all sphere components of $u$ and $K$ is a compact operator. On the other hand, following the proof of Theorem \ref{ChartsDisks}, it follows that there exists $R$ a complex representations of the underlying isotropy group, so that $\oplus_i(\ker(D_{u_i^{sphere}}\bar\partial_J+K)\oplus R)$ is a complex representation of $O(d)$. Similarly as in the proof of Theorem \ref{ChartsDisks}, after fixing an $O(d)$-equivariant connection on $\T(\underline{x}, \beta)$ and a non-zero vector $v\in T_{(t, l)}Gr_{k+1}$, the result follows after taking parallel transport of $R$ along the horizontal lift of $v$.   
\end{proof}

\subsubsection{Compactification}
For $\underline{x}=\{x_0, \dots, x_k\}$ finite set of critical points of $f$ and $\beta\in\pi_2(X, L)$ the $\textit{virtual}$ dimension of the moduli space $\mathcal{M}(\underline{x}, \beta)$ is given by
\begin{equation}
    \mu(\underline{x}, \beta)=\mu(x_0)-\Sigma_{i=1}^k\mu(x_i)+\mu(\beta)-nk+(k-2)
\end{equation}
where $n=\dim_\R L$ and $\mu(x_j), \mu(\beta)$ are the Morse index and Maslov index of $x_j, \beta$ respectively. \\
We are interested in the cases when $\mu(\underline{x}, \beta)=0,1$.
In the case when $\mu(\underline{x}, \beta)=1$, Theorem \ref{ThickeningofTheMainComponent} provides us with a Global Kuranishi chart $(G, \T, \E, s)$ where $s^{-1}(0)/G$ is not compact. Indeed, for a sequence $((t_n,l_n),u_n)\in s^{-1}(0)/G$ the failure of finding a convergent subsequence in $s^{-1}(0)/G$ is due to the possibility of having $l_n(e_n)\rightarrow\infty$ as $n\rightarrow\infty$, where $e_n\in C^1_{int}(T_n)$ and $T_n$ is the underlying tree of $t_n$. Indeed, in the case when $\mu(\underline{x}, \beta)=1$, a compactification $\M(\underline{x}, \beta)$ of $\mathcal{M}(\underline{x}, \beta)$ is given by $\textit{adding}$ all possible ways where $(t_n,l_n)\in Gr_{k+1}$ $\textit{breaks}$ once. 
\begin{definition}\label{Codimension1boundaryofFloerMorse}
    For $\underline{x}$ finite set critical points of $f$ and $\beta\in\pi_2(X, L)$ we denote by $\M(\underline{x}, \beta)$ to be the topological space, for the sequential topology, definied inductively by $\M(\underline{x}, \beta):=\mathcal{M}(\underline{x}, \beta)\cup\partial\mathcal{M}(\underline{x}, \beta)$, where
\[
\partial \mathcal M(\underline x, \beta) = \bigcup \M(\underline{x}_1, \beta_1)\times \M(\underline{x}_2, \beta_2)
\]
where the union is taken over $\beta_1+\beta_2=\beta$ so that $\beta_i$ can be represented by a pseudo-holomorphic curve and $\underline{x}_1\cup \underline{x}_2=\underline x \cup\{y\}$ where $y\in\mathrm{Crit}(f)$ so that $\mu(\underline{x}_1, \beta_1)+\mu(\underline{x}_2, \beta_2)=\mu(\underline{x}, \beta)-1$.
\end{definition}


\begin{theorem}\label{NormalComChartMorse}
    Let $\underline{x}=\{x_0, \dots,x_k\}$ be a finite set of critical points of $f$ and $\beta\in\pi_2(X, L)$ be a relative-spherical class. Then for a generic choice of $\tilde f:\overline{Gr}_{k+1}\rightarrow C^\infty(L)$ of smooth map as in Remark \ref{smoothOnCompactified} with the property that $\tilde f|_{Gr_{k+1}}\in C^\infty(\widetilde{Gr}_{k+1}\times L)$ as in Definition \ref{domDepMorse}:
    \begin{enumerate}
        \item The moduli space $\M(\underline{x}, \beta)$ admits an oriented normally complex Global Kuranishi chart $(G_\beta, \T(\underline{x}, \beta), \E_{\underline x, \beta}, s_{\underline{x}, \beta})$, together with a smooth submersion $\T(\underline{x}, \beta)\rightarrow \overline{Gr}_{k+1}$.
        \item Each factor of $\partial\mathcal M(\underline{x}, \beta)$ as in Definition \ref{Codimension1boundaryofFloerMorse} admits an oriented normally complex Global Kuranishi chart $(G_{\beta_i}, \T(\underline{x}_i, \beta_i), \E_{\underline{x}_i, \beta_i}, s_{\underline{x}_i, \beta_i})$ together with a smooth submersion $\T(\underline x_i, \beta_i)\rightarrow \overline{Gr}_{k_i+1}$ so that:
        \begin{enumerate}
        \item $G_{\beta_1}\times G_{\beta_2}\hookrightarrow G_{\beta}$.
            \item The following diagram commutes:
        \[\begin{tikzcd}
	{\mathcal{E}_{\underline x_1, \beta_1}\oplus\mathcal{E}_{\underline x_2, \beta_2}} && {\mathcal{E}_{\underline x,\beta}} \\
	\\
	{\mathcal{T}(\underline x_1, \beta_1)\times\mathcal{T}(\underline x_2, \beta_2)} && {\mathcal{T}(\underline x, \beta)} \\
	\\
	{\overline{Gr}_{k_1+1}\times\overline{Gr}_{k_2+1}} && {\overline{Gr}_{k+1}}
	\arrow[hook, from=1-1, to=1-3]
	\arrow[from=1-1, to=3-1]
	\arrow[from=1-3, to=3-3]
	\arrow["{s_{\underline x_1, \beta_1}\oplus s_{\underline x_2, \beta_2}}", bend left = 30pt, from=3-1, to=1-1]
	\arrow[hook, from=3-1, to=3-3]
	\arrow[from=3-1, to=5-1]
	\arrow["{s_{\underline x, \beta}}"', bend right = 30pt, from=3-3, to=1-3]
	\arrow[from=3-3, to=5-3]
	\arrow["glue"', from=5-1, to=5-3]
\end{tikzcd}\]
and compatible with the normal complex structure and the orientation.
\end{enumerate}
    \end{enumerate}
\end{theorem}
\begin{remark}
    By $\textit{compatible with the normal structure}$ we mean, that if we denote by $J^N$ the complex structure on the normal bundle of the $G_\beta$-fixed point set of $\T(\underline x, \beta)$ and by $J_i^N$ the complex structure on the $G_{\beta_i}$-fixed point set of $\T(\underline x_i, \beta_i)$, then the pullback of $J^N$ by the second horizontal arrow in the second item of the statement is equal to $J^N_1\oplus J_2^N$.
\end{remark}
\begin{proof}
First we consider the case when $\mu(\underline{x}, \beta)=1$ and argue by upward two-layered induction on the energy $\w(\beta)$ and $k$ the number of leaves (minus $1$).\\
    For the base case we adhere to Theorem \ref{ModuliSpaceofMorseAsFamily} when $\beta=0$ or Theorem \ref{ThickeningAsFamily} when $\beta$ is the minimal non-zero class represented by a pseudo-holomorphic disk.\\
    Now assume that $\beta_1+\beta_2=\beta$ and $k_1+k_2=k+1$, and for $i=1,2$ we have $\tilde f_i\in C^\infty(\widetilde{Gr}_{k_i+1}\times L)$ so that $\mathcal{M}(\underline{x}_i, \beta_i; \tilde f_i)$ is given a Global Kuranishi chart as in Theorem \ref{ThickeningofTheMainComponent}, with the understanding that $\underline{x}_1\cup\underline{x}_2=\underline{x}\cup\{y\}$ for some $y\in\mathrm{Crit}(f)$ and that $\beta_i$ can be trivial. Also fix one of the possible $k_2$-gluing maps $Gr_{k_1 +1}\times Gr_{k_2 +1}\rightarrow \overline{Gr}_{k+1}$, as the condition $\mu(\underline{x}, \beta)=1$ allows us to only deal with the case of one $\textit{breaking}$. Since each $\tilde f_i$ satisfies the conditions of Definition \ref{domDepMorse}, in particular $\tilde f_i\equiv f$ in a neighborhood of each vertex, it follows by Whitney's Extension Theorem \cite{whitney1934analytic} that we can define $\tilde f:\widetilde{Gr}_{k+1}\times L\rightarrow\R$ so that the pullback of $\tilde f$ by the gluing map and after restricting to the $Gr_{k_i +1}$-factor, agrees with $\tilde f_i$. This construction defines a closed non-empty subset of the Banach space $C^N(\widetilde{Gr}_{k+1}\times L)$. Now using Theorem \ref{ThickeningAsFamily} and Theorem \ref{ThickeningofTheMainComponent} and after applying the above argument $k_2$-times, it follows that for a generic choice of $\tilde f\in C^\infty(\widetilde{Gr}_{k+1}\times L)$, we have the first item and the first part of the second item in the above statement.\\
    Now fix $\tilde f$ as above and suppose that $\mu(\underline{x}_1, \beta_1)+\mu(\underline{x}_2, \beta_2)=\mu(\underline{x}, \beta)-1=0$ (we still need to consider the main component of moduli spaces of negative expected dimension). We construct the map $\T(\underline{x}_1, \beta_1)\times\T(\underline{x}_2, \beta_2)\hookrightarrow\T(\underline{x}, \beta)$ corresponding to one of the gluing map $Gr_{k_1+1}\times Gr_{k_2+1}\rightarrow \overline{Gr}_{k+1}$ as follows: \\
    For $i=1,2$, let $((t_i, l_i), u_i)\in \T(\underline{x}_i, \beta_i)$ and $y\in\mathrm{Crit}(f)$ corresponding to the root of $t_1$ and the exterior vertex of $t_2$ where we glue the underlying trees. Take a small enough neighborhoods so that on the closure of these neighborhoods, we have $\tilde f_i\equiv f$ for $i=1,2$, with the understanding that $\tilde f_i$ is the pullback $\tilde f$ by the gluing map restricted to $Gr_{k_i+1}$. On such neighborhoods, denote by $\gamma_1:[s_1, \infty)\rightarrow L$ and $\gamma_2:(-\infty, s_2]\rightarrow L$ the corresponding Morse trajectory of each of the $u_1, u_2$ respectively. Now standard gluing results provides us with an embedding $\mathcal M(\gamma_1(s_1), y)\times \mathcal M(y, \gamma_2(s_2))\times [R, \infty)\hookrightarrow\mathcal{M}(\gamma_1(s_1), \gamma_2(s_2))$. In return, if $ U_i\subset\T(\underline{x}_i, \beta_i)$ are neighborhoods of $((t_i, l_i), u_i)$ and from our choice of neighborhoods around each corresponding vertex as above, we get an induced embedding $U_1\times U_2\times [R,\infty)\hookrightarrow \T(\underline{x}, \beta)$ compatible with the associated gluing map $Gr_{k_1+1}\times Gr_{k_2+1}\rightarrow \overline{Gr}_{k+1}$, giving us a commutative diagram:
    \[\begin{tikzcd}
	{\mathcal{T}(\underline{x}_1, \beta_1)\times\mathcal{T}(\underline x_2, \beta_2)} && {\mathcal{T}(\underline x, \beta)} \\
	\\
	{Gr_{k_1+1}\times Gr_{k_2+1}} && {\overline{Gr}_{k+1}}
	\arrow[hook, from=1-1, to=1-3]
	\arrow[from=1-1, to=3-1]
	\arrow[from=1-3, to=3-3]
	\arrow[from=3-1, to=3-3]
\end{tikzcd}\]
Now we further apply Theorem \ref{relativeLashof} on the above diagram, with the understanding that $G_{\beta_i}, G_\beta$ acts trivially on $Gr_{k_i+1}, \overline{Gr}_{k+1}$ respectively, to ensure that the embedding is a smooth embedding compatible with $G_{\beta_1}\times G_{\beta_2}\hookrightarrow G_\beta$.
On the other hand, the commutativity of
\[\begin{tikzcd}
	{\mathcal{E}_{\underline x_1, \beta_1}\oplus\mathcal{E}_{\underline x_2, \beta_2}} && {\mathcal E_{\underline x, \beta}} \\
	\\
	{\mathcal T(\underline x_1, \beta_1)\times\mathcal{T}(\underline x_2, \beta_2)} && {\mathcal T(\underline x, \beta)}
	\arrow[hook, from=1-1, to=1-3]
	\arrow[from=1-1, to=3-1]
	\arrow[from=1-3, to=3-3]
	\arrow[hook, from=3-1, to=3-3]
\end{tikzcd}\]
follows from construction of the obstruction bundles as in Definition \ref{ObsMorseFloerTraj}.\\

As for the compatibility of the normal structure, we once agian argue by induction over $\w(\beta)$ and $k$, and stabilize our constructed Global Kuranishi charts accordingly.\\
The base case, follows from orientability of the moduli space $\mathcal M(\underline x; \tilde f)$ in the case when $\beta=0$. In the case when $\beta$ is the minimal class represented by a pseudo-holomorphic disk, we adhere to Theorem \ref{ThickeningofTheMainComponent}.\\
As for the inductive step, let $((t_i, l_i), u_i)\in \T(x_i, \beta_i)$ and assume that its isotropy group is non-trivial. Thus the non-trivial sub-representations of $T_{((t_i, l_i), u_i)}\T(\underline{x}_i, \beta_i)$ are of the form $\ker(D_{u^{sphere}_i}\bar\partial_J+K_i)$ where $u_i^{sphere}$ denotes the sphere components of $u_i$ and $K_i$ is a compact linear operator. By the same argument as in the first proof of Theorem \ref{ChartsDisks}, we can find a complex $G_i$-representation $R_i$ so that $\ker(D_{u^{sphere}_i}\bar\partial_J+K_i)\oplus R_i$ is a non-trivial complex representation of the underlying isotropy group of $((t_i, l_i), u_i)$. On the other hand, if we denote by $((t, l), u)$ the image of $((t_1, l_1), u_1), ((t_2, l_2), u_2)$ under the above embedding, then the non-trivial sub-representations of $T_{((t, l),u)}\T(\underline{x}, \beta)$ is of the form $\ker(D_{u^{sphere}_1}\bar\partial_J+K_1)\oplus\ker(D_{u^{sphere}_2}\bar\partial_J+K_2) $. Therefore, after using the same construction as in the proof of Lemma \ref{finiteAppScheme}, we can find a $G$-representation $R$ together with an equivariant injective linear map $R_1\oplus R_2\hookrightarrow R$ respecting the diagonal embedding $G_{\beta_1}\times G_{\beta_2}\hookrightarrow G_\beta$, it follows that, $\ker(D_{u^{sphere}_1}\bar\partial_J+K_1)\oplus\ker(D_{u^{sphere}_2}\bar\partial_J+K_2) \oplus R$ is a non-trivial complex representation of the underlying isotropy group of $((t, l), u)$. Finally, after fixing an $O(d)$-connection on $\T(\underline{x}, \beta)$, the result follows after taking parallel transport of $R$ along horizontal lifts of vectorfields in the normal direction of $Gr_{k_1+1}\times Gr_{k_2+1}\rightarrow \overline{Gr}_{k+1}$.\\

As for the orientation, we start by fixing a triangulation of $L$ and after possibly taking subdivisions, we can assume that we have a one-to-one correspondence between $\mathrm{Crit}(f)$ and the barycenters of such triangulation. We also assume that the unstable manifolds of $f$ are the interiors of such triangulation. With this mind, we first argue by induction over the cell-structure of $Gr_{k+1}$ as follows:\\
Let $t_0$ be the unique tree of a single interior vertex and $(k+1)$-leaves, then with the preparation above, we can use \cite{fukaya2010lagrangian} to give an orientation on $\T(t_0, \underline{x}, \beta)$. Now assume that we have oriented all $\T(t,, \underline{x}, \beta)$ and let $t'>t$ with the understanding that $t$ is formed from $t'$ after shrinking a single interior edge. Using Proposition \ref{ExtensionofDomDepMorse} and Remark $4.24$ of \cite{charest2015floer}, we orient $\T(t', \underline{x},\beta)$. Lastly, we argue by induction on energy and the number of leaves, and following Proposition $8.1.14$ of \cite{fukaya2010lagrangian} to get our desired result.\\

To finish the proof, we argue by induction over the number of breakings. The base case is precisely Theorem \ref{ThickeningofTheMainComponent}. While the inductive step is an analogue of the above argument with the right change of notation. 
\end{proof}

\section{Fukaya Algebra over $\Z$}
In this section we prove our main result, namely Theorem \ref{MainResult}. We start by recalling the necessary definitions needed as found in \cite{Fukaya2009}.
We denote by $\Pi(L):=\pi_2(X, L)/\sim$ where $\alpha\sim\beta$ if and only if $\w(\alpha)=\w(\beta)$ and $\mu(\alpha)=\mu(\beta)$ where $\mu:\pi_2(X, L)\rightarrow 2\Z$ is the Maslov index of $L$. \\

\subsection{Filtered $A_{n, K}$-algebra}
In this section, we fix $(f, g)$ a Morse-Smale pair on $L$.\\
We give the monoid $\Pi(L)\times\mathbb{N}$ a total preorder $\precsim$ as follows: \\
Given $\beta\in \Pi(L)$ we set 
\begin{equation*}
    ||\beta||:=\sup\{n:\beta=\beta_1+\dots+\beta_n\}+\lceil\omega(\beta)\rceil-1
\end{equation*}
where $\beta_i\neq 0$ and each can be represented by a pseudo-holomorphic curve,
\begin{equation*}
    ||\beta||=-1
\end{equation*}
if $\beta=0\in\Pi(L)$.\\
We say $(\beta_1, k_1)\prec(\beta_2, k_2)$ if and only if 
\begin{equation*}
    ||\beta_1||+k_1<||\beta_2||+k_2
\end{equation*}
or
\begin{equation*}
    ||\beta_1||+k_1=||\beta_2||+k_2 \text{ and } ||\beta_1||<||\beta_2||.
\end{equation*}
We write $(\beta_1, k_1)\sim(\beta_2, k_2)$ if 
\begin{equation*}
     ||\beta_1||+k_1=||\beta_2||+k_2 \text{ and } ||\beta_1||=||\beta_2||.
\end{equation*}
Working with equivalence classes modulo $\sim$, namely 
\begin{equation*}
    [(\beta,k)]:=\{(\beta', k')\in \Pi(L)\times\mathbb N:||\beta'||=||\beta||, k'=k\},
\end{equation*}
$\precsim$ defines a partial order on $\Pi(L)\times\mathbb N/\sim$.\\

For $n,k, K\geq 0$ integers and $\beta\in\Pi(L)$, we will abuse notation and write $(\beta, k)\precsim(n, K)$ if $(||\beta||, k)\precsim(n,K)$. 
We denote by, $\operatorname{Crit}(f)$ the finite set of critical points of $f$ and $CM(f;\Z):=\Z\langle\operatorname{Crit}(f)\rangle$ the Morse complex on $L$ generated by critical points of $f$.
\begin{definition}
    We set $CM(f; \Lambda):=CM(f; \Z)\widehat{\otimes}_{\Z}\Lambda$ where $\widehat{\otimes}$ is the completion of the tensor product with respect to the $T$-adic topology.
\end{definition}
\begin{definition}
    A filtered, gapped $A_{n, K}$-algebra structure on $CM(f; \Lambda)$ is a collection of $\Lambda$-linear maps
    \begin{equation*}
        m_{k}:=\Sigma_{\Pi(L)}m_{k,\beta}T^{\w(\beta)}e^{\frac{\mu(\beta)}{2}}:CM(f; \Lambda)[1]^{\otimes k}\rightarrow :CM(f; \Lambda)[1]
    \end{equation*}
    for every $(\beta, k)\precsim(n, K)$, satisfying:
    \begin{equation*}
    \Sigma_{\beta_1+\beta_2=\beta, k_1+k_2=k+1}\Sigma_i(-1)^{*_i}m_{k_2, \beta_2}(x_0, \dots, x_{i-1}, m_{k_1, \beta_1}(x_i, \dots, x_{i+k_1}), x_{i+k_1+1},\dots, x_{k})=0
\end{equation*}
where $*_i=\Sigma_{j=1}^{i}|x_j|'$ and $|x_j|'=|x_j|-1$.
\end{definition}
\begin{theorem}\label{AnK-structure}
     The Morse complex $CM(f; \Lambda):=CM(f; \Z)\widehat{\otimes}_{\Z}\Lambda$ carries the structure of a curved, gapped, filtered $A_{n, K}$-algebra $\{m_k\}_{k\geq0}$.
\end{theorem}
\begin{proof}
    Given $\beta\in\Pi(L)$ and a pure tensor $x_1\otimes\dots\otimes x_k$ of critical points of $f$ we consider the associated derived orbifold chart $(\mathcal{U}, \mathcal{E}, s)$ of the Global Kuranishi chart of $\mathcal{M}(\underline{x}, \beta)$, where $\underline{x}=\{x_0, x_1, \dots, x_k\}$ as given in Theorem \ref{NormalComChartMorse}. We set 
    \begin{equation*}
        m_{k, \beta}(x_1\otimes\dots\otimes x_k):=\Sigma_{x_0\in \mathrm{Crit}(f)} (-1)^*|s_\epsilon^{-1}(0)\cap\mathcal{U}_{free}|x_0
    \end{equation*} 
     where the sum is taken over all $x_0\in \mathrm{Crit}(f)$ such that $\mu(\underline{x}, \beta)=0$ and $||s-s_\epsilon||_{C^0}<\epsilon$ is an FOP-perturbation of $s$ as in Proposition \ref{AbsFOPapp}.
    We set
    \begin{equation*}
        m_k:=\Sigma_{\beta}m_{k, \beta} T^{\omega(\beta)}e^{\frac{\mu(\beta)}{2}}
    \end{equation*} and extend $\Lambda$-linearily.
    Now the proof follows from the codimension 1 boundary description of  $\mathcal{T}(\underline{x}, \beta)$ for $\mu(\underline x, \beta)=1$, provided that we have chosen our FOP-perturbations consistently. In order to do so, we argue by upward induction over $\omega(\beta)$ the energy and $k$ the number of leaves, using the partial order $\precsim$.\\
    
    The base case follows from Proposition \ref{AbsFOPapp} and Proposition \ref{relExtFOP}. As for the inductive step, we have to consider $\mathcal{T}(\underline{x}, \beta)$ where the expected dimension $\mu(\underline{x}, \beta)=1$ and construct FOP-sections on its associated derived orbifold chart, which we denote by $(\mathcal{U},\mathcal{E}, s)$, that agrees with the one already constructed on its codimesion 1 corners. Indeed, consider a component of the codimesion 1 corner of $\mathcal{M}(\underline{x}, \beta)$, this can be given a Global Kuranishi chart 
    \begin{equation*}
        (G_1\times G_2, \mathcal{T}(\underline{x}_1, \beta_1)\times\mathcal{T}(\underline{x}_2, \beta_2), \E_{\underline{x}_1, \beta_1}\oplus \E_{\underline{x}_2, \beta_2}, s_{\underline x_1, \beta_1}\oplus s_{\underline{x}_2, \beta_2})
    \end{equation*}
   where by the inductive hypothesis we are assuming that each $s_i:=s_{\underline{x}_i, \beta_i}/G_i$ is a strongly transverse FOP-section, for $i=1,2$. By Proposition \ref{sumFOP}, it follows that the associated orbibundle section of $s_1\oplus s_2$ is strongly transverse. Now after group enlargement by $G_1\times G_2\hookrightarrow G$ and stabilization, we get an equivalent Global Kuranishi chart, where we denote its associated derived orbifold chart by $(\mathcal{U}', \mathcal{E}', s')$ together with an open embedding of $(\mathcal{U}', \mathcal{E}')$ into $(\mathcal{U}, \mathcal{E})$. Now the inductive step follows from Proposition \ref{relExtFOP}, Proposition \ref{FOPstabilization} and Theorem \ref{IntegralCycles}.
\end{proof}
\subsection{$A_\infty$-algebra}
Now we are ready to prove the main result of this paper.
\begin{corollary}\label{Ainfty}
    The $A_{n,K}$-algebra on $CM(f, \Lambda)$ as constructed in Theorem \ref{AnK-structure}, can be extended to an $A_\infty$-algebra structure on $CM(f, \Lambda)$.
\end{corollary}
This result follows from Lemma $7.2.69$ of \cite{fukaya2010lagrangian} which we recall here.
\begin{lemma}[7.2.69 \cite{fukaya2010lagrangian}]
  If $C$ has a structure of filtered, $\mathfrak{G}$-gapped $A_{n, K}$-algebra for any non-negative integers $n, K$ such that for each $(n, K)\precsim(n', K')$ the restriction of the filtered, $\mathfrak{G}$-gapped $A_{n', K'}$-algebra structure coincides with the filtered, $\mathfrak{G}$-gapped $A_{n, K}$-algebra structure. Then 
  \begin{equation*}
      m_k=\Sigma_{\beta\in\mathfrak{G}}m_{k, \beta}T^{\w(\beta)}e^\frac{\mu(\beta)}{2}
  \end{equation*}
  is a filtered, $\mathfrak{G}$-gapped $A_\infty$-algebra on $C$.
\end{lemma}
Before we give a proof of the above Corollary, we need the following notions found in \cite{Fukaya2009}.
Let $\Bar{C}$ be a graded group and $1\leq m_1\leq m_2$ be two natural numbers. We set
\begin{equation*}
    B_{m_1\dots m_2}\Bar{C}[1]:=\oplus_{k=1}^{m_2}B_{k}\Bar{C}[1]/\oplus_{k=1}^{m_1-1}B_{k}\Bar{C}[1],
\end{equation*}
where $B_k\Bar{C}[1]:=\Bar{C}[1]^{\otimes k}$.
Consider a collection of degree $1$ group morphisms $\{\Bar{m}_1, \dots, \Bar{m}_{K}\}$
\begin{equation*}
    \Bar{m}_k:B_{k}\Bar{C}[1]\rightarrow\Bar{C}[1]
\end{equation*}
We have an induced morphisms
\begin{equation*}
    \hat{\Bar{m}}_k:B_{1\dots K}\Bar{C}[1]\rightarrow B_{1\dots K}\Bar{C}[1]
\end{equation*}
and define 
\begin{equation*}
    \hat{\Bar{d}}_{1\dots K}:=\Sigma_{k=1}^K \hat{\Bar{m}}_k
\end{equation*}
\begin{definition} [$A_K$-structure]
    The collection  $\{\Bar{m}_1, \dots, \Bar{m}_{K}\}$ defines a structure of an $A_{K}$-algebra on $\Bar{C}$ if
    \begin{equation*}
         \hat{\Bar{d}}_{1\dots K}\circ  \hat{\Bar{d}}_{1\dots K}=0
    \end{equation*}
\end{definition}
\begin{remark}
    If $C$ is an $A_K$-algebra then $C$ is an $A_{K'}$-algebra for every $K'\leq K$.
\end{remark}
Next we define the notion of $A_K$-algebra morphisms.
Let $(\Bar{C}_i, \Bar{m})$ where $i=1, 2$, be two $A_K$-algebras. Consider a collection of group morphisms
    \begin{equation*}
        \Bar{f}_k:B_{1\dots K}\Bar{C}_1[1]\rightarrow\Bar{C}_2[1]
    \end{equation*}
of degree zero. Such collection induces coalgebra morphisms 
\begin{equation*}
    \hat{\Bar{f}}_k:B_{1\dots K}\Bar{C}_1[1]\rightarrow B_{1\dots K}\Bar{C}_2[1]
\end{equation*}
We set
\begin{equation*}
    \hat{\Bar{f}}_{1\dots K}:=\Sigma_{k=1}^{K}\hat{\Bar{f}}_k
\end{equation*}
\begin{definition}[$A_K$-algebra morphisms]
    The collection $\{\Bar{f}_1, \dots, \Bar{f}_K\}$ is said to be an $A_K$-algebra morphism if 
    \begin{equation*}
        \hat{\Bar{f}}_{1\dots K}\circ \hat{\Bar{d}}_{1\dots K}= \hat{\Bar{d}}_{1\dots K}\circ\hat{\Bar{f}}_{1\dots K}
    \end{equation*}
\end{definition}

\begin{remark}
    Composition of $A_K$-algebra morphisms is an $A_K$-algebra morphism.
\end{remark}
In order to define the notion of homotopy of $A_K$-algebras or $A_\infty$-algebras in general, we have to give a model of $[0, 1]\times C=:\mathcal{C}$. 
\begin{definition}[Algebraic definition of $\mathcal{C}$]
    Let $\mathcal{C}$ be an $A_{K}$-algebra over $\Z$ together with $A_{K}$-morphisms:
    \begin{enumerate}
        \item $\overline{\operatorname{incl}}:\Bar{C}\rightarrow\Bar{\mathcal{C}}$
         \item $\overline{\operatorname{Eval_{s}}}:\Bar{\mathcal{C}}\rightarrow\Bar{C}$, where $s=0, 1$.
    \end{enumerate}
    The above data is said to be a model of $[0, 1]\times\Bar{C}$ if and only if 
    \begin{enumerate}
        \item $\overline{\operatorname{incl_{K}}}:B_{1\dots K}\Bar{C}[1]\rightarrow\Bar{\mathcal{C}}[1]$ is zero.
        \item $\overline{\operatorname{Eval_{s}}}\circ\overline{\operatorname{incl}}=id_{\Bar{C}}$ for $s=0, 1$.
        \item  $\overline{\operatorname{incl}}$ and $\overline{\operatorname{Eval_{s}}}$ for $s=0, 1$ are homotopy equivalences.
        \item The morphism $\overline{\operatorname{Eval_{0}}}\oplus\overline{\operatorname{Eval_{1}}}:\Bar{\mathcal{C}}\rightarrow\Bar{C}\oplus\Bar{C}$ is surjective.
    \end{enumerate}
\end{definition}

\begin{proof}[Proof of Corollary \ref{Ainfty}]
Denote by $\hbar_\omega:=\inf\{\omega(\beta)>0:\mathcal{M}(\beta)\neq\emptyset\}$ and note that, by Gromov's Compactness, $\hbar_\omega>0$. Let $0<\epsilon<\hbar_\omega$ and using Lemma \ref{Integral symplectic form}, we choose $\Omega\in\Omega^2(X)$ so that there exists a fixed positive integer $N$ such that $|N\Omega-\omega|<\epsilon$ in the $C^0$-norm. For such $\Omega$, we use Lemma \ref{LineBundleToFrame}, to get a Hermitian line bundle $(\mathcal L, \nabla)\rightarrow X$ to frame our curves. 
     Now by the second statement of Theorem \ref{IntegralCycles}, it follows that for each $(n, K)\precsim (n', K')$ the restriction of the $A_{n', K'}$-algebra on $CM(f;\Lambda)$ of Theorem \ref{AnK-structure} is exactly the $A_{n, K}$-algebra on $CM(f; \Lambda)$ also given by Theorem \ref{AnK-structure}.\\
    In order to apply Lemma $7.2.69$ of \cite{fukaya2010lagrangian}, we yet have to give a model of $[0,1]\times C$ in our setting. To this end and to avoid notational clutter, we set $C:=CM(f; \Lambda)$ and $\Bar{C}:=CM(f; \Z)$ the $\Z$-reduction of $C$. As our ground ring is a torsion ring, the model that we will use for $[0, 1]\times C$ is
\begin{equation*}
    C^{[0, 1]}:=C\oplus C[-1]\oplus C
\end{equation*}
Now define $\mathcal{I}, \mathcal{I}_0, \mathcal{I}_1:C\rightarrow C^{[0, 1]}$ by 
\begin{equation*}
    \mathcal{I}_0(x):=(x, 0, 0)
\end{equation*}
\begin{equation*}
    \mathcal{I}_1(x):= (0, 0, x)
\end{equation*}
\begin{equation*}
    \mathcal{I}(x):= (0, x, 0)
\end{equation*}
and note that both $\mathcal{I}_0, \mathcal{I}_1$ are degree preserving while $\mathcal{I}$ is of degree 1.
Now we extend  $\mathcal{I}, \mathcal{I}_0, \mathcal{I}_1$ linearly to $B_{1\dots K}\Bar{C}[1]:=\oplus_{i=1}^KCM(f; \Z)[1]^{\otimes i}$ and give $C^{[0, 1]}$ the structure of an $A_K$-algebra as follows:
\begin{equation*}
    \mathfrak{m}_1(\mathcal{I}_0(x)):=\mathcal{I}_0(m_1(x))+(-1)^*\mathcal{I}(x)
\end{equation*}
\begin{equation*}
    \mathfrak{m}_1(\mathcal{I}_1(x)):=\mathcal{I}_1(m_1(x))-(-1)^*\mathcal{I}(x)
\end{equation*}
\begin{equation*}
    \mathfrak{m}_1(\mathcal{I}(x)):=\mathcal{I}(m_1(x))
\end{equation*}
where $*=\deg(x)-1$. On the other hand, for $1\leq k,l\leq K$ such that $k+l\leq K$ and $\underline{x}, \underline{z}$ are pure tensor of elements of $\Bar{C}[1]$ of lengths $k, l$ respectively and $y\in\Bar{C}[1]$, we define 
\begin{equation*}
    \mathfrak{m}_{k+l+1}(\mathcal{I}_0(\underline{x}), \mathcal{I}(y), \mathcal{I}(\underline{z})):=(-1)^{\deg z-1}\mathcal{I}(m_{k+l+1}(\underline{x}, y, \underline{z}))
\end{equation*}
\begin{equation*}
    \mathfrak{m}_k(\mathcal{I}_0(\underline{x})):=\mathcal{I}_0(m_k(\underline{x}))
\end{equation*}
\begin{equation*}
    \mathfrak{m}_l(\mathcal{I}_1(\underline{z})):=\mathcal{I}_1(m_l(\underline{z}))
\end{equation*}
Moreover, for $x, y, z\in \Bar{C}[1]$ we set
\begin{equation*}
    \operatorname{Eval_{s=0}}(x, y, z):=x
\end{equation*}
\begin{equation*}
    \operatorname{Eval_{s=1}}(x, y, z):=z
\end{equation*}
\begin{equation*}
    \operatorname{incl}(x):=(x, 0, x)
\end{equation*}
\begin{equation*}
    \mathfrak{m}_0(1):=\operatorname{incl}(m_0(1))
\end{equation*}

    Then as per \cite{Fukaya2009} it follows that $(C^{[0, 1]}, \mathfrak{m}_k)$ is an $A_K$-algebra and $\operatorname{incl}$, $\operatorname{Eval_{s}}$ for $s=0, 1$ are $A_K$-algebra morphisms. Moreover, this data forms a model for $[0, 1]\times C$.

Now consider $\Tilde{f}:(-2, 2)\times L\rightarrow\R$ be given by
\begin{equation*}
    \Tilde{f}(t, x):= f(x)+\sigma(t)
\end{equation*}
where $\sigma:(-2, 2)\rightarrow\R$ is a smooth function with exactly $3$ critical points at $t= -1, 0, 1$ such that the indices at these critical points are $0, 1, 0$ respectively. Observe that, $CM(\Tilde{f}; \Lambda)$ forms a model for $[0, 1]\times C$. Indeed, with the discussion above in mind, it suffice to show that we have a $\Lambda$-module isomorphism between $CM(\Tilde{f}; \Lambda)$ and $C\oplus C[-1]\oplus C$. We define such isomorphism on the critical points of $\Tilde{f}$ by
    \begin{equation*}
        (x, 0, 0) \mapsto (x, -1)
    \end{equation*}
    \begin{equation*}
        (0, x, 0)\mapsto (x, 0)
    \end{equation*}
    \begin{equation*}
        (0, 0, x)\mapsto (x, 1)
    \end{equation*}
    and extend $\Lambda$-linearly.\\
    Now the result follows from Lemma $7.2.69$ of \cite{fukaya2010lagrangian}, after noting that from our choice of $\epsilon$ as above, the $\omega$-valuation on $\Lambda$ denoted by $\nu_\omega$ satisfies $\nu_\omega\equiv N\nu_\Omega$ and hence the $T$-adic topology on $\Lambda$ given the $\Omega$-valuation $\nu_\Omega$ is the same as the one given by the $\omega$-valuation $\nu_\omega$.
\end{proof}

\appendix
\section{Smoothing Theory}
\subsection{Lashof's Smoothing Theorem}
In this section we recall all the necessary definitions and results needed in order to give a smooth structure on our Global Kuranishi charts. We follow closely the discussion and the notations as presented in \cite{abouzaid2021complex}, \cite{bai2022arnold}. \\
Let $X$ be a topological space.
\begin{definition}[Microbundles]
   A rank $n$ microbundle on $X$ is the data of:
    \begin{equation*}
        X\xlongrightarrow{s} E\xlongrightarrow{p} X
    \end{equation*}
    where $E$ is a topological space and $s, p$ are continuous maps satisfying the following conditions:
    \begin{enumerate}
        \item $p\circ s=id_X$.
        \item for every $x\in X$ there exists $U_x\subseteq X$ an open neighborhood of $x$ and $ V_x\subseteq E$ an open neighborhood of $s(x)$ together with a homeomorphism 
        \begin{equation*}
            h_x:U_x\times\mathbb{R}^n\rightarrow V_x
        \end{equation*}
        so that $p\circ h_x=pr_{U_x}$ is the natural projection map to $U_x$ and $h_x|_{U_x\times\{0\}}=s$.
    \end{enumerate}
\end{definition}

\begin{example}\label{MicEx}
    The trivial bundle $X\times\mathbb{R}^n\rightarrow X$ over $X$, together with $s:X\rightarrow X\times\mathbb{R}^n$ given by $s(x)=(x, 0)$ is a rank $n$ microbundle on $X$.\\
    In general, any vectorbundle over a topological space $E\rightarrow X$ is a microbundle after considering $s:X\rightarrow E$ as the zero section.
\end{example}

\begin{definition}[Morphisms of Microbundles] \label{MorMicro}
    Let $X\xlongrightarrow{s_i}E_i\xlongrightarrow{p_i}X$ be microbundles on $X$ for $i=1, 2$, not necessarily of the same rank. A microbundle morphism $\phi: E_1\rightarrow E_2$ is the equivalence class of a continuous map $\phi:U_1\rightarrow E_2$ defined over a neighborhood $U_1\subseteq E_1$ of $s_1(X)$ which commutes with the maps $s_i, p_i$ for $i=1, 2$, respectively.\\
    We say two such morphisms $\phi_1, \phi_2$ are equivalent if there exists an open neighborhood $U\subseteq U_1\cap U_2\subseteq E_1$ containg $s_1(X)$ such that $\phi_1\equiv\phi_2$ on $U$.
\end{definition}

\begin{definition}[Isomorphisms of Microbundles]
    Following the same notation as in Definition \ref{MorMicro}. We say $\phi:E_1\rightarrow E_2$ is an isomorphism of microbundles if $\phi:U_1\rightarrow\phi(U_1)$ is a homeomorphism and $\phi(U_1)\supseteq s_2(X)$.
\end{definition}
\begin{definition}[Isotopy of Isomorphisms of Microbundles]
    Let $\phi_0, \phi_1:E_1\rightarrow E_2$ be two isomorphisms of microbundles
    \begin{equation*}
        X\longrightarrow E_i\longrightarrow X
    \end{equation*}
    on $X$, for $i=1, 2$.\\
    An isotopy between $\phi_0, \phi_1$ is $[0, 1]$-family of microbundle isomorphisms $\Tilde{\phi}_t:E_1\rightarrow E_2$ where $t\in[0, 1]$ so that $\Tilde{\phi}_0$ is equivalent to $\phi_0$ and $\Tilde{\phi}_1$ is equivalent to $\phi_1$.
\end{definition}

\begin{definition}[The Tangent Microbundle]
    The tangent microbundle of $X$ denoted by $T_\mu X$ is the data:
    \begin{equation*}
        X\xlongrightarrow{\Delta}X\times X\xlongrightarrow{pr_1}X
    \end{equation*}
    where $\Delta(x):=(x, x)$ is the diagonal inclusion and $pr_1(x, y):=x$ is the projection onto the first factor.
\end{definition}

\begin{remark}\cite{abouzaid2021complex}
    If $X$ is a topological manifold that admits a smooth structure, then its tangent microbundle $T_\mu X$ is isomorphic as microbundles to its tangent bundle $TX$.
\end{remark}

We also have well-defined notions of direct sum, pullbacks and restrictions of microbundles as in the case of vectorbundles.

\begin{definition}[Direct Sum of Microbundles]
    Given $X\xlongrightarrow{s_i}E_i\xlongrightarrow{p_i}X$ two microbundles over $X$, where $i=1, 2$, denote by $\Delta:X\rightarrow X\times X$ the diagonal inclusion. Let,
    \begin{equation*}
        E_{12}:=(p_1\times p_2)^{-1}(\Delta(X))\subseteq E_1\times E_2
    \end{equation*}
    The direct sum of $E_1, E_2$ denoted by $E_1\oplus E_2$ is the microbundle
    \begin{equation*}
        X\xlongrightarrow{(s_1\times s_2)\circ\Delta}E_{12}\xlongrightarrow{p_1\times p_2}\Delta(X)\cong X
    \end{equation*}
\end{definition}

\begin{definition}[Pullback of Microbundles]\label{pullMicr}
    Given $X\xlongrightarrow{s}E\xlongrightarrow{p}X$ microbundle over $X$ and $f:Y\rightarrow X$ a continuous map, we define the pullback microbundle by $f$ to $Y$, to be 
    \begin{equation*}
        Y\xlongrightarrow{f^*s}f^*E\xlongrightarrow{f^*p}Y
    \end{equation*}
    where
    \begin{equation*}
        f^*E:=\{(e, y)\in E\times Y: p(e)=f(x)\}\subseteq E\times Y
    \end{equation*}
    \begin{equation*}
        f^*s(y):=(s(f(y)), y)
    \end{equation*}
    \begin{equation*}
        f^*p(e, y):=y
    \end{equation*}
\end{definition}

\begin{definition}[Restriction of Microbundles]
    Following the same notation as in Definition \ref{pullMicr} and if in particular $Y\subseteq X$ and $f:Y\hookrightarrow X$ is the canonical inclusion, then we define the restriction of $E$ to $Y$, denoted by $E|_Y$, as $f^*E$.
\end{definition}

For applications, namely to give a smooth structure on our Global Kuranishi charts, we are interested in the $G$-equivariant case, where $G$ is a compact topological group.

\begin{definition} [$G$-microbundle]
    Suppose that $X$ is a $G$-space. A $G$-microbundle on $X$ is a microbundle $X\xlongrightarrow{s}E\xlongrightarrow{p}X$ such that $E$ is a $G$-space and both $s, p$ are $G$-equivariant, namely $g.s(x)=s(g.x)$ and $g.p(e)=p(g.e)$ for any $g\in G$, $x\in X$ and $e\in E$.
\end{definition}

\begin{definition}[Morphisms of $G$-microbundles]
    Let $X\xlongrightarrow{s_i}E_i\xlongrightarrow{p_i}X$ be two $G$-microbundles over $X$, for $i=1, 2$. A morphism of $G$-microbundles $\phi:E_1\rightarrow E_2$ is a $G$-equivariant microbundle morphism. Similarly, an isomorphism of $G$-microbundles is a $G$-equivariant isomorphism of $G$-microbundles.
\end{definition}

We also have the same operations as above in the $G$-equivariant case. We only mention the direct sum case.

\begin{definition}[Direct Sum of $G$-microbundles]
    Let $X\xlongrightarrow{s_i}E_i\xlongrightarrow{p_i}X$ be two $G$-microbundles over $X$, where $i=1, 2$. We define the direct sum of $E_1, E_2$ as $G$-microbundles to be the direct sum $E_1\oplus E_2$ as microbundles where the $G$-action is the diagonal one, namely $g.(e_1, e_2):=(g.e_1, g.e_2)$.
\end{definition}

In what follows, we will be interested in answering whether a given $G$-microbundle over a $G$-space is induced by a $G$-equivariant vectorbundle. Indeed, given $G$-equivariant vectorbundle $E\rightarrow X$, we can associate a $G$-microbundle on $X$, as in Example \ref{MicEx}, namely by 
\begin{equation*}
    X\xlongrightarrow{s=0}E\xlongrightarrow{p} X
\end{equation*}
To this end we introduce the following terminology.
\begin{definition}[$G$-vectorbundle Lift]
    Given  $X\xlongrightarrow{s}E\xlongrightarrow{p} X$ a $G$-microbundle. A $G$-vectorbundle lift of $E$ is a $G$-equivariant vectorbundle $V\rightarrow X$ on $X$, such that the associated $G$-microbundle of $V$ is isomorphic to $E$ as $G$-microbundles.
\end{definition}

The following result of Lashof \cite{lashof2006stable}, gives us sufficient conditions to give a topological manifold together with an action of a compact Lie group, a smooth structure after stabilizing by some $G$-representation.

\begin{theorem} [Theorem $6.8$ \cite{lashof2006stable}]
\label{lashof} 
    Suppose that, $X$ is a topological manifold and $G$ is a compact Lie group acting continuously and locally-linear on $X$. Assume in addition that,
    \begin{enumerate}
        \item There are only finitely many orbit types.
        \item The tangent microbundle $T_\mu X$ admits a $G$-vectorbundle lift.
    \end{enumerate}
    Then there exists a $G$-representation $V$ such that $X\times V$ admits a $G$-equivariant smooth structure.
\end{theorem}
We call the $G$-equivariant smooth structure on $X\times V$, a $\textbf{$G$-smoothing of $X$}$.
\begin{remark}
    In the case when $G$ is the trivial group, the above theorem recovers a result of Milnor \cite{milnor1964microbundles}.
\end{remark}

\begin{definition} [Stable $G$-isotopy of Smooth Structures] \label{stableIsoDefn}
    Let $X$ be a topological manifold together with a compact Lie group $G$ acting continuously on $X$. For $i=1, 2$, let $V_i$ be a $G$-representation such that $X\times V_i$ admits a $G$-equivariant smooth structure. We say that the smooth structures on $X\times V_i$ are stably $G$-isotopic, if there exists $G$-representations $V'_i$ so that $V_1\oplus V'_1\cong V_2\oplus V'_2$ as $G$-representations and there exists an isotopy through $G$-equivariant homeomorphisms between $id_{X\times (V_1\oplus V'_1)}$ and a diffeomorphism on $X\times (V_2\oplus V'_2)$.
\end{definition}

\begin{definition}[Stable $G$-isotopy of Vectorbundle Lifts]
The setting is as in Definition \ref{stableIsoDefn}. For $i=1, 2$, let $V_i$ be a $G$-representation and assume that $T_\mu X\oplus V_i$ (direct sum is taken as $G$-microbundles) has a $G$-vectorbundle lift $E_i$. Such lifts are said to be stably $G$-isotopic, if there exists $G$-representations $V'_i$ so that $V_1\oplus V'_1\cong V_2\oplus V'_2$ as $G$-representations and $E_1\oplus V'_1\cong E_2\oplus V'_2$ as $G$-equivariant vectorbundles.   
\end{definition}

Stable $G$-isotopy of smooth structures and stable $G$-isotopy of vectorbundle lifts are a equivalence relations. In fact, we have the following correspondence between their equivalence classes.

\begin{theorem}[Theorem $6.8$ \cite{lashof2006stable}]
    The setting is as in Theorem \ref{lashof} and we further assume that $X$ satisfies the conditions in \ref{lashof}. Then we have a bijection between equivalence classes of stable $G$-isotopies of smooth structures on $X$ and stable $G$-isotopies of vectorbundle lifts of $T_\mu X$.
\end{theorem}
In fact, what we will be using throughout this paper is the $\textit{relative-version}$ of Lashof's Theorem \ref{lashof}.
\begin{theorem} [Relative Lashof's Theorem, Theorem $B.3$\cite{bai2022arnold}]\label{relativeLashof}
    The setting is as in Theorem \ref{lashof} and in addition, we suppose that $C\subseteq X$ is a $G$-invariant closed subset and $C\subset U$ is a $G$-invariant open neighborhood of $C$. Denote by $E\rightarrow X$ a $G$-vectorbundle lift of the tangent microbundle $T_\mu X$. Furthermore, we assume that
    \begin{enumerate}
    \item $E|_U\rightarrow U$ is a $G$-vectorbundle lift of the tangent microbundle $T_\mu U$.
        \item $V_1$ is a $G$-representation and $U\times V_1$ is equipped with a $G$-equivariant smooth structure.
    \end{enumerate}
    Then there exists $C\subset U'\subseteq U$ a $G$-invariant neighborhood of $C$, a $G$-representation $V$ and $G$-equivariant smooth structure on $X\times V$ such that the product $G$-equivariant smooth structure on $(U'\times V_1)\times V$ is diffeomorphic to the product $G$-equivariant smooth structure on $(U'\times V)\times V_1$.
\end{theorem}
\subsection{Topological Submersions}
Throughout this section, we fix a compact Lie group $G$ and topological manifolds $\T$ and $\F$, together with a surjective continuous map $\pi:\T\rightarrow\F$. We follow closely the discussion as in \cite{abouzaid2021complex}, \cite{bai2022arnold} and start by introducing the following notation. Given any $W\subseteq\T$ and $y\in\F$, we set $W|_y:=\pi^{-1}(y)\cap W$.

\begin{definition}[Product Neighborhood]\label{ProNbhd}
    Let $x\in\T$ and $y:=\pi(x)\in\F$. A product neighborhood of $x$ in $\T$ is an open neighborhood $x\in W\subseteq\T$ together with a homeomorphism $\iota:W\rightarrow W|_y\times\pi(W)$ satisfying the following conditions:
    \begin{enumerate}
        \item $\pi\circ\iota^{-1}$ is the projection map onto $\pi(W)$.
        \item $\iota|_{W|_y}:W|_y\rightarrow W|_y\times\{y\}$ is the identity map.
    \end{enumerate}
\end{definition}

\begin{definition}[Topological Submersion]\label{TopSubmersion}
    The setting is as in Definition \ref{ProNbhd}. $\pi$ is said to be a topological submersion, if every $x\in\T$ admits a product neighborhood. 
\end{definition}

\begin{example}
    In the case when, $\T$ and $\F$ are smooth manifolds and $\pi:\T\rightarrow\F$ is a submersion. Then $\pi$ is a topological submersion.
\end{example}

\begin{lemma}
    Suppose that $\pi:\T\rightarrow\F$ is a topological submersion and let $y\in\F$. Then, $\pi^{-1}(y)$ is a topological manifold.
\end{lemma}
\begin{proof}
    The collection of product neighborhoods of $x\in\pi^{-1}(y)$ forms an atlas on $\pi^{-1}(y)$.
\end{proof}

\begin{definition}[Vertical Tangent Microbundle of a Topological Submersion]
    Suppose that $\pi:\T\rightarrow\F$ is a topological submersion. We define $T^{vt}_\mu \T(\pi)$ the vertical tangent microbundle of $\pi$ by
    \begin{equation*}
        \T\xlongrightarrow{\Delta}\T\times_\F\T\xlongrightarrow{\pi\times\pi}\T
    \end{equation*}
    where $\T\times_\F\T$ is the fibre product given by $\pi$ and $\Delta:\T\rightarrow\T\times\T$ is the diagonal map, and $(\pi\times\pi)(x, x')=\pi(x)\equiv\pi(x')$.
\end{definition}

For applications, we are interested in the $G$-equivariant case. From now on we assume that $\T$ and $\F$ are topological $G$-manifolds and $\pi:\T\rightarrow\F$ is a $G$-equivariant map. Moreover, for $x\in\T$ or $x\in\F$ we denote by $G_x\leq G$ the isotropy subgroup of $x$.\\
We extend all the above definitions naturally to the equivariant case. \\
Indeed, notice that if the pair $W, \iota$ form a product neighborhood of $x\in\T$ and $g\in G$ then, the pair $g.W:=\{g.x':x'\in W\}$ and $g_*\iota(x'):=g.\iota(g^{-1}.x')$ is a product neighborhood of $g.x$, with the understanding that $g$ acts diagonally on $W|_x\times \pi(W)$. With this in mind, we introduce the following definition.

\begin{definition}[$G_x$-invariant Product Neighborhood]
    Given $x\in\T$, a $G_x$-invariant neighborhood of $x$ is a product neighborhood of $x$ given by a pair $x\in W$ and $\iota$ as in Definition \ref{ProNbhd}, such that $g.W=W$ and $g^*\iota\equiv\iota$ for every $g\in G_x$.
\end{definition}

\begin{definition} [$G$-equivariant Topological Submersion]
    $\pi:\T\rightarrow\F$ is said to be a $G$-equivariant topological submersion if every $x\in\T$ admits a $G_x$-invariant product neighborhood.
\end{definition}

In our setting, we are interested in the case where $G$ acts $\textit{locally linear}$ on both $\T$ and $\F$.

\begin{definition}[Locally Linear $G$-action]\label{locallyLinearAction}
    A topological $G$-manifold $M$ is said to be locally linear, if for every $p\in M$ there exists a locally Euclidean neighborhood of $p\in M$ such that, the stabilizer $G_p\leq G$ of $p$ acts linearly on this chart.
\end{definition}

\begin{definition}[Fibre-wise Locally Linear $G$-equivariant Topological Submersion]
    Suppose that $\T$ and $\F$ are locally linear topological $G$-manifolds and $\pi:\T\rightarrow\F$ is a $G$-equivariant topological submersion. $\pi$ is said to be fibre-wise locally linear, if for every $y\in\F$, the stabilizer subgroup $G_y$ acts locally linearly on the topological manifold $\pi^{-1}(y)\subseteq\T$.
\end{definition}

\begin{example}
    Suppose that $p:V\rightarrow B$ is $G$-vectorbundle over topological $G$-manifold $B$. Then, $p$ is fibre-wise locally linear $G$-equivariant topological submersion. 
\end{example}

\subsection{Fibre-wise Smooth Structure}
In this section, we suppose that $G$ is a compact Lie group, $\T$ is a topological $G$-manifold and $\F$ is a smooth $G$-manifold. In addition, we suppose that $\pi:\T\rightarrow\F$ is a $G$-equivariant topological submersion.

\begin{definition} [$C^1_{loc}$-compatible Product Neighborhoods]
    Let $y_i\in\F$ and $W_i, \iota_i$ for $i=1, 2$ be two product neighborhoods of $\pi$. Such neighborhoods are said to be $C^1_{loc}$-compatible, if for every $x\in W_1\cap W_2$, there exists a product neighborhood $W\subseteq\F$ of $x$ and $\iota:W\rightarrow W|_y\times\pi(W)$ where $y=\pi(x)$ and the $\pi(W)$-family of maps
    \begin{equation*}
        \eta_v:W|_y\rightarrow W_i|_{y_i}
    \end{equation*}
    given by $w\mapsto pr_i(\iota((\iota|_{W|_v})^{-1}(w))$ for each $v\in\pi(W)$ are all smooth and vary continuously with respect to the $C^1_{loc}$-topology, where $pr_i$ is the projection map to $W_i$ for $i=1, 2$.
\end{definition}
\begin{definition}[Fibre-wise Smooth $C^1_{loc}$ $G$-bundle] \label{Fibre-wise C1}
    $\pi:\T\rightarrow\F$ is said to be a fibre-wise smooth $C^1_{loc}$ $G$-bundle, if there exists a collection of points $x_i\in\T$ together with a collection of $C^1_{loc}$-compatible $G_{x_i}$-invariant product neighborhoods $W_i, \iota_i$ whose domains cover $\T$.
\end{definition}
In the case of $\pi:\T\rightarrow\F$ is a fibre-wise smooth $C^1_{loc}$ $G$-bundle, we have a well-defined $\textit{vertical tangent bundle}$ of $\pi$, denoted by $T^{vt}\T$.

\begin{definition}[Vertical Tangent Bundle of a Fibre-wise Smooth $C^1_{loc}$ $G$-bundle]
    Using the same notation as in Definition \ref{Fibre-wise C1}, the vertical tangent bundle $T^{vt}\T$ of $\pi$ is the $G$-equivariant vectorbundle whose restriction to $W_i$ of a product neighborhood $W_i, \iota_i$ is the pullback of $T(W|_{y_i})$ to $W_i$ via the projection map to $\pi(W_i)$.
\end{definition}

\begin{lemma}[Lemma $4.29$ \cite{abouzaid2021complex}]
    Suppose that $\pi:\T\rightarrow\F$ is a fibre-wise smooth $C^1_{loc}$ $G$-bundle. Then $T^{vt}\T$ the vertical tangent bundle is a $G$-equivariant lift of $T^{vt}_{\mu}\T$ the vertical tangent microbundle of $\pi$. 
\end{lemma}
\begin{proposition}[Corollary $4.26$ \cite{abouzaid2021complex}]\label{microBundleLift}
    Suppose that $\F$ is locally linear topological $G$-manifold and $\pi:\T\rightarrow\F$ is a fibre-wise locally linear $G$-equivariant topological submersion. Then there exists a morphism of microbundles $P:T_\mu \T\rightarrow T_\mu^{vt}(\T)$ such that $P|_{T_\mu^{vt}(\T)}:T_\mu^{vt}(\T)\rightarrow T_\mu^{vt}(\T)$ is the identity map. Moreover, $P$ induces a $G$-equivariant isomorphism of mictrobundles 
    \begin{equation*}
        P\oplus\tau:T_\mu \T\rightarrow T_\mu^{vt}(\T)\oplus\pi^*(T_\mu\F)
    \end{equation*}
    where $\tau:T_\mu \T\rightarrow\pi^*(T_\mu \F)$ is given by
    \begin{equation*}
        \tau(p, w):=(p, (\pi(p), \pi(w)))\in\pi^*(T_\mu \F).
    \end{equation*}
\end{proposition}
\section{Fukaya-Ono-Parker Perturbations}
Following \cite{bai2022integral}, \cite{bai2022arnold} we use the language of derived orbifold charts. In order to define the notion of an FOP-perturbation we need the notion of normally complex derived orbifold charts and the notion of straightening on derived orbifold charts. We recall all the necessary definitions and results in this section that are used in this paper.\\
We note that all our orbifolds are assumed to be effective orbifolds.
\begin{definition}[Derived Orbifold Chart]\label{DerOrbi}
    A derived orbifold chart is a triple $(\mathcal{U}, \E, s)$ where $\E\rightarrow\mathcal{U}$ is an orbibundle and $s:\mathcal{U}\rightarrow\E$ is a continuous section.
    \begin{enumerate}
        \item $(\mathcal{U}, \E, s)$ is said to be a compact derived orbifold chart if $s^{-1}(0)\subseteq\mathcal{U}$ is compact.
        \item $(\mathcal{U}, \E, s)$ is said to be an oriented derived orbifold chart if both $\mathcal{U}$ and $\E$ are oriented.
    \end{enumerate}
\end{definition}
\begin{remark}
    As all our orbifolds are effective, it follows that, the subset of all points of trivial isotropy forms an open and dense set. 
\end{remark}
\begin{definition}
    A derived orbifold chart $(\mathcal{U}, \E, s)$ is said to be normally complex if the orbibundle $\E\rightarrow\mathcal{U}$ is normally complex as in Definition \ref{normallyComplexOrbibundle}.
\end{definition}
We list the operations on derived orbifold charts that are used in this paper.
\begin{definition}\cite{bai2022arnold}
    \begin{enumerate}
        \item An open embedding from a derived orbifold $(\mathcal{U}, \E, s)$ to $(\mathcal{U}', \E', s')$ consists of an open embedding $\phi:\mathcal{U}\rightarrow\mathcal{U}'$ of orbifolds together with an orbibundle isomorphism $\hat{\phi}$ such that the following diagram commutes:
\[\begin{tikzcd}
	{\mathcal{E}} && {\mathcal{E}'|_{\phi(\mathcal{U})}} \\
	\\
	{\mathcal{U}} && {\phi(\mathcal{U})}
	\arrow["{\hat{\phi}}", from=1-1, to=1-3]
	\arrow[from=1-1, to=3-1]
	\arrow[from=1-3, to=3-3]
	\arrow["s", bend left = 30 pt, from=3-1, to=1-1]
	\arrow["\phi"', from=3-1, to=3-3]
	\arrow["{s'|_{\phi(\mathcal{U})}}"', bend right = 30 pt, from=3-3, to=1-3]
\end{tikzcd}\]
and $(s')^{-1}(0)\subseteq \phi(\mathcal{U})$.
\item For $i=1,\dots,k$, let $(\mathcal{U}_i, \E_i, s_i)$ be a derived orbifold chart. Their product is the derived orbifold chart $(\mathcal{U}_1\times\dots\times\mathcal{U}_k, \pi_1^*\E_1\oplus\dots\oplus\pi_k^*\E_k, s_1\oplus\dots\oplus s_k)$, where $\pi_i:\mathcal{U}_1\times\dots\times\mathcal{U}_k\rightarrow\mathcal{U}_i$ is the natural projection map.
       \item Let $(\mathcal{U}, \E, s)$ be a derived orbifold chart and $\pi_\F:\F\rightarrow\mathcal{U}$ be an orbibundle. The stabilization of $(\mathcal{U}, \E, s)$ by $\F$ is the derived orbifold chart $(\F, \pi^*_\F\E\oplus\pi^*_\F\F, \pi^*_\F s\oplus\tau_\F)$ where $\tau_\F:\F\rightarrow\pi^*_\F\F$ is the tautological section.
    \end{enumerate}
\end{definition}
As in \cite{fukaya2001floer}, \cite{bai2022integral}, \cite{bai2022arnold}, one has to fix a metric and a connection on orbibundles in order to identify tubular neighborhoods of fixed point sets with some disk bundle of their normal bundles in a way compatible with the isotropy stratification. 
\begin{definition}[Straightening Structure]
    A straightening of a derived orbifold chart $(\mathcal{U}, \mathcal{E}, s)$ is a choice of a Riemannian metric $g$ on $\mathcal{U}$ and a connection $\nabla$ on $\mathcal{E}$ satisfying the following conditions:
    \begin{enumerate}
        \item for each chart $(\Gamma, U, E)$ the pullback metric $g_U$ satisfies the following condition. Given a subgroup $G\subset\Gamma$, we require that on a neighborhood of $U^G$ the ambient metric $g_U$ agrees with the bundle metric on $NU^G$ given by $g_U$ via the exponential map along the normal directions.
        \item for each chart $(\Gamma, U, E)$ we require that the pullback connection $\nabla^E$ on $E$ satisfies the following condition. Given a subgroup $G\subset\Gamma$ and using the exponential map given by $g_U$, we identify a neighborhood of $U^G$ by a neighborhood of the zero-section of $NU^G$. Then, after also identifying $E|_{NU^G}$ with the pullback of $E|_{U^G}$ by $NU^G\rightarrow U^G$ and the parallel transport along normal geodesics using $\nabla^E$, the connection $\nabla^E$ agrees with pullback connection of the restriction $\nabla^E$ to $U^G$.
    \end{enumerate}
\end{definition}
\begin{lemma}[Lemma $3.15, 3.20$ \cite{bai2022integral}, Lemma $2.2$ \cite{bai2022arnold}]
    A compact derived orbifold chart $(\mathcal{U}, \E, s)$ admits a straightening in a neighborhood of $s^{-1}(0)$. 
\end{lemma}

\subsection{Universal Zero-Locus}
In \cite{fukaya2001floer} the authors showed us that in order to get an integer-valued Euler cycle from an orbibundle, one has to construct sections whose graphs are transverse to each isotropy strata of a certain algebraic variety.
Following \cite{fukaya2001floer}, \cite{bai2022integral} and \cite{bai2022arnold} we start the discussion in the linear case. \\

Let $\Gamma$ be a finite group and $V, W$ be two finite-dimensional complex $\Gamma$-representations. For a positive integer $d$, we denote by $\operatorname{Poly}^\Gamma_d(V, W)$ the finite-dimensional complex vectorspace of all $\Gamma$-equivariant complex polynomials from $V$ to $W$ of degree at most $d$. We have a natural evaluation map
\begin{align*}
    ev: V\times\operatorname{Poly}^\Gamma_d(V, W) \ & \to W \\
    (v, P) \ & \mapsto P(v)
\end{align*}
where we denote its zero-locus by $Z^\Gamma_d:=ev^{-1}(0)$.
\begin{proposition}[Theorem $3.3$ \cite{bai2022integral}, Theorem $2.5$ \cite{bai2022arnold}, Lemma $4.7$ \cite{parker2013integral}]\label{StratificationOfUni}
    There exists a unique Whitney stratification on $Z^\Gamma_d$ such that:
    \begin{enumerate}
    \item It is $\Gamma$-invariant.
    \item It is induced from a Whitney pre-stratification on $Z^\Gamma_d$ whose strata are algebraic varieties.
    \item It is invariant under any $\Gamma$-equivariant diffeomorphism of $V\times\operatorname{Poly}^\Gamma_d(V, W)$ preserving $Z^\Gamma_d$.
        \item If $G$ is a subgroup of $\Gamma$ and if we denote by $V_G\subseteq V$ the set of all $v\in V$ whose isotropy groups are exactly $G$, then for each $x\in Z^\Gamma_d\cap (V_{G}\times \operatorname{Poly}^\Gamma_d(V, W))$, the germ through $x$ is contained in $V_{G}\times \operatorname{Poly}^\Gamma_d(V, W)$. Moreover, such stratification is the minimal one, with respect to inclusion, satisfying these conditions. 
    \end{enumerate}
\end{proposition}
Now suppose that $B$ is a smooth manifold for which $\Gamma$ acts trivially on and let $\pi_V,\pi_W:V, W\rightarrow B$ be two $\Gamma$-equivariant complex vectorbundles. Then we have an induced complex vectorbundle $\operatorname{Poly}^\Gamma_d(V, W)\rightarrow B$ whose fibre over $b\in B$ is $\operatorname{Poly}^\Gamma_d(V_b, W_b)$. Similarly, we have a bundle map $ev:V\oplus \operatorname{Poly}^\Gamma_d(V, W)\rightarrow W$ covering the identity, where we denote its kernel by $\mathcal{Z}^\Gamma_d$. $\mathcal{Z}^\Gamma_d\rightarrow B$ is a $\Gamma$-equivariant complex vectorbundle whose fibre over $b\in B$ is $Z^\Gamma_d$.
\begin{proposition}\label{UniversalZeroLocus}Lemma $5$, SubLemma \cite{fukaya2001floer}
    For $d$ large enough, the vectorbundle $\mathcal{Z}^\Gamma_d\rightarrow B$ satisfies the following property:\\
    If $G$ is a subgroup of $\Gamma$ and if we denote by $V_G\subseteq V$ the set of all elements $v\in V$ whose isotropy group is exactly $G$, then $\mathcal{Z}^\Gamma_d\cap (V_G\times\operatorname{Poly}^\Gamma_d(V, W))$ is a smooth algebraic manifold. In particular, its of even real dimension.
\end{proposition}
\begin{proposition}\label{CanonicalWhitney Stratification} $\text{Theorem 3.3}$ \cite{bai2022integral}, $\text{Theorem 2.5}$ \cite{bai2022arnold}
The setting is as in Proposition \ref{UniversalZeroLocus}. There exists a unique Whitney stratification on $\mathcal{Z}^\Gamma_d$ such that for every $b\in B$, such stratification restricts to the one as in Proposition \ref{StratificationOfUni} on $Z^\Gamma_d$.
\end{proposition}
We call such Whitney stratification, the $\textbf{canonical Whitney Stratification}$ on $\mathcal{Z}^\Gamma_d$.
\subsection{FOP-sections}
Following the same discussion and notation as above, suppose that in addition $V$ is equipped with a $\Gamma$-invariant metric and let $\epsilon>0$. Consider $V_\epsilon\subset V$ an $\epsilon$-disk subbundle of $V$. We are interested in approximating smooth sections of $\pi_V^*W\rightarrow V_\epsilon$. To this end let $s:V_\epsilon\rightarrow\pi^*_V W$ be a $\Gamma$-equivariant section.
\begin{definition}\cite{bai2022integral}, \cite{bai2022arnold}, \cite{parker2013integral}\label{FOPsec}
    $s$ is said to be an $\textbf{FOP}$-section of degree at most $d$ if for every $(b, v)\in V_\epsilon$ we can find a $\Gamma$-invariant neighborhood $U\subseteq V_\epsilon$ of $(b, v)$ and a $\Gamma$-equivariant bundle map $f:U\rightarrow\operatorname{Poly}^\Gamma_d(V, W)$ such that
    \begin{equation*}
        s(b',v')=f(b',v')(v')
    \end{equation*}
    for every $(b',v')\in U$.
\end{definition}
We call the bundle map $f$ in Definition \ref{FOPsec}, a $\textbf{local lift}$ of $s$.
\begin{definition}\cite{bai2022integral}, \cite{bai2022arnold}
    Suppose that $s$ is an FOP-section of degree at most $d$. $s$ is said to be strongly-transverse at $(b,v)\in V_\epsilon$ if there exists $f:U\rightarrow\operatorname{Poly}^\Gamma_d(V, W)$ a local lift of $s$ around $(b, v)$ such that the graph map
  \begin{align*}
     \ U & \to V\times\operatorname{Poly}^\Gamma_d(V, W) \\
    (b', v') \ & \mapsto ((b', v'), f(b', v'))
\end{align*}
is transverse to each strata of the canonical Whitney stratification on $\{0\}\oplus\mathcal{Z}^\Gamma_d$.
\end{definition}
\begin{definition} \cite{bai2022integral}, \cite{bai2022arnold}
Given $(\mathcal{U}, \mathcal{E}, s)$ a straightened normally complex derived orbifold chart. We say:
\begin{enumerate}
    \item $s$ is an FOP section if for every bundle chart $(\Gamma, U, E)$ where $s:U\rightarrow E$ is a $\Gamma$-equivariant continuous section and after identifying a tubular neighborhood of $U^\Gamma$ with a disk subbundle $N_\epsilon U^\Gamma$ of $NU^\Gamma\rightarrow U^\Gamma$ and identifying
    \begin{equation*}
        E|_{N_{\epsilon}U^\Gamma}\cong\pi^*_{NU^\Gamma}\dot{E}^\Gamma\oplus\pi^*_{NU^\Gamma}\check{E}^\Gamma
    \end{equation*}
    as direct sum of trivial and non-trivial $\Gamma$-representations respectively, using the straightened structure, then we can write $s|_{N_\epsilon U^\Gamma}=(\dot{s}, \check{s})$. We require $\check{s}$ to be an FOP-section as in Definition \ref{FOPsec}.
    \item Suppose that $s$ is an FOP-section. We say $s$ is strongly transverse at $x\in\mathcal{U}$ if the following condition hold. If $(\Gamma, U, E)$ is an orbibundle chart centered at $x$ and after identifying $U$ with a disk subbundle $N_\epsilon U^\Gamma$ of the normal bundle $NU^\Gamma\rightarrow U^\Gamma$ and decomposing 
    \begin{equation*}
        E|_{N_{\epsilon}U^\Gamma}\cong\pi^*_{NU^\Gamma}\dot{E}^\Gamma\oplus\pi^*_{NU^\Gamma}\check{E}^\Gamma
    \end{equation*}
    as above and after writing $s=(\dot{s}, \check{s})$. We consider a bundle map $f:N_\epsilon U^\Gamma\rightarrow\operatorname{Poly}^\Gamma_{d}(NU^\Gamma, \check{E}^\Gamma)$ lifting $\check{s}$ and require the induced bundle map
    \begin{equation*}
        (\dot{s}, \operatorname{graph}(f)):N_\epsilon U^\Gamma\rightarrow\dot{E}^\Gamma\oplus(NU^\Gamma\oplus\operatorname{Poly}^G_d(NU^\Gamma, \check{E}^\Gamma))
    \end{equation*}
    to be transverse to the subbundle $\{0\}\times\mathcal{Z}^\Gamma_d(NU^\Gamma, \check{E}^\Gamma)$ with respect to the canonical Whitney stratification as in \ref{CanonicalWhitney Stratification}.
\end{enumerate}
\end{definition}
 \begin{proposition}\label{FOPstabilization}
    Suppose that $(\mathcal{U}, \mathcal{E}, s)$ is a normally complex derived orbifold chart equipped with a straightening and $s$ is a strongly-transverse FOP-section. Let $\pi_\F:\mathcal{F}\rightarrow\mathcal{U}$ be a complex orbibundle and $\tau_\F:\mathcal{F}\rightarrow\pi_\F^*\mathcal{F}$ be the tautological section. Then, the pullback metric on $\F$ together with the induced connection on $\pi_\F^*\E\rightarrow\F$ and the trivial connection on $\pi^*_\F\F\rightarrow\F$ gives a straightening on $\pi^*_\F\mathcal{E}\oplus\pi_\F^*\mathcal{F}\rightarrow\F$ for which the section $\pi_\F^*s\oplus\tau_\F:\mathcal{F}\rightarrow\pi^*_\F\mathcal{E}\oplus\pi_\F^*\mathcal{F}$ is a strongly transverse FOP-section.
\end{proposition}
\begin{proposition}\label{sumFOP}$\text{Corollary 2.4}$ \cite{bai2022arnold}.
    For $i=1,2$ let $(\mathcal{U}_i,\E_i,s_i)$ be a normally complex derived orbifold chart equipped with a straightening. Suppose that $s_i$ is a strongly-transverse FOP-section for $i=1,2$. Then, $s_1\oplus s_2:\mathcal{U}_1\times\mathcal{U}_2\rightarrow\pi_1^*\E_1\oplus\pi_2^*\E_2$ is a strongly-transverse FOP-section for the induced straightening on $\pi_1^*\E_1\oplus\pi_2^*\E_2\rightarrow\mathcal{U}_1\times\mathcal{U}_2$, where $\pi_i:\mathcal{U}_1\times\mathcal{U}_2\rightarrow\mathcal{U}_i$ is the projection map.
\end{proposition}
Now we state the results that allow us to approximate orbibundle sections by a strongly-transverse FOP-sections.
\begin{proposition}\label{AbsFOPapp}$\text{Proposition 1}$ \cite{fukaya2001floer},$\text{Propositin 6.4}$\cite{bai2022integral}
    Given a compact derived orbifold chart $(\mathcal{U}, \E, s)$ equipped with a straightening and $\epsilon>0$. Then there exists an orbibundle section $s_\epsilon:\mathcal{U}\rightarrow\E$ and a precompact neighborhood $O$ of $s^{-1}(0)$ such that:
    \begin{enumerate}
        \item $||s-s_\epsilon||_{C^0(O)}<\epsilon$.
        \item $s_\epsilon|_O:O\rightarrow\E|_O$ is a strongly-transverse FOP-section.
    \end{enumerate}
\end{proposition}
\begin{proposition}\label{relExtFOP}$\text{Proposition 1}$ \cite{fukaya2001floer}, $\text{Proposition 2.11}$ \cite{bai2022arnold} 
    Suppose that $(\mathcal{U},\E,s)$ is a compact derived orbifold chart equipped with a straightening, $O$ is a precompact neighborhood of $s^{-1}(0)$ and $\epsilon>0$. In addition, suppose that $K\subseteq\mathcal{U}$ is a compact set and $O'$ is a precompact neighborhood of $K$ such that $s|_{O'}:O'\rightarrow\E|_{O'}$ is a strongly transverse FOP-section on a neighborhood of $\overline{O}\cap K$. Then there exists an orbibundle section $s_\epsilon:\mathcal{U}\rightarrow\E$ such that:
    \begin{enumerate}
        \item $s_\epsilon$ is a strongly-transverse FOP-section on a neighborhood of $\overline{O}$.
        \item $s\equiv s_\epsilon$ on a neighborhood of $K$.
        \item $||s-s_\epsilon||_{C^0(O')}<\epsilon$.
    \end{enumerate}
\end{proposition}
Now we are ready to state the main theorem of this appendix which will allow us to extract integral Euler cycles in the setting of a normally complex orbibundle.
\begin{theorem}\label{IntegralCycles}$\text{Proposition 1}$ \cite{fukaya2001floer}, $\text{Lemma 4.20, Proposition 4.23}$ \cite{parker2013integral}, $\text{Theorem 1.1}$ \cite{bai2022integral}, $\text{Proposition 2.11}$ \cite{bai2022arnold} 
    Suppose that, $(\mathcal{U}, \E, s)$ is an oriented normally complex derived orbifold chart equipped with a straightening and let $\epsilon>0$. Denote by $\mathcal{U}_{free}\subseteq\mathcal{U}$ the set of all points with trivial isotropy and let $s_\epsilon:\mathcal{U}\rightarrow \E$ be as in Proposition \ref{AbsFOPapp}. Then,
    \begin{enumerate}
        \item $s_\epsilon^{-1}(0)\cap\mathcal{U}_{free}\hookrightarrow\mathcal{U}$ is an oriented pseudo-cycle of real dimension equal to $\dim_\R\mathcal{U}-\operatorname{rk}_\R\E$. 
        \item The homology class $[s_\epsilon^{-1}(0)\cap\mathcal{U}_{free}]\in H_*(\mathcal{U}; \Z)$ is an invariant of the normally complex derived orbifold chart $(\mathcal{U}, \E, s)$. In particular, it is independent of the choice of straightening on $\E\rightarrow\mathcal{U}$ and $s_\epsilon$.
    \end{enumerate}
\end{theorem}

 \nocite{Rabah2024} \nocite{Hirschi:2022twt} \nocite{mak2025orbifoldfloertheoryglobal}
\bibliographystyle{alpha}
\bibliography{ref}

\begin{thebibliography}{FOOO21b}

\bibitem[Abo10]{abouzaid2010geometric}
Mohammed Abouzaid.
\newblock A geometric criterion for generating the fukaya category.
\newblock {\em Publications Math{\'e}matiques de l'IH{\'E}S}, 112:191--240, 2010.

\bibitem[Abo11]{abouzaid2011topological}
Mohammed Abouzaid.
\newblock A topological model for the fukaya categories of plumbings.
\newblock {\em Journal of Differential Geometry}, 87(1):1--80, 2011.

\bibitem[ALR07]{adem2007orbifolds}
Alejandro Adem, Johann Leida, and Yongbin Ruan.
\newblock {\em Orbifolds and stringy topology}, volume 171.
\newblock Cambridge University Press, 2007.

\bibitem[AMS21]{abouzaid2021complex}
Mohammed Abouzaid, Mark McLean, and Ivan Smith.
\newblock {Complex Cobordism, Hamiltonian loops and Global Kuranishi Charts}.
\newblock {\em arXiv preprint arXiv:2110.14320}, 2021.

\bibitem[AMS24]{abouzaid2024gromov}
Mohammed Abouzaid, Mark McLean, and Ivan Smith.
\newblock {Gromov-Witten Invariants in Complex and Morava-Local K-Theories}.
\newblock {\em Geometric and Functional Analysis}, pages 1--87, 2024.

\bibitem[BV06]{boardman2006homotopy}
John~Michael Boardman and Rainer~M Vogt.
\newblock {\em Homotopy invariant algebraic structures on topological spaces}, volume 347.
\newblock Springer, 2006.

\bibitem[BX22a]{bai2022arnold}
Shaoyun Bai and Guangbo Xu.
\newblock {Arnold Conjecture over Integers}.
\newblock {\em arXiv preprint arXiv:2209.08599}, 2022.

\bibitem[BX22b]{bai2022integral}
Shaoyun Bai and Guangbo Xu.
\newblock An integral euler cycle in normally complex orbifolds and z-valued gromov-witten type invariants.
\newblock {\em arXiv preprint arXiv:2201.02688}, 2022.

\bibitem[CW15]{charest2015floer}
Fran{\c{c}}ois Charest and Chris~T Woodward.
\newblock {Floer Theory and Flips}.
\newblock {\em arXiv preprint arXiv:1508.01573}, 2015.

\bibitem[FO96]{fukaya1996arnold}
Kenji Fukaya and Kaoru Ono.
\newblock Arnold conjecture and gromov-witten invariant for general symplectic manifolds.
\newblock {\em The Arnoldfest (Toronto, ON, 1997)}, 24:173--190, 1996.

\bibitem[FO97]{fukaya1997zero}
Kenji Fukaya and Yong-Geun Oh.
\newblock {Zero-loop Open Strings in the cotangent bundle and Morse homotopy}.
\newblock {\em Asian Journal of Mathematics}, 1(1):96--180, 1997.

\bibitem[FO01]{fukaya2001floer}
Kenji Fukaya and Kaoru Ono.
\newblock Floer homology and gromov-witten invariant over integer of general symplectic manifolds--summary--.
\newblock In {\em Taniguchi Conference on Mathematics Nara'98}, volume~31, pages 75--92. Mathematical Society of Japan, 2001.

\bibitem[FOOO09]{Fukaya2009}
Kenji Fukaya, Yong-Geun Oh, Hiroshi Ohta, and Kaoru Ono.
\newblock {\em Lagrangian Intersection Floer Theory: Anomaly and Obstruction (Part 1)}, volume 46.1 of {\em AMS/IP Studies in Advanced Mathematics}.
\newblock American Mathematical Society / International Press, 2009.

\bibitem[FOOO10]{fukaya2010lagrangian}
Kenji Fukaya, Yong-Geun Oh, Hiroshi Ohta, and Kaoru Ono.
\newblock {\em Lagrangian intersection Floer theory: anomaly and obstruction, Part II}, volume~2.
\newblock American Mathematical Soc., 2010.

\bibitem[FOOO16]{fukaya2016exponential}
Kenji Fukaya, Yong-Geun Oh, Hiroshi Ohta, and Kaoru Ono.
\newblock Exponential decay estimates and smoothness of the moduli space of pseudoholomorphic curves.
\newblock {\em arXiv preprint arXiv:1603.07026}, 2016.

\bibitem[FOOO21a]{fukaya2021construction}
K~Fukaya, YG~Oh, H~Ohta, and K~Ono.
\newblock Construction of kuranishi structures on the moduli spaces of pseudo holomorphic disks: I, to appear in surveys in differential geometry.
\newblock {\em arXiv preprint arXiv:1710.01459}, 2021.

\bibitem[FOOO21b]{fukaya1710construction}
K~Fukaya, YG~Oh, H~Ohta, and K~Ono.
\newblock {Construction of Kuranishi structures on the moduli spaces of pseudo holomorphic disks: I, to appear in Surveys in Differential Geometry}.
\newblock {\em arXiv preprint arXiv:1710.01459}, 2021.

\bibitem[FOOO24]{fukaya2024corrigendum}
Kenji Fukaya, Yong-Geun Oh, Hiroshi Ohta, and Kaoru Ono.
\newblock Corrigendum of" construction of kuranishi structures on the moduli spaces of pseudo holomorphic disks i, surveys in differential geometry xxii (2018), 133-190".
\newblock {\em arXiv preprint arXiv:2403.19683}, 2024.

\bibitem[Giv00]{givental2000wdvv}
Alexander Givental.
\newblock On the wdvv equation in quantum k-theory.
\newblock {\em Michigan Mathematical Journal}, 48(1):295--304, 2000.

\bibitem[Gro85]{gromov1985pseudo}
Mikhael Gromov.
\newblock Pseudo holomorphic curves in symplectic manifolds.
\newblock {\em Inventiones mathematicae}, 82(2):307--347, 1985.

\bibitem[Hir23]{hirschi2023properties}
Amanda Hirschi.
\newblock Properties of gromov-witten invariants defined via global kuranishi charts.
\newblock {\em arXiv preprint arXiv:2312.03625}, 2023.

\bibitem[HS24]{Hirschi:2022twt}
Amanda Hirschi and Mohan Swaminathan.
\newblock {Global Kuranishi charts and a product formula in symplectic Gromov{\textendash}Witten theory}.
\newblock {\em Selecta Math.}, 30(5):87, 2024.

\bibitem[KL06]{katz2006enumerative}
Sheldon Katz and Chiu-Chu~Melissa Liu.
\newblock Enumerative geometry of stable maps with lagrangian boundary conditions and multiple covers of the disc.
\newblock {\em Geometry \& Topology Monographs}, 8:1--47, 2006.

\bibitem[KO00]{kwon2000structure}
Daesung Kwon and Yong-Geun Oh.
\newblock {Structure of the Image of (pseudo)-holomorphic Discs with Totally Real Boundary Condition}.
\newblock {\em Communications in Analysis and Geometry}, 8(1):31--82, 2000.

\bibitem[Las06]{lashof2006stable}
Richard Lashof.
\newblock Stable g-smoothing.
\newblock In {\em Algebraic Topology Waterloo 1978: Proceedings of a Conference Sponsored by the Canadian Mathematical Society, NSERC (Canada), and the University of Waterloo, June 1978}, pages 283--306. Springer, 2006.

\bibitem[Laz11]{lazzarini2011relative}
Laurent Lazzarini.
\newblock {Relative Frames on J-holomorphic Curves}.
\newblock {\em Journal of Fixed Point Theory and Applications}, 9:213--256, 2011.

\bibitem[LLY97]{lian1997mirror}
Bong Lian, Kefeng Liu, and Shing-Tung Yau.
\newblock Mirror principle i.
\newblock {\em arXiv preprint alg-geom/9712011}, 1997.

\bibitem[Mes16]{mescher2016perturbed}
Stephan Mescher.
\newblock Perturbed gradient flow trees and {$A_\infty$}-algebra structures on morse cochain complexes.
\newblock {\em arXiv preprint arXiv:1611.07916}, 2016.

\bibitem[Mil64]{milnor1964microbundles}
John Milnor.
\newblock Microbundles: Part i.
\newblock {\em Topology}, 3:53--80, 1964.

\bibitem[MS12]{mcduff2012j}
Dusa McDuff and Dietmar Salamon.
\newblock {\em J-holomorphic curves and symplectic topology}, volume~52.
\newblock American Mathematical Soc., 2012.

\bibitem[MSS25]{mak2025orbifoldfloertheoryglobal}
Cheuk~Yu Mak, Sobhan Seyfaddini, and Ivan Smith.
\newblock Orbifold floer theory for global quotients and hamiltonian dynamics, 2025.

\bibitem[MW10]{ma2010geometric}
Sikimeti Ma’u and Chris Woodward.
\newblock Geometric realizations of the multiplihedra.
\newblock {\em Compositio Mathematica}, 146(4):1002--1028, 2010.

\bibitem[Par13]{parker2013integral}
Brett Parker.
\newblock Integral counts of pseudo-holomorphic curves.
\newblock {\em arXiv preprint arXiv:1309.0585}, 2013.

\bibitem[Par16]{pardon2016algebraic}
John Pardon.
\newblock An algebraic approach to virtual fundamental cycles on moduli spaces of pseudo-holomorphic curves.
\newblock {\em Geometry \& Topology}, 20(2):779--1034, 2016.

\bibitem[Par22]{pardon2022enough}
John Pardon.
\newblock Enough vector bundles on orbispaces.
\newblock {\em Compositio Mathematica}, 158(11):2046--2081, 2022.

\bibitem[Rab]{Rabah2024}
Mohamad Rabah.
\newblock {Lagrangian Floer Theory over Z}.
\newblock Manuscript in preparation.

\bibitem[Rez22]{rezchikov2022integralarnoldconjecture}
Semon Rezchikov.
\newblock Integral arnol'd conjecture, 2022.

\bibitem[Sie99]{siebert1999gromov}
N~Siebert.
\newblock Gromov-witten invariants for general symplectic manifolds, new trends in algebraic geometry (warwick, 1996), 375-424, cambrdge uni, 1999.

\bibitem[Sta63]{stasheff1963homotopy}
James~Dillon Stasheff.
\newblock Homotopy associativity of h-spaces. ii.
\newblock {\em Transactions of the American Mathematical Society}, 108(2):293--312, 1963.

\bibitem[Whi34]{whitney1934analytic}
Hassler Whitney.
\newblock Analytic extensions of differentiable functions defined in closed sets.
\newblock {\em Transactions of the American Mathematical Society}, 36(1):63--89, 1934.

\bibitem[Zer17]{zernik2017moduli}
Amitai~Netser Zernik.
\newblock {Moduli of Open Stable Maps to a Homogeneous Space}.
\newblock {\em arXiv preprint arXiv:1709.07402}, 2017.

\end{thebibliography}
\end{document}